\documentclass[11pt]{mcom-l}
\usepackage{amssymb,amsfonts,amsmath,epsfig,hyperref}
\usepackage[margin=1.57in]{geometry}
\usepackage{mathrsfs}
\usepackage{subfigure}
\renewcommand{\theequation}{\thesection.\arabic{equation}}
\allowdisplaybreaks
\newcommand{\nn}{\nonumber}
\def\refe#1{(\ref{#1})}

\def\d{\,{\rm d}}

\title{ 
\bf A new approach for numerical simulation of 
the time-dependent 
Ginzburg--Landau equations
} 

\author{Buyang Li$^\dagger\,\,\,$ \and\,\,  Zhimin Zhang$^\ddagger$}

\thanks{This work was supported in part by the National Natural Science Foundation of China
(NSFC) under grants No. 11301262, No. 11471031, No. 91430216,
and the US National Science Foundation (NSF) through grants DMS-1115530 and DMS-1419040. \newline
\indent $^\dagger$Department of Mathematics,
Nanjing University, Nanjing, China. 
(buyangli@nju.edu.cn) \newline
\indent $^\ddagger$Beijing Computational Science Research Center,
Beijing, China.\newline
\indent $^\ddagger$Department of Mathematics, Wayne State University, Detroit, USA
(ag7761@wayne.edu).
}

\thanks{}

\begin{document}

\maketitle

\begin{abstract}
\small  
We introduce a new approach for
finite element simulations of the time-dependent 
Ginzburg--Landau equations (TDGL) in a general
curved polygon, possibly with reentrant corners.
Specifically, we reformulate the TDGL
into an equivalent system of equations 
by decomposing the magnetic potential 
to the sum of its divergence-free 
and curl-free parts, respectively. 
Numerical simulations of vortex dynamics show that,
in a domain with reentrant corners, 
the new approach is much more stable and accurate 
than the old approaches of solving the TDGL directly 
(under either the temporal gauge or the Lorentz gauge);
in a convex domain, 
the new approach gives comparably accurate 
solutions as the old approaches.

\end{abstract}



\pagestyle{myheadings}
\thispagestyle{plain}
\markboth{}{}

\section{Introduction}
\setcounter{equation}{0}

Based on the Ginzburg--Landau theory of superconductivity \cite{GL},
the macroscopic state of a superconductor  
is described by the complex-valued order parameter $\psi$,
the real scalar-valued electric potential $\phi$,
and the real vector-valued magnetic potential ${\bf A}$.
In the nondimensionalization form, 
the order parameter satisfies that $0\leq |\psi|^2\leq 1$, 
where $|\psi|^2=0$ corresponds to 
the normal state and $|\psi|^2=1$ corresponds to 
the superconducting state, 
and $0<|\psi|^2<1$ represents an intermediate state 
between the normal and superconducting states.
If the superconductor occupies a long cylinder in the $x_3$-direction 
with a finite cross section and the external magnetic field 
is ${\bf H}=(0,0,H)$, then 
the order parameter $\psi$ and the magnetic potential
${\bf A}=(A_1,A_2)$
are governed by the 
time-dependent Ginzburg--Landau equations (TDGL) 
\begin{align}
&\eta\frac{\partial \psi}{\partial t}
+ \bigg(\frac{i}{\kappa} \nabla + \mathbf{A}\bigg)^{2} \psi
 + (|\psi|^{2}-1) \psi + i\eta\kappa\psi\phi = 0,
\label{GLLPDEq1}\\[5pt]
&\frac{\partial \mathbf{A}}{\partial t} 
+ \nabla\times(\nabla\times{\bf A})
+\nabla \phi+  {\rm Re}\bigg[\psi^*\bigg(\frac{i}{\kappa} \nabla 
+ \mathbf{A}\bigg) \psi\bigg] =  \nabla\times H,
\label{GLLPDEq2}
\end{align}
in the two-dimensional cross sectional domain $\Omega$,
where $\eta$ is the normalized conductivity,  
$\kappa$ is the Ginzburg-Landau parameter, 
and $\psi^*$ denotes the complex conjugate of $\psi$.
Discovered by Schmid \cite{Schmid} and
derived by Gor'kov and Eliashberg \cite{GE} 
from the microscopic principles,
the TDGL was widely accepted for 
 simulation of transient behaviors and
vortex motions of superconductors 
\cite{FUD91,LMG91}. 
Variables of physical interest in this model are
the superconducting density $|\psi|^2$, 
the magnetic induction field $B=\nabla\times{\bf A}$, and
the electric field ${\bf E}=\partial_t{\bf A}+\nabla\phi$. 
The boundary conditions are 
\begin{align}
&\biggl(\frac{i}{\kappa}\nabla\psi +
\mathbf{A}\psi\biggl) \cdot{\bf n} = 0
&\mbox{on}\,\,\,\,\partial\Omega , \\
&B= H  
&\mbox{on}\,\,\,\,\partial\Omega ,\\
&{\bf E}\cdot{\bf n}=0 
&\mbox{on}\,\,\,\,\partial\Omega ,
\label{OrBDC3}
\end{align}
where $\mathbf{n}$ denotes the unit outward normal vector on the boundary
$\partial\Omega$.
Detailed description of the physics of superconductivity phenomena 
can be found in the review articles \cite{CHO,DGP92} and the books
\cite{Gennes,Tinkham}. 
Here, in a two-dimensional domain, we use the notations
\begin{align*}
&\nabla\times {\bf A}
=\frac{\partial A_2}{\partial x_1}-\frac{\partial A_1}{\partial x_2},
\qquad 
\nabla\cdot {\bf A}=\frac{\partial A_1}{\partial x_1}
+\frac{\partial A_2}{\partial x_2},\\
&\nabla\times H=\bigg(\frac{\partial H}{\partial x_2},\,
-\frac{\partial H}{\partial x_1}\bigg),\quad 
\nabla\psi=\bigg(\frac{\partial \psi}{\partial x_1},\,
\frac{\partial \psi}{\partial x_2}\bigg).
\end{align*}

The TDGL requires an additional gauge condition 
to determine the solution uniquely
 \cite{Abrikosov88,CHL}.
Via a gauge transformation
\begin{align}
&\psi=\widetilde\psi e^{i\kappa\chi},
\qquad
{\bf A}=\widetilde {\bf A}+\nabla\chi,
\qquad
\phi=\widetilde\phi-\frac{\partial \chi}{\partial t} ,
\end{align}
the two solutions $(\psi,{\bf A},\phi)$ and 
$(\widetilde\psi,\widetilde{\bf A},\widetilde\phi)$ 
are equivalent in producing the physical variables, e.g.
superconducting density, magnetic induction and
electric field. 
As a consequence, solving the TDGL under different gauges
is theoretically equivalent in calculating
the quantities of physical interest.
However,
solving the TDGL under different gauges 
is not equivalent computationally. 
It is important to use a gauge under which 
the numerical solution is stable and accurate. 

A widely used gauge in numerical simulations 
is the temporal gauge $\phi=0$; see \cite{FUD91,LMG91}.
Numerical simulations of the TDGL 
under the temporal gauge
have been done in many works
with either finite element or
finite difference methods;
see \cite{ASPM11,GKLLP96,MH98,
RPCA04,VMB03,WA02}. 
Under the temporal gauge,
\refe{GLLPDEq1}-\refe{GLLPDEq2} reduce to 
\begin{align}
&\eta\frac{\partial \psi}{\partial t} + \biggl(\frac{i}{\kappa}\nabla +
\mathbf{A}\biggl)^{2} \psi + (|\psi|^{2}-1) \psi = 0,
\label{TPDEq1}\\[5pt]
&\frac{\partial \mathbf{A}}{\partial t} +
\nabla\times(\nabla\times{\bf A}) 
+ {\rm Re}\bigg[\psi^*\bigg(\frac{i}{\kappa} \nabla 
+ \mathbf{A}\bigg) \psi\bigg]  = \nabla\times  H ,
\label{TPDEq2}
\end{align}
and the boundary conditions reduce to 
\begin{align}
&\biggl(\frac{i}{\kappa}\nabla\psi +
\mathbf{A}\psi\biggl) \cdot{\bf n} = 0 
&\mbox{on}\,\,\,\, \partial\Omega , \\
&\nabla\times{\bf A}= H 
&\mbox{on}\,\,\,\, \partial\Omega , \\
&{\bf A}\cdot{\bf n}=0
&\mbox{on}\,\,\,\, \partial\Omega . 
\label{TBDC}
\end{align}
Since $\int_\Omega|\nabla\times{\bf A}|^2\d x$ is not equivalent to 
$\int_\Omega|\nabla{\bf A}|^2\d x$,
the equation \refe{TPDEq2} is degenerate parabolic.
Due to the degeneracy and the nonlinear structure,
both theoretical analysis and numerical approximation of 
\refe{TPDEq1}-\refe{TBDC} are difficult. 
In a smooth domain $\Omega$, 
existence and uniqueness of a solution for this system
were proved in \cite{Du94}.
Finite element approximations of 
\refe{TPDEq1}-\refe{TBDC} and convergence of the numerical solutions 
have been reviewed in \cite{Du05} and an alternating Crank--Nicolson 
schemes was proposed in \cite{MH98}.
Some implicit, explicit and implicit-explicit 
time discretization schemes were studied in \cite{GKL02}.  
For the finite element approximations, 
error estimates were carried out for a regularized
problem, by adding a term $-\epsilon\nabla(\nabla\cdot{\bf A})$
to the equation \refe{TPDEq2}. 
Depending on the parameter $\epsilon$, 
convergence rate of the numerical solution to the
exact solution cannot be expressed explicitly.
Although an explicit convergence rate was proved
in \cite{Yang}, the strong regularity assumption on the solution 
restrict the problem to a smooth domain 
without corners.
In a domain with reentrant corners,
well-posedness of \refe{TPDEq1}-\refe{TBDC} 
remains open and convergence of 
the numerical solution is not known yet.

To overcome the difficulties caused by degeneracy, 
the Lorentz gauge 
$\phi=-\nabla\cdot{\bf A}$ was introduced in \cite{CH94}
for the simulation of TDGL.
Under the Lorentz gauge,
\refe{GLLPDEq1}-\refe{GLLPDEq2} 
reduce to 
\begin{align}
&\eta\frac{\partial \psi}{\partial t} 
+ \bigg(\frac{i}{\kappa} \nabla + \mathbf{A}\bigg)^{2} \psi
 + (|\psi|^{2}-1) \psi 
 -i\eta \kappa \psi \nabla\cdot{\bf A} = 0,
\label{PDE1}\\[5pt]
&\frac{\partial \mathbf{A}}{\partial t} 
+ \nabla\times(\nabla\times{\bf A})
-\nabla(\nabla\cdot{\bf A}) 
+  {\rm Re}\bigg[\psi^*\bigg(\frac{i}{\kappa} \nabla 
+ \mathbf{A}\bigg) \psi\bigg] =  \nabla\times H, 
\label{PDE2}
\end{align}
with the boundary conditions 
\begin{align}
&\frac{i}{\kappa}\nabla\psi  \cdot{\bf n} = 0 
&\mbox{on}\,\,\,\, \partial\Omega , \\
&\nabla\times{\bf A}= H
&\mbox{on}\,\,\,\, \partial\Omega , \\
&{\bf A}\cdot{\bf n}=0
&\mbox{on}\,\,\,\, \partial\Omega . 
\label{TBDC2}
\end{align}
The equation \refe{PDE2} is parabolic
without degeneracy, as
$\|\nabla\times{\bf A}\|_{L^2}^2+\|\nabla\cdot{\bf A}\|_{L^2}^2$
is equivalent to $\|\nabla {\bf A}\|_{L^2}^2$ for any 
${\bf A}\in {\bf H}^1_{\rm n}(\Omega):=\{{\bf a}\in H^1(\Omega)^2: 
{\bf a}\cdot{\bf n}=0~\,\mbox{on}~\,\partial\Omega\}$.
In a bounded smooth domain,
existence and uniqueness of solution for 
\refe{PDE1}-\refe{TBDC2}
were proved by Chen et al. \cite{CHL}.
Error estimates of the FEM were presented in 
\cite{CD01} with a backward 
Euler scheme and presented in \cite{GLS} 
with a linearized Crank--Nicolson scheme. 
Besides, the regularized TDGL under temporal gauge are 
approximately in the form of \refe{PDE1}-\refe{PDE2}; see
\cite{Mu97}.
If the domain contains a reentrant corner, then 
the magnetic potential may not be in $L^2(0,T;{\bf H}^1_{\rm n}(\Omega))$
and well-posedness of the TDGL remains open
in this case.

Overall, convergence of the numerical solution is not guaranteed 
under either gauge 
if the domain contains reentrant corners.
Meanwhile, correct numerical approximation of the TDGL 
in domains with reentrant corners are important for physicists 
to study the effects of surface defects in superconductivity  
\cite{ASPM11,BKP05,VMB03}, which was often 
done by solving \refe{TPDEq1}-\refe{TBDC} 
or \refe{PDE1}-\refe{TBDC2}  
with the finite element method (FEM).
We believe that the magnetic potential ${\bf A}$
may not be in $L^2(0,T;{\bf H}^1_{\rm n}(\Omega))$ in a domain
with reentrant corners, and the finite element solutions
of \refe{PDE1}-\refe{TBDC2} may converge to 
an incorrect solution.
Moreover, the incorrect numerical solution of 
${\bf A}$ may pollutes the numerical solution
of $\psi$ through the coupling of the equations 
and lead to wrong approximation of the physical
quantity $|\psi|$.

In this paper, we introduce a new approach
to simulate the TDGL in a curved polygon which
may contain reentrant corners.
Specifically, we reformulate 
the TDGL into an equivalent system of equations 
whose solutions are in $L^2(0,T;H^1(\Omega))$, and propose
a simple numerical scheme to solve the reformulated system.
We shall demonstrate the efficiency of the new approach
via numerical simulations, comparing the numerical results
with the numerical solutions of \refe{TPDEq1}-\refe{TBDC} 
and \refe{PDE1}-\refe{TBDC2} 
by using the same triangulation and finite element space.
We will see that, in a domain with reentrant corners, 
the numerical solution of 
\refe{TPDEq1}-\refe{TBDC} is unstable and the 
numerical solution of \refe{PDE1}-\refe{TBDC2} 
is incorrect, while our new approach leads
to stable and accurate numerical solutions.
Existence and uniqueness of 
solutions for the reformulated system
and its equivalence to the original TDGL system 
are proved in a separate paper \cite{LZ2}.

\section{A new approach}
\setcounter{equation}{0}
It is well known that any vector field is a sum of  
a divergence-free vector field and 
a curl-free vector field \cite{BCNS}.
If we assume that 
${\bf A}\in L^2(\Omega)\times L^2(\Omega)$ and 
$\nabla\cdot{\bf A}\in L^2(\Omega)\times L^2(\Omega)$,
then the magnetic potential has the decomposition
\begin{align}\label{DecpA}
&{\bf A} = \nabla\times u + \nabla v,
\end{align}
where $u$ and $v$ are the solutions of 
\begin{align*}
\left\{\begin{array}{ll}
-\Delta u=\nabla\times{\bf A} &\mbox{in}\,\,\,\,\Omega,\\
u=0&\mbox{on}\,\,\,\,\partial\Omega,
\end{array}
\right.
\end{align*}
and
\begin{align*}
\left\{\begin{array}{ll}
\Delta v=\nabla\cdot{\bf A} &\mbox{in}\,\,\,\,\Omega,\\
\partial_nv=0 &\mbox{on}\,\,\,\,\partial\Omega,
\end{array}
\right.
\end{align*}
respectively.
This decomposition is consistent with the boundary condition
${\bf A}\cdot{\bf n}=0$, which is a consequence of
$u=0$ and $\partial_nv=0$ on $\partial\Omega$.
Similarly, we have the decomposition
\begin{align}\label{DecpRe}
&{\rm Re}\bigg[\psi^*\bigg(\frac{i}{\kappa} \nabla 
+ \mathbf{A}\bigg) \psi\bigg]=\nabla\times p + \nabla q ,
\end{align}
where $p$ and $q$ are the solutions of 
\begin{align*}
\left\{\begin{array}{ll}
\displaystyle 
\Delta p =-\nabla\times \bigg({\rm Re}\bigg[
\psi^*\bigg(\frac{i}{\kappa} \nabla 
+ \mathbf{A}\bigg) \psi\bigg]\bigg) &\mbox{in}\,\,\,\,\Omega,\\
p =0 &\mbox{on}\,\,\,\,\partial\Omega,
\end{array}\right.
\end{align*}
and
\begin{align*}
\left\{\begin{array}{ll}
\displaystyle 
\Delta q =\nabla\cdot \bigg({\rm Re}\bigg[
\psi^*\bigg(\frac{i}{\kappa} \nabla 
+ \mathbf{A}\bigg) \psi\bigg]\bigg) 
&\mbox{in}\,\,\,\,\Omega,\\
\partial_n q =0 &\mbox{on}\,\,\,\,\partial\Omega,
\end{array}\right.
\end{align*}
respectively.
With \refe{DecpA} and \refe{DecpRe},
the equation \refe{PDE2} reduces to
$$
\nabla\times\bigg(
\frac{\partial u}{\partial t} -\Delta u-H+p\bigg)
+\nabla \bigg(\frac{\partial v}{\partial t} -\Delta v
+q \bigg)
=  0 .
$$
In the above equation, the divergence-free and curl-free parts 
must vanish simultaneously.
Thus we can reformulate \refe{PDE1}-\refe{PDE2} as
\begin{align}
&\eta\frac{\partial \psi}{\partial t} 
+ \bigg(\frac{i}{\kappa} \nabla + {\bf A}\bigg)^{2} \psi
 + (|\psi|^{2}-1) \psi-i\eta \kappa \psi \nabla\cdot{\bf A} = 0,
\label{RFPDE1}\\[5pt]
&\Delta p =-\nabla\times \bigg({\rm Re}\bigg[
\psi^*\bigg(\frac{i}{\kappa} \nabla 
+ \mathbf{A}\bigg) \psi\bigg]\bigg)\label{RFPDEp}\\[5pt]
&\Delta q =\nabla\cdot \bigg({\rm Re}\bigg[
\psi^*\bigg(\frac{i}{\kappa} \nabla 
+ \mathbf{A}\bigg) \psi\bigg]\bigg)\label{RFPDEq}\\[5pt]
&\frac{\partial u}{\partial t} -\Delta u
=  H-p ,\label{RFPDEu}\\[5pt]
&\frac{\partial v}{\partial t} -\Delta v
=  -q ,
\label{RFPDEv}
\end{align}
with the boundary conditions  
\begin{align}
&\partial_n\psi =0 &\mbox{on}\,\,\,\,\partial\Omega,\\
&p=0 &\mbox{on}\,\,\,\,\partial\Omega,\\
&\partial_nq =0 &\mbox{on}\,\,\,\,\partial\Omega,\\
&u=0 &\mbox{on}\,\,\,\,\partial\Omega,\\
&\partial_nv =0 &\mbox{on}\,\,\,\,\partial\Omega,
\label{RFbc}
\end{align}
and the initial conditions 
\begin{align}
& \psi(\cdot,0) = \psi_{0} , \quad u(\cdot,0) =
u_{0} , \quad v(\cdot,0) =
v_{0} ,\quad
 \mathrm{in}\,\,\, \Omega \, ,
\label{RFinit}
\end{align}
where $u_0$ and $v_0$ are defined by
$$
\left\{\begin{array}{ll}
-\Delta u_0=\nabla\times{\bf A}_0 &\mbox{in}~~\Omega,\\
u_0=0 &\mbox{on}~~\partial\Omega,
\end{array}\right.
\qquad\mbox{and}\qquad
\left\{\begin{array}{ll}
\Delta v_0=\nabla\cdot{\bf A}_0 &\mbox{in}~~\Omega,\\
\partial_nv_0=0 &\mbox{on}~~\partial\Omega,
\end{array}\right.
$$
with $\int_\Omega v_0(x)\d x=0$. 
From this new system of equations,
one can solve the order parameter $\psi$ and
find the magnetic potential ${\bf A}=\nabla\times u+\nabla v$
by solving $u$ and $v$.

Unlike the system \refe{PDE1}-\refe{TBDC2} whose solution
${\bf A}$ may not be in $L^2(0,T;{\bf H}_n^1(\Omega))$,
the reformulated system \refe{RFPDE1}-\refe{RFinit} consists of 
heat equations and Poisson's equations
whose solutions are in $L^2(0,T;H^1(\Omega))$ in a 
arbitrary curved polygon. 
Thus the new system is easier to solve 
than the original system of equations.
Here we propose a simple linearized and decoupled FEM
to solve \refe{RFPDE1}-\refe{RFinit} based on the
backward Euler time-stepping scheme.

For a given triangulation of the domain $\Omega$, 
we let $\mathcal{V}_{h}^{r}$ be the space 
of complex-valued globally continuous piecewise polynomials of degree $r\geq 1$ 
subject to the triangulation, 
let $V_{h}^{r}$ be the subspace of $\mathcal{V}_{h}^{r}$ consisting
of real-valued functions, 
and let $ \mathring V_{h}^{r}$ be the subspace of $V_h^r$
consisting of functions which are zero on the boundary. 
Let $0=t_0<t_1<\cdots<t_N=T$ be a uniform partition of the time interval with
define $\tau = T/N$.
For the given $\psi_h^n$, ${\bf A}_h^n$, $u_h^n$, $v_h^n$, 
we first calculate $\psi^{n+1}_h\in {\mathcal  V}_{h}^{r}$, 
by solving the equation
\begin{align}
&\bigg(\frac{\psi^{n+1}_h-\psi^{n}_h}{\tau}, \varphi\bigg) +\bigg
(\bigg(\frac{i}{\kappa}\nabla\psi^{n+1}_h + \mathbf{A}^{n}_h\psi^{n+1}_h\bigg)
,\bigg(\frac{i}{\kappa}\nabla \varphi + \mathbf{
A}^{n}_h \varphi\bigg) \bigg) \nn\\
&\qquad\qquad\qquad +\big((|\psi^{n}_h|^{2}-1) \psi^{n+1}_h, \varphi\big)
+\big(i\eta\kappa {\bf A}_h^n,\nabla((\psi_h^{n+1})^* \varphi)\big)= 0,
&\forall\,\varphi\in{\mathcal  V}^r_h,
\label{FEMEq1}
\end{align}
and define 
$$
F_h^{n+1}= {\rm Re}\bigg[(\psi^n_h)^*\bigg(\frac{i}{\kappa}\nabla\psi^{n+1}_h 
+ \mathbf{A}^n_h \psi^{n+1}_h\bigg)\bigg] .
$$
Then we look for $p^{n+1}_h,u^{n+1}_h\in \mathring V_h^{r}$ and
$q^{n+1}_h,v^{n+1}_h\in V_h^{r}$ satisfying the equations
\begin{align}
&(\nabla p^{n+1}_h,\nabla \xi)=\big(F_h^{n+1},
\nabla\times \xi\big),
&\forall\,\xi \in \mathring V^r_h,
\label{FEMEq2}\\
&(\nabla q^{n+1}_h,\nabla \zeta)
=\big(F_h^{n+1},\nabla \zeta\big)
&\forall\,   \zeta\in V^r_h,
\label{FEMEq3}\\
&\bigg(\frac{u^{n+1}_h-u_h^n}{\tau}, \theta\bigg) 
+\big(\nabla u^{n+1}_h,\nabla\theta\big)=(f^{n+1}-p^{n+1}_h,\theta),
&\forall\, \theta\in \mathring V^r_h,
\label{FEMEq4}\\
&\bigg(\frac{v^{n+1}_h-v_h^n}{\tau}, \vartheta\bigg) 
+\big(\nabla v^{n+1}_h,\nabla\vartheta\big)
=(-q^{n+1}_h,\vartheta),
&\forall\,\vartheta\in V^r_h,
\label{FEMEq5}
\end{align} 
and set 
${\bf A}_h^{n+1}=\nabla\times u_h^{n+1}+\nabla v_h^{n+1}$. 
At the initial time step, $u^0_h\in\mathring V^r_h$ 
and $v^0_h\in V^r_h$ can be solved from
\begin{align}
&(\nabla u^{0}_h,\nabla \xi)=\big({\bf A}_0,
\nabla\times \xi\big)
,\quad\forall~\xi\in\mathring V_h^r ,
\label{FEMEqu0}\\
&(\nabla v^{0}_h,\nabla \zeta)
=\big({\bf A}_0,\nabla\cdot \zeta\big)
,\quad~\, \forall~\zeta\in V_h^r ,
\label{FEMEqv0}
\end{align}
and $\psi^0_h$ can be chosen as 
the Lagrange interpolation of $\psi_0$. 
At each time step, one only needs to solve a system of linear
equations.

In the next section, we demonstrate the efficiency 
of the proposed scheme via numerical simulations, by comparing
the results with the numerical solutions of
\refe{TPDEq1}-\refe{TBDC} and 
\refe{PDE1}-\refe{TBDC2}.
\medskip

{\it Remark 2.1}$\quad$
If the magnetic induction $B=\nabla \times {\bf A}$
and electric field ${\bf E}=\partial_t{\bf A}-\nabla(\nabla\cdot{\bf A})$
are also desired, one can solve
\begin{align}
&\frac{\partial w}{\partial t} 
-\Delta w+  \nabla\times {\rm Re}\bigg[\psi^*\bigg(\frac{i}{\kappa} \nabla 
+ \mathbf{A}\bigg) \psi\bigg] = -\frac{\partial H}{\partial t}   \label{RFPDEw} 
\end{align}
additionally, with $w=\nabla\times {\bf A}-H$, which 
is derived by considering the
curl of \refe{PDE2}. Then one has
\begin{align}
&B=w+H,\\
&{\bf E} = -\nabla\times w
-{\rm Re}\bigg[\psi^*\bigg(\frac{i}{\kappa} \nabla 
+ \mathbf{A}\bigg) \psi\bigg] .
\end{align}
A fully discrete scheme for solving \refe{RFPDEw} is given by
\begin{align}
&\bigg(\frac{w^{n+1}_h-w_h^n}{\tau}, \chi\bigg) 
+\big(\nabla w^{n+1}_h,\nabla\chi\big)
=-(F^{n+1}_h,\nabla\times\chi)-(\partial_tH^{n+1},\nabla\times\chi),
&\forall\,\chi\in \mathring V^r_h,
\end{align} 
which can be solved with \refe{FEMEq1}-\refe{FEMEq5} together.
Then the magnetic induction and electric field can be 
approximated by
\begin{align}
&B_h^{n+1}=w_h^{n+1}+H^{n+1},\\
&{\bf E}_h^{n+1} = -\nabla\times w_h^{n+1}
-F_h^{n+1}.
\end{align}

In this paper, we focus on numerical simulation of the 
superconductivity density $|\psi|^2$. \medskip

\section{Numerical simulations}\label{femmethod}
\setcounter{equation}{0}
$\quad$

In this section, we present numerical simulations of the
vortex dynamics in domains with or without reentrant corners,
and compare the numerical solutions given by the different approaches
by using the same triangulation and time-step size,
with the backward Euler scheme for time discretization.\bigskip

\begin{figure}[htp]
\centering
\subfigure[$h=1/16$] 
{\includegraphics[height=1.15in,width=1.6in]{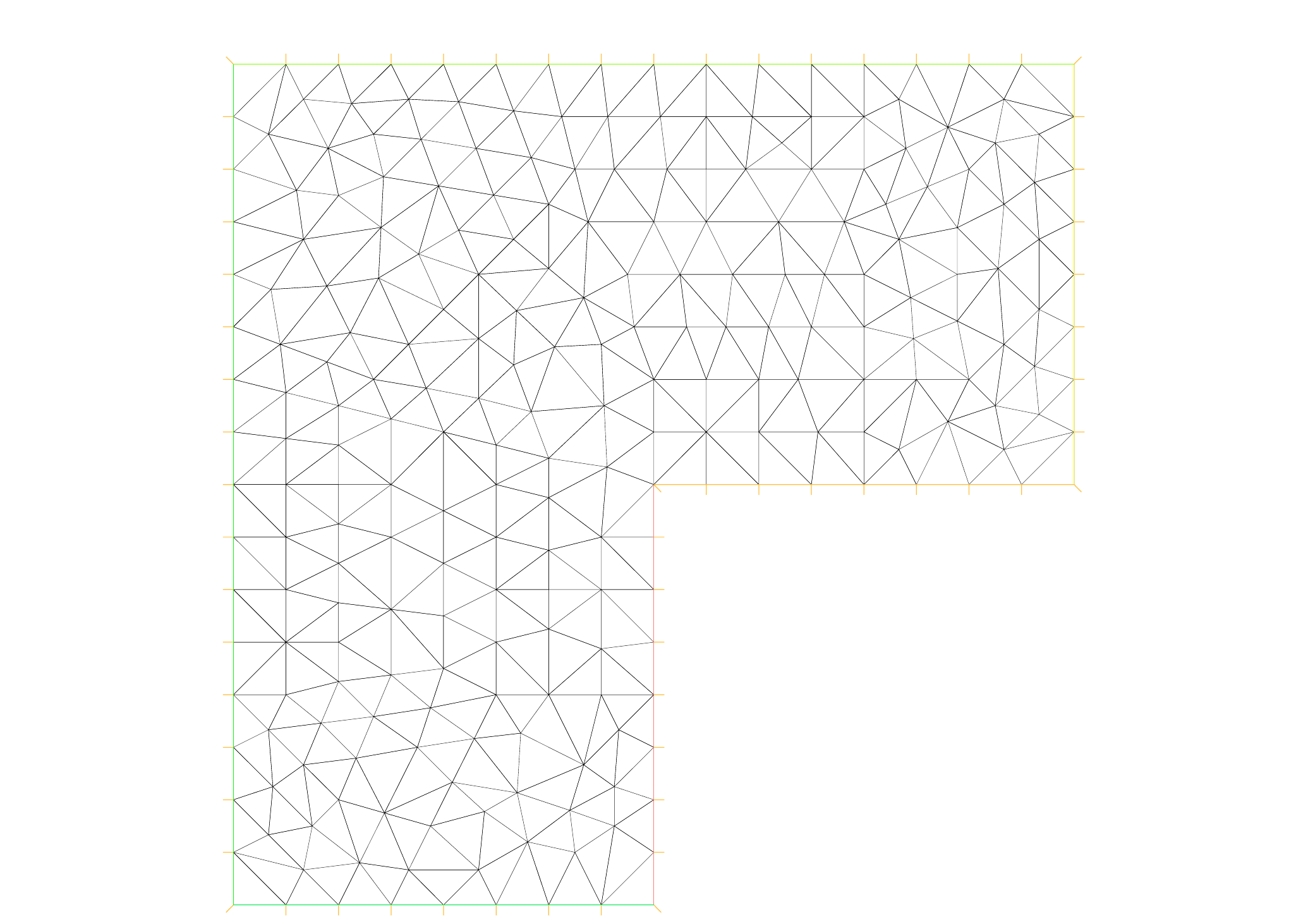}}
\subfigure[$h=1/32$] 
{\includegraphics[height=1.15in,width=1.6in]{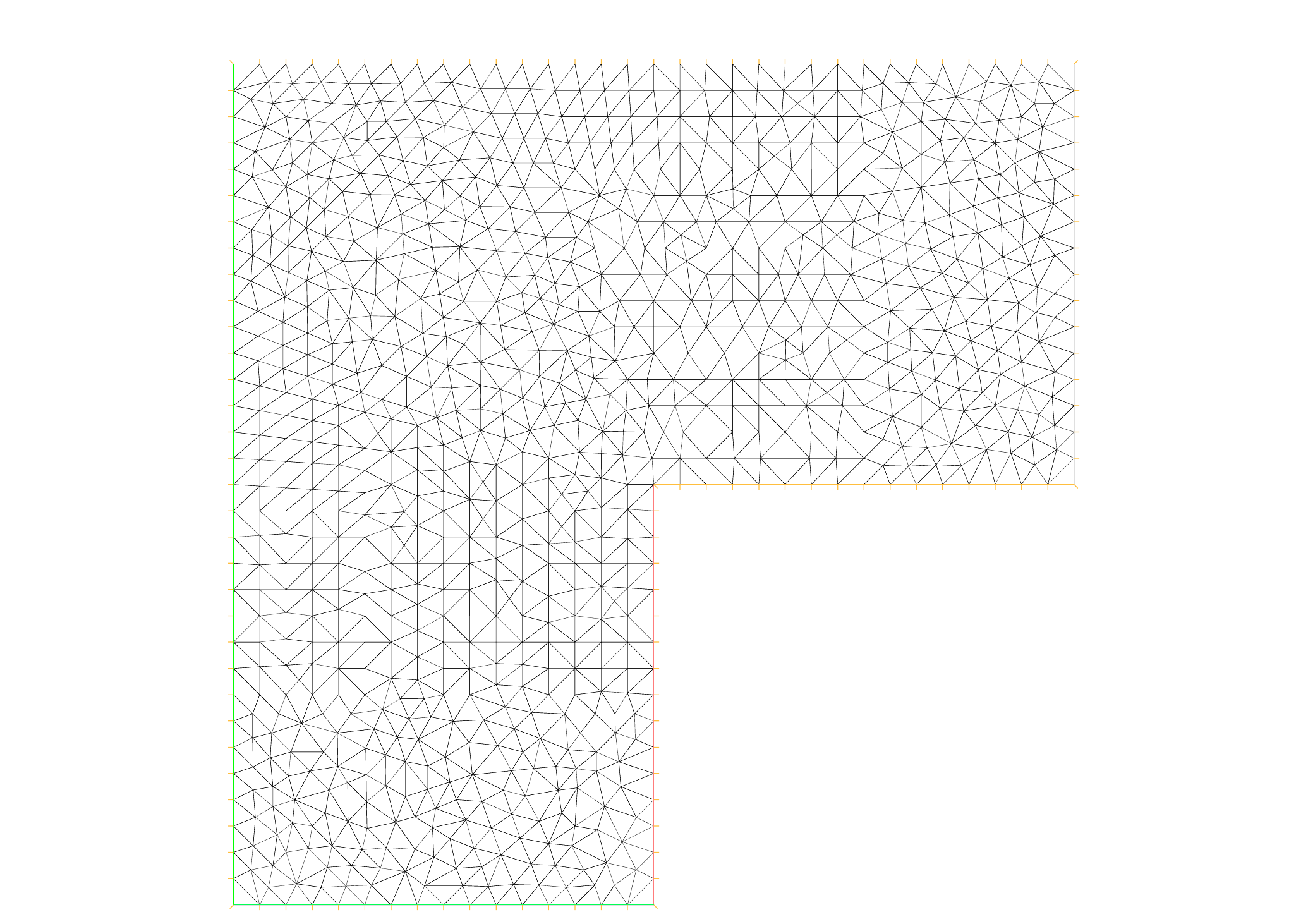}}
\subfigure[$h=1/64$]
{\includegraphics[height=1.15in,width=1.6in]{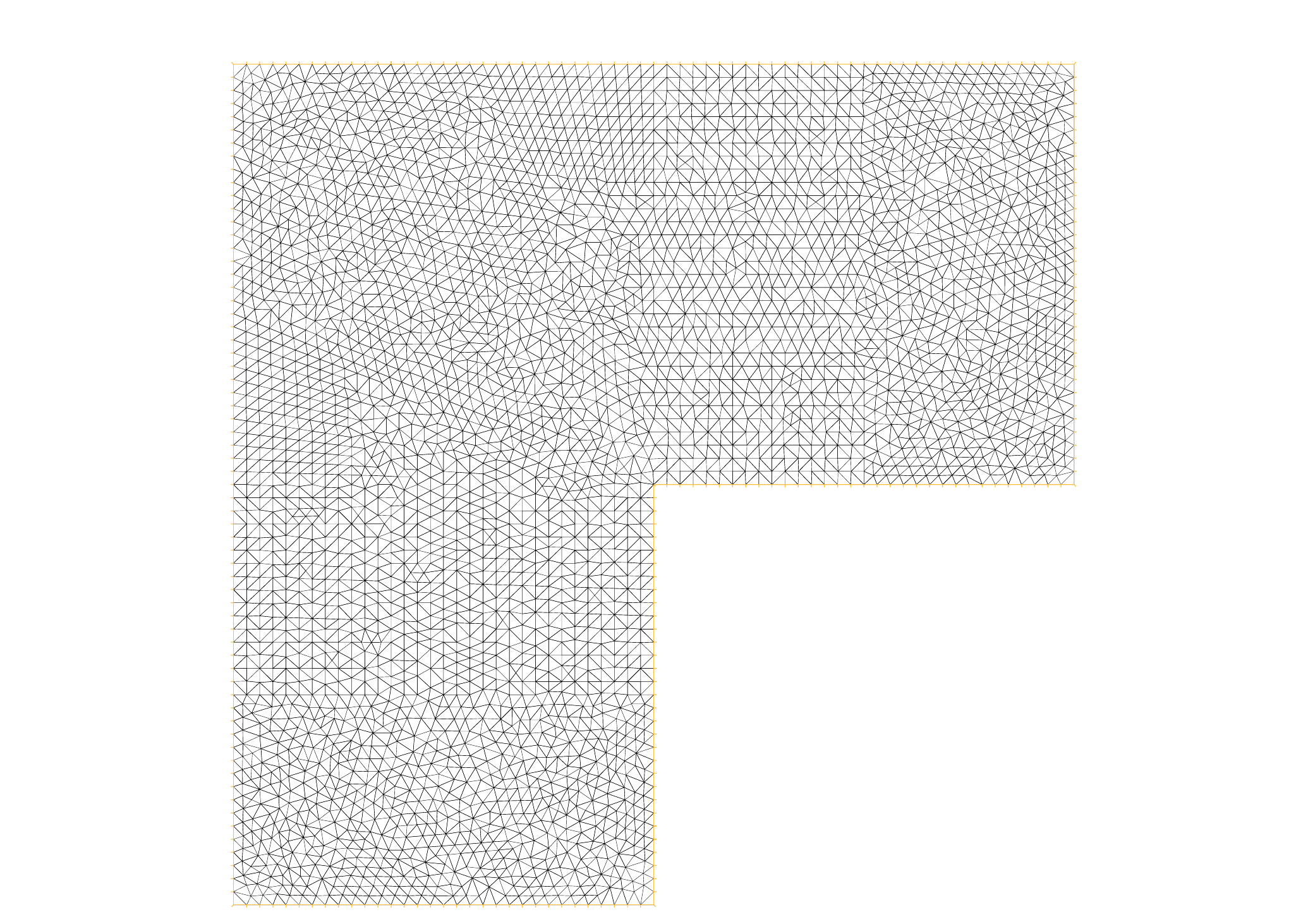}}
\caption{Quasi-uniform triangulations of the L-shape domain.}
\label{LshapeD}\vspace{10pt}
\subfigure[$t=5$] 
{\includegraphics[height=1.15in,width=1.6in]{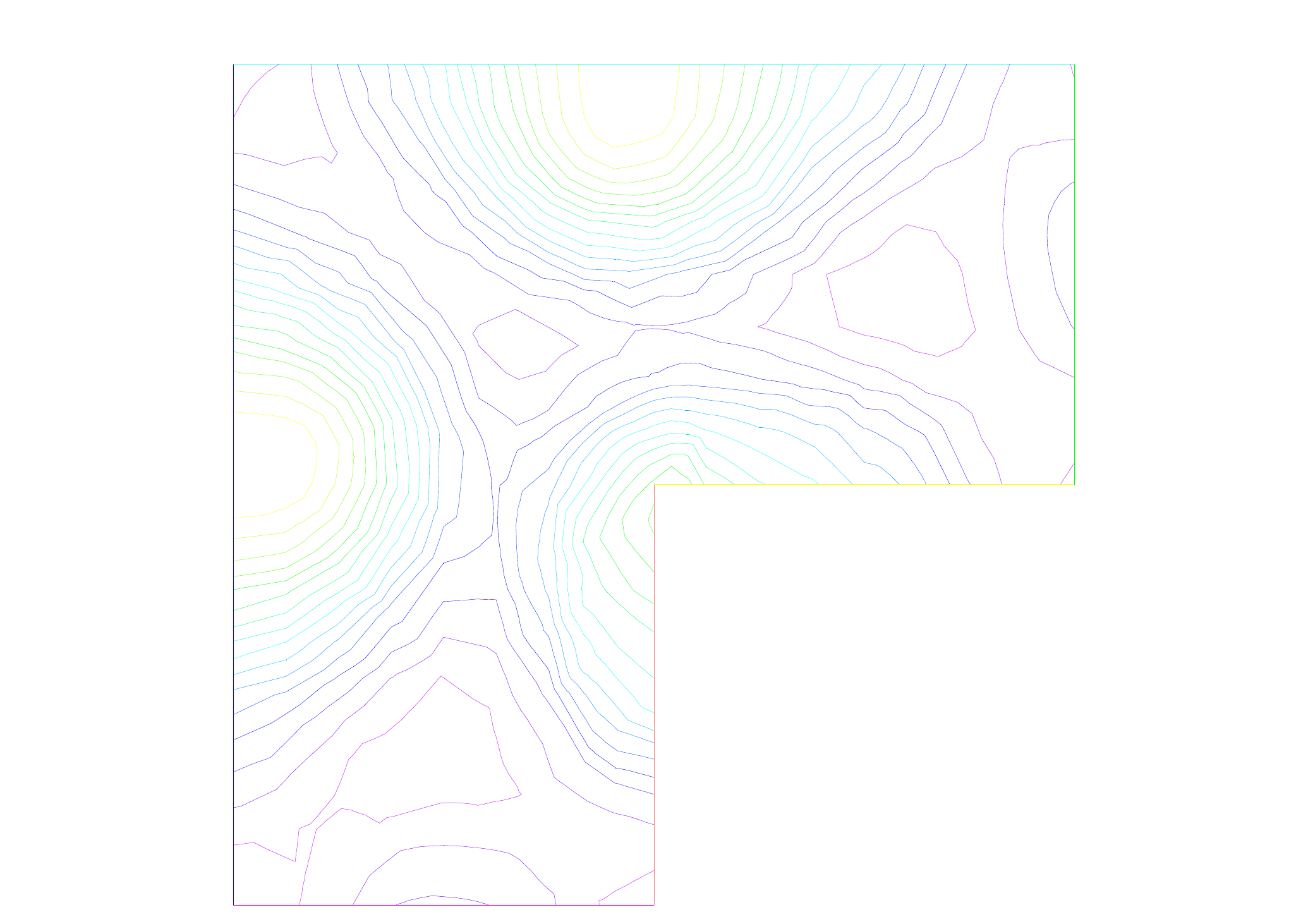}}
\subfigure[$t=20$] 
{\includegraphics[height=1.15in,width=1.6in]{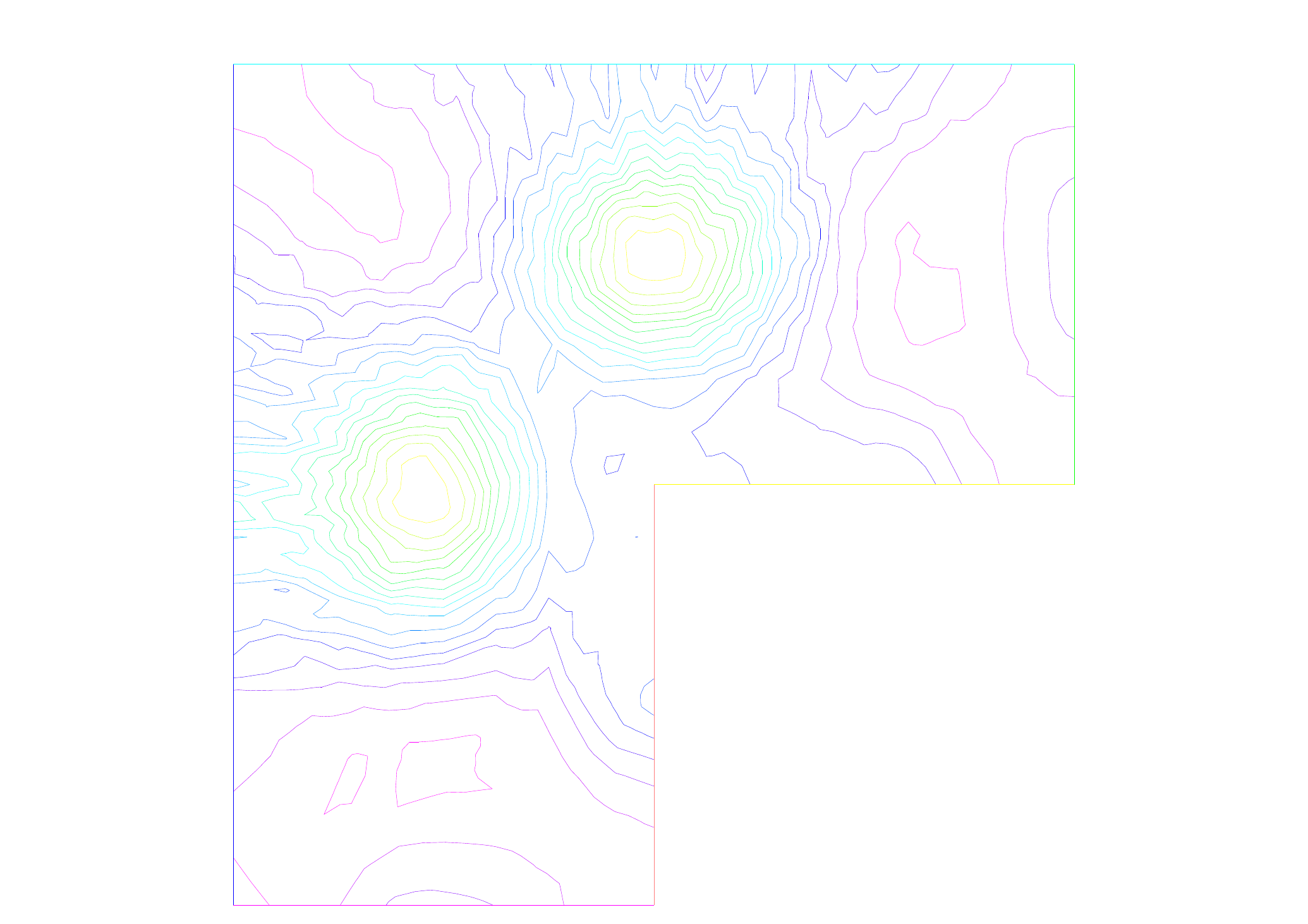}}
\subfigure[$t=40$]
{\includegraphics[height=1.15in,width=1.6in]{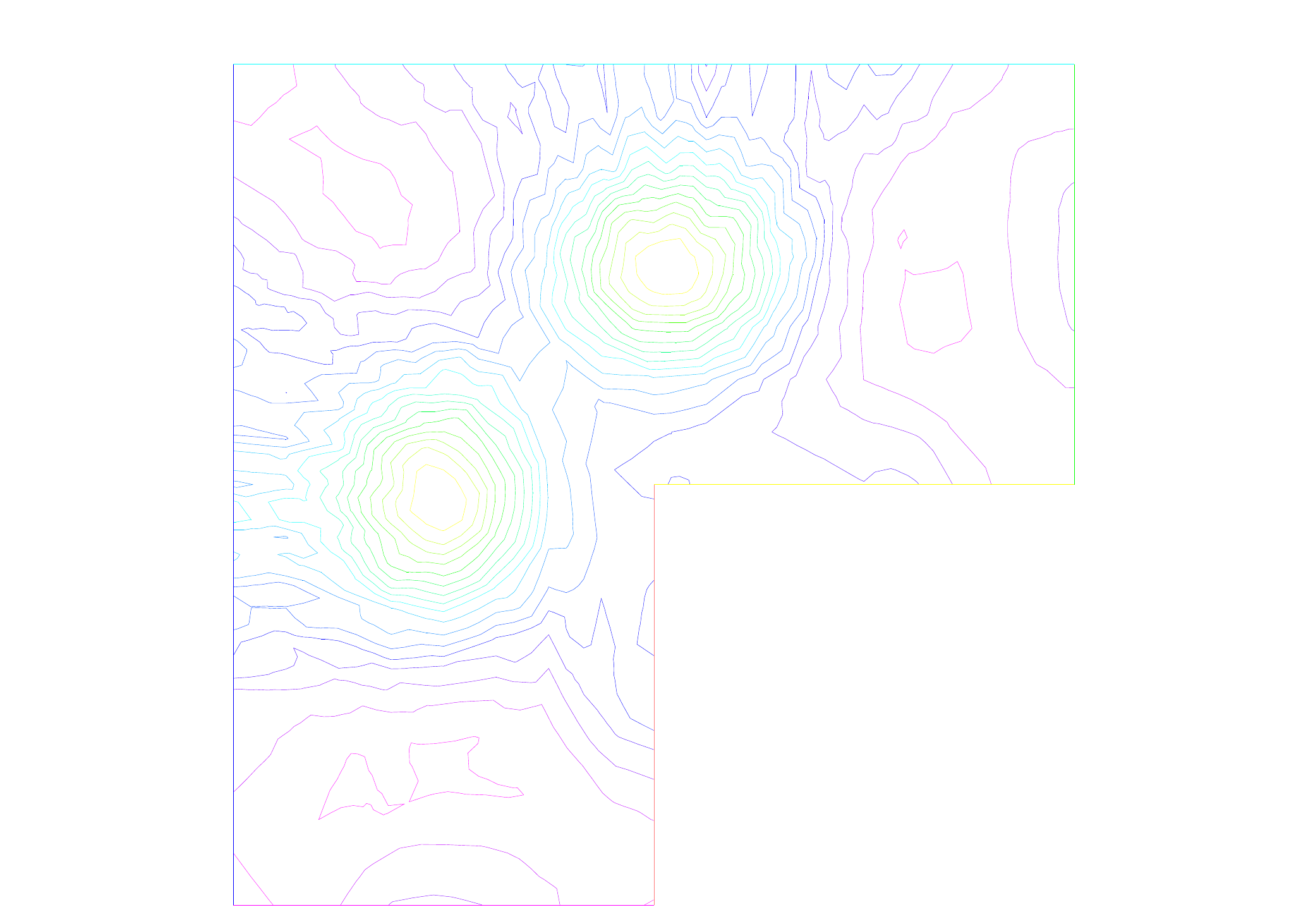}}
\caption{Contour of $|\psi|^2$ by solving 
the TDGL under the temporal gauge
with $h=1/16$.}
\label{LShp1}\vspace{10pt}
\subfigure[$t=5$] 
{\includegraphics[height=1.15in,width=1.6in]{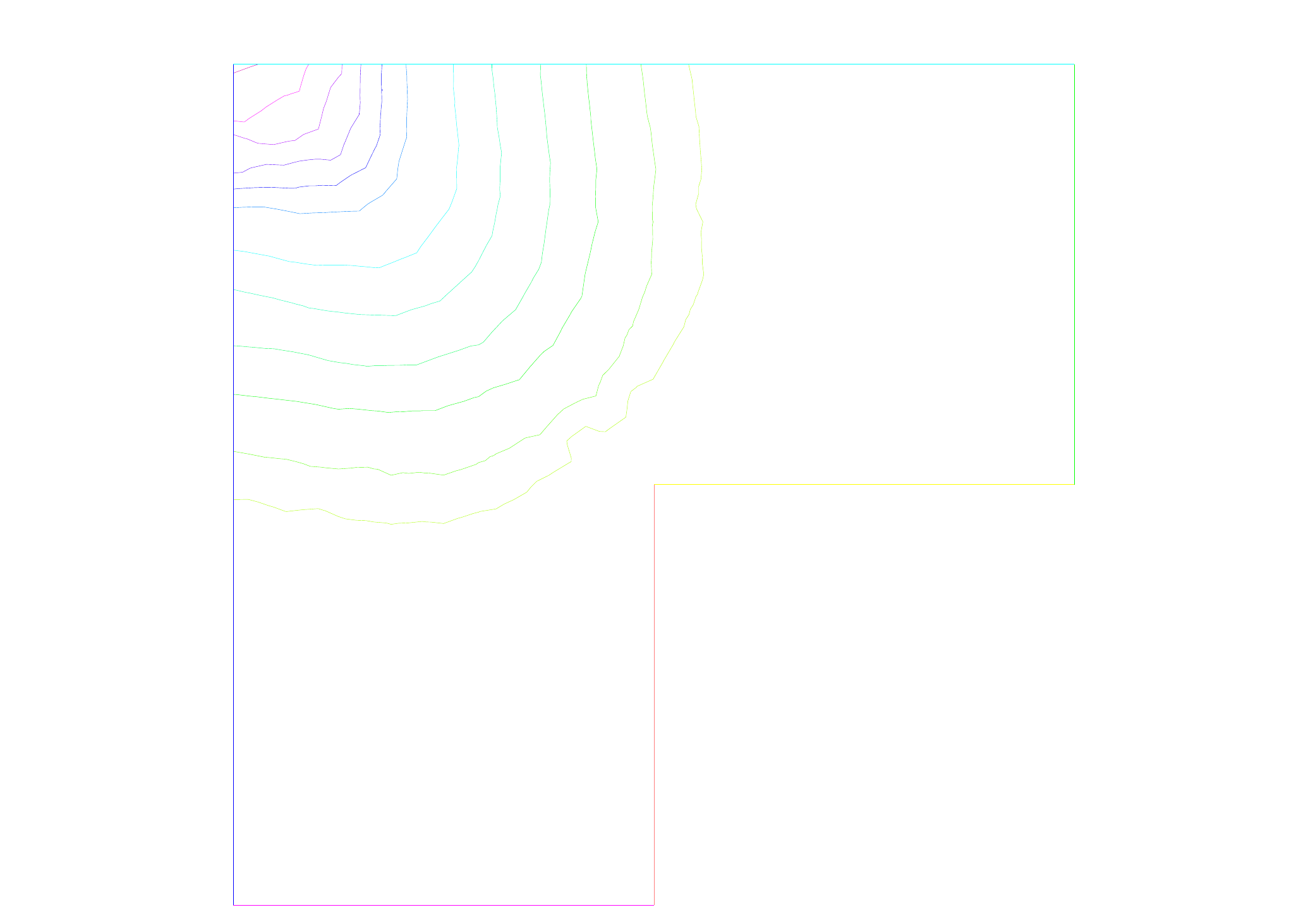}}
\subfigure[$t=20$] 
{\includegraphics[height=1.15in,width=1.6in]{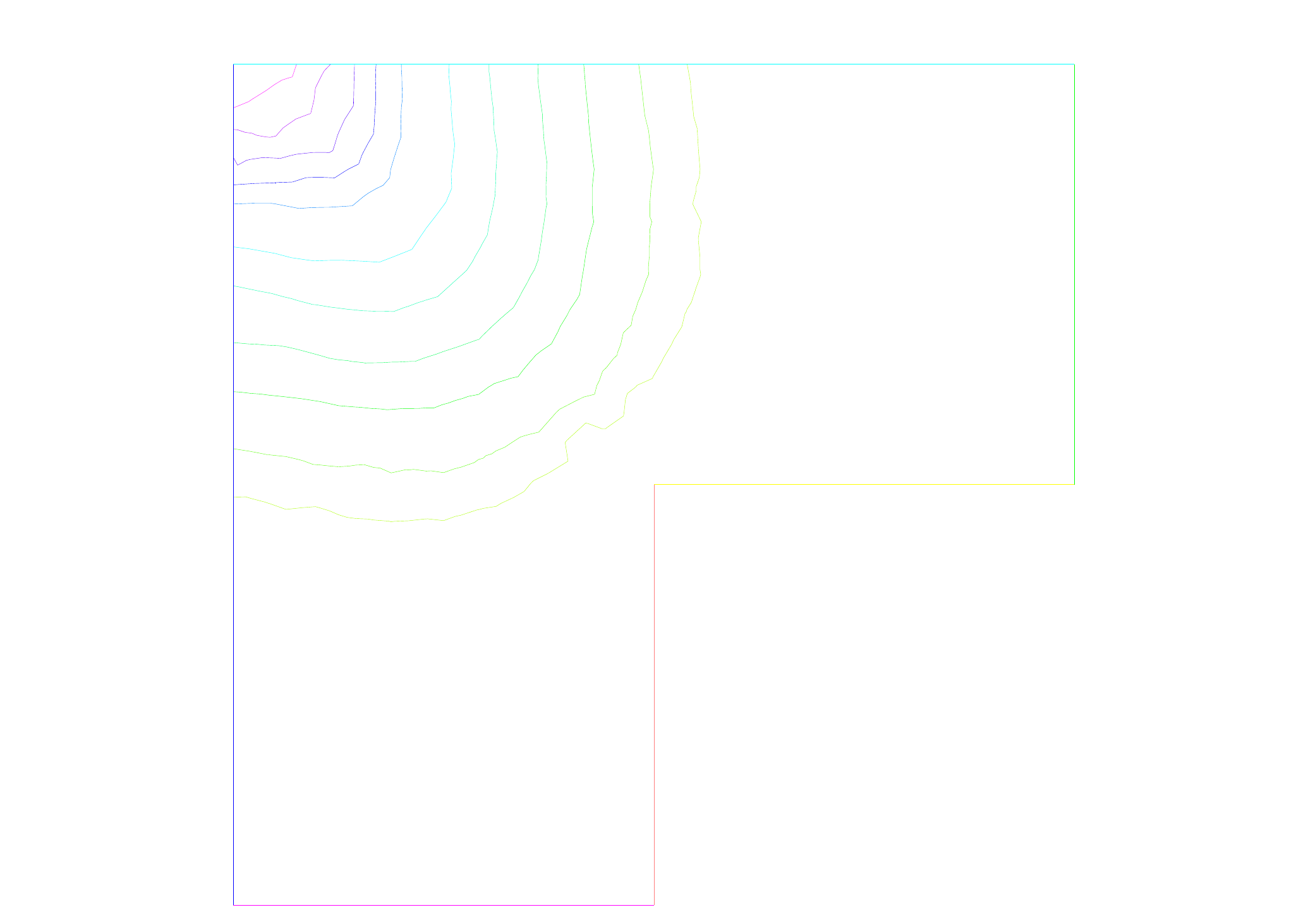}}
\subfigure[$t=40$] 
{\includegraphics[height=1.15in,width=1.6in]{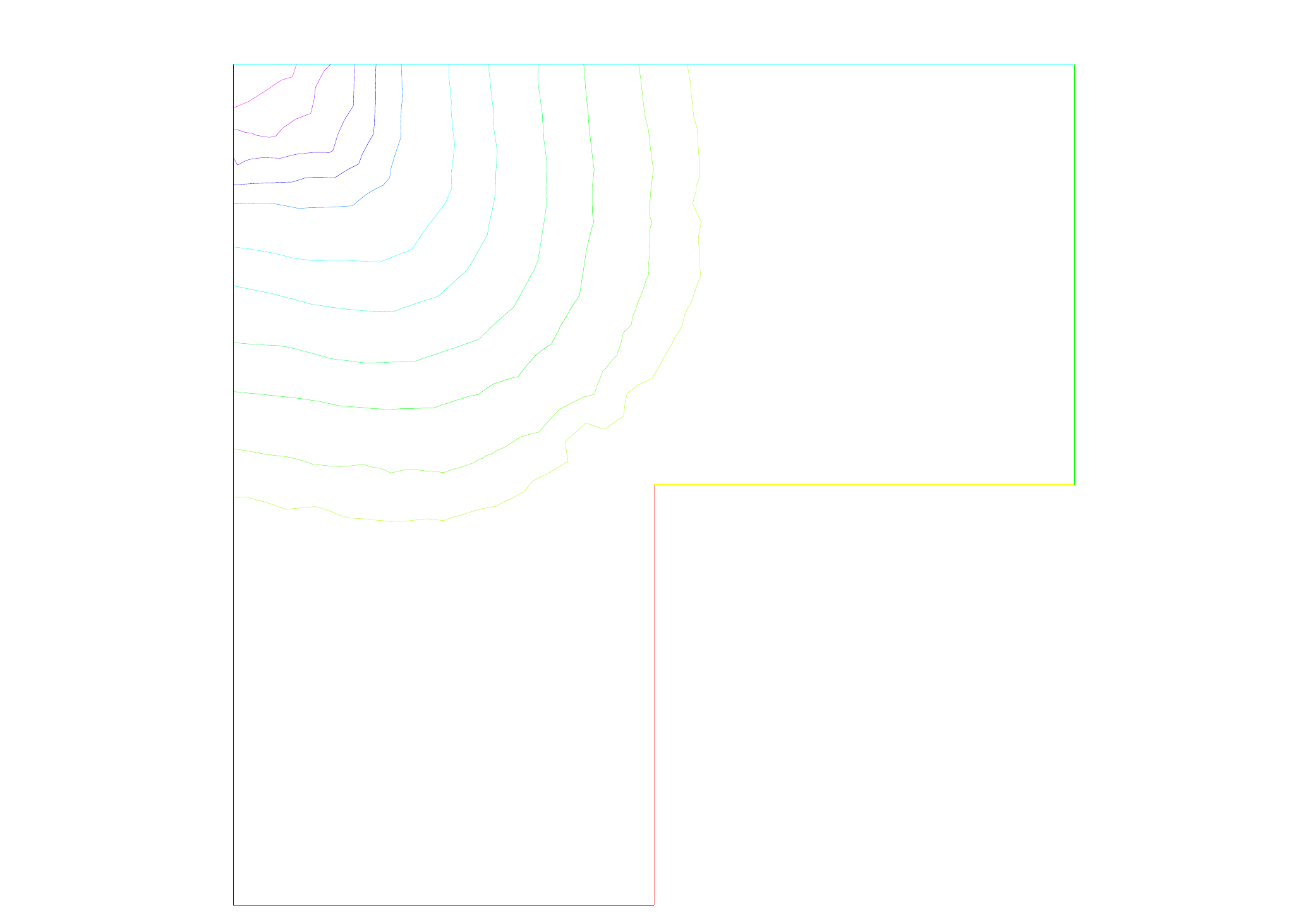}}
\caption{Contour of $|\psi|^2$ by solving 
the TDGL under the Lorentz gauge
with $h=1/16$.}
\label{LShp2}\vspace{10pt}
\subfigure[$t=5$] 
{\includegraphics[height=1.15in,width=1.6in]{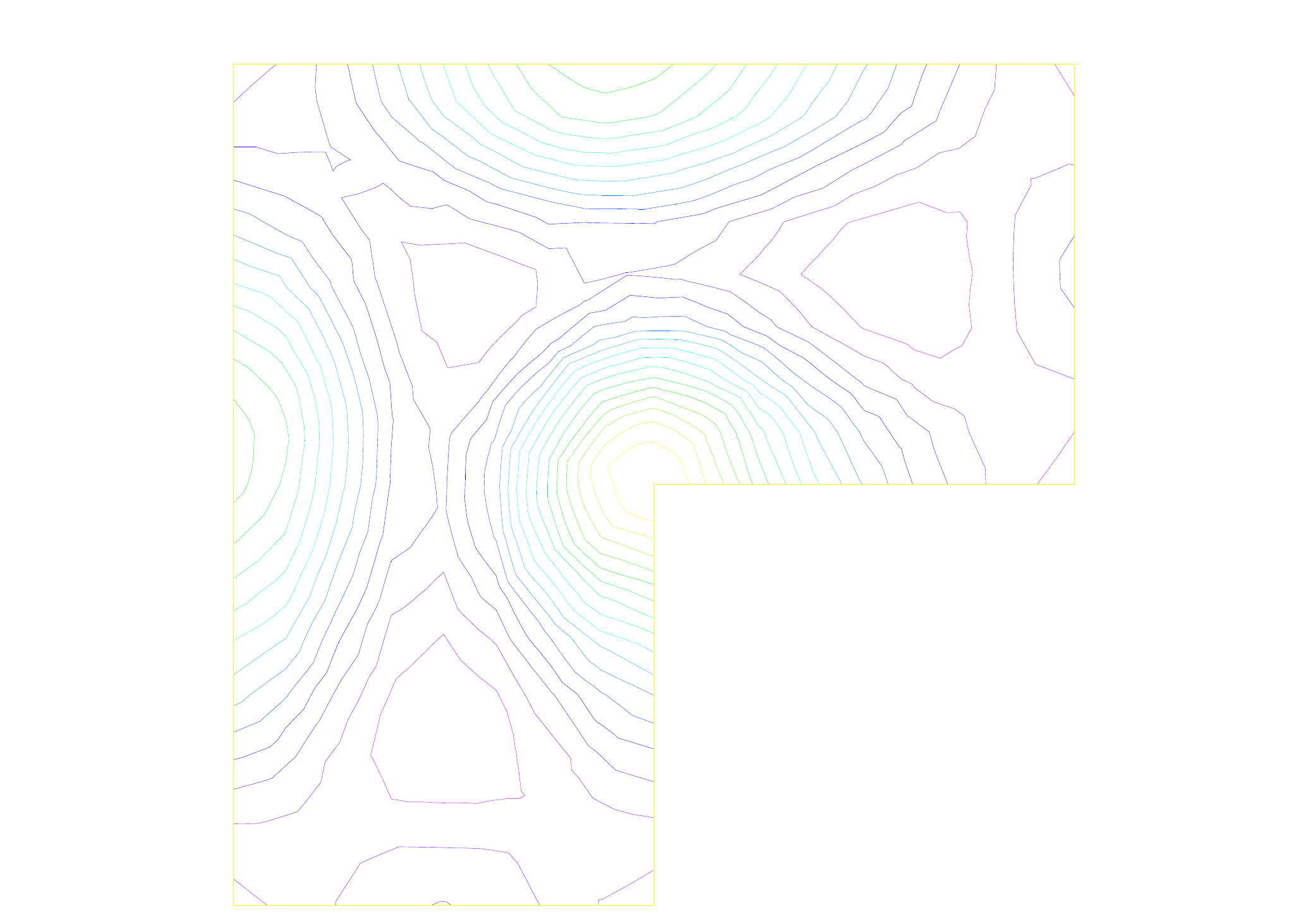}}
\subfigure[$t=20$] 
{\includegraphics[height=1.15in,width=1.6in]{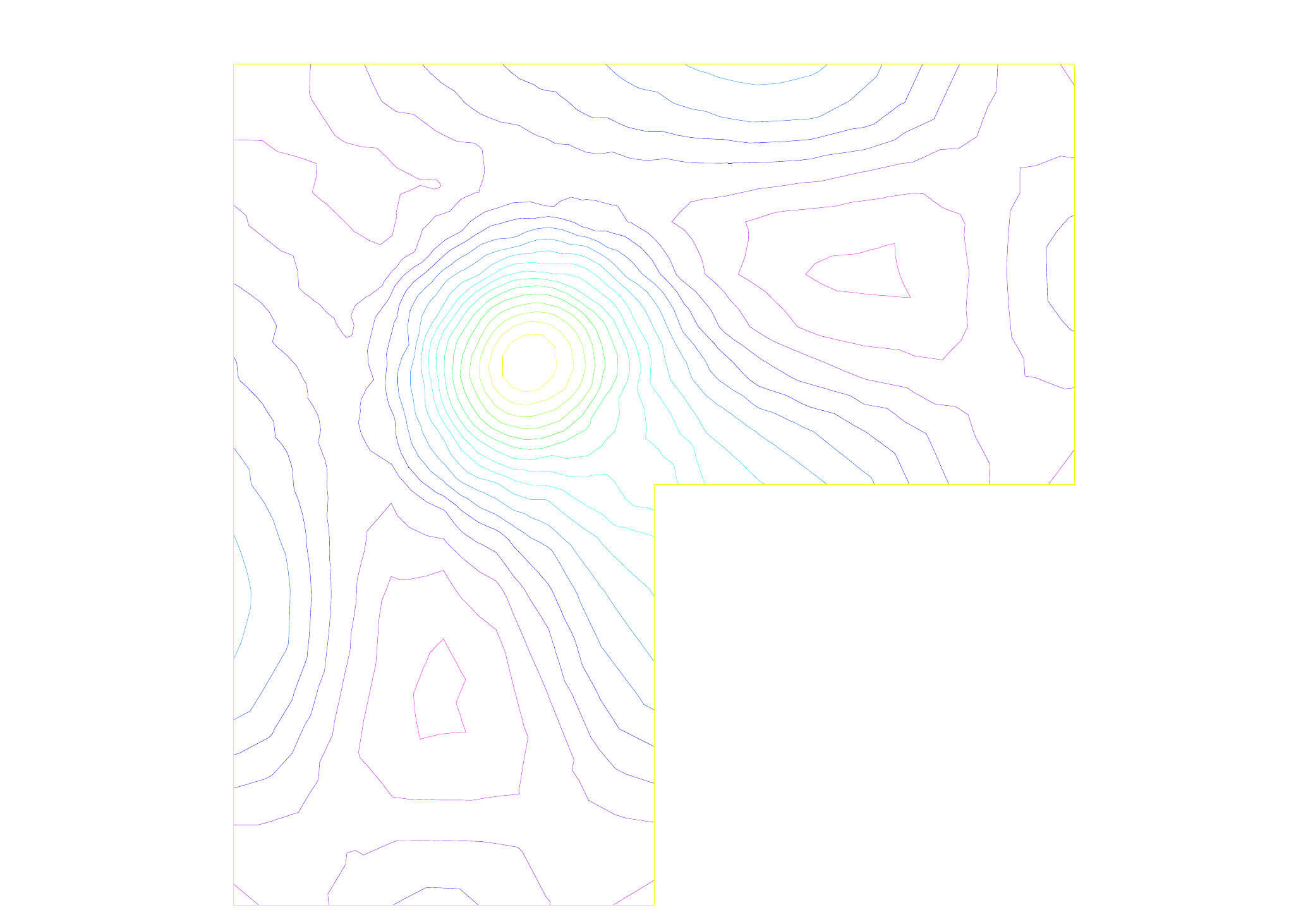}}
\subfigure[$t=40$] 
{\includegraphics[height=1.15in,width=1.6in]{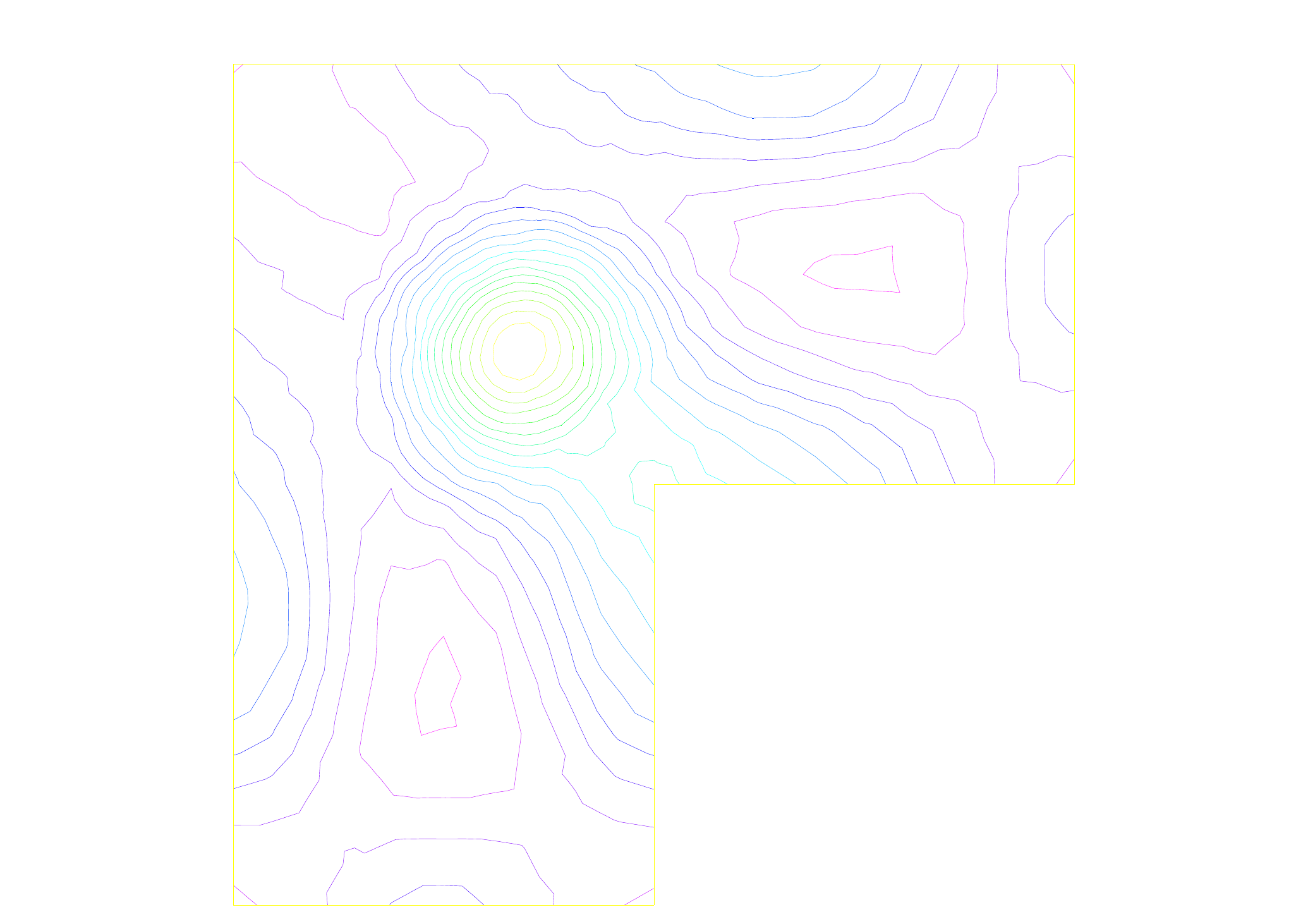}}
\caption{Contour of $|\psi|^2$ computed by the new approach
with $h=1/16$.}
\label{LShp3}

\end{figure}

{\it Example 3.1}$\quad$
Firstly, we simulate the the vortex dynamics 
in an L-shape domain 
whose longest side has unit length, centered at the origin,
with 
$\eta=1$, $\kappa=10$ and  
$$
\psi_0 = 0.6 + 0.8i, \qquad {\bf A}_0 = (0,0), \qquad
f = 5  .
$$
The L-shape domain is triangulated quasi-uniformly, 
as shown in Figure \ref{LshapeD}, 
with $M$ nodes per unit length on each side, 
and we denote by $h=1/M$ for simplicity.
 
\begin{figure}[htp]
\centering
\subfigure[$t=5$] 
{\includegraphics[height=1.15in,width=1.6in]{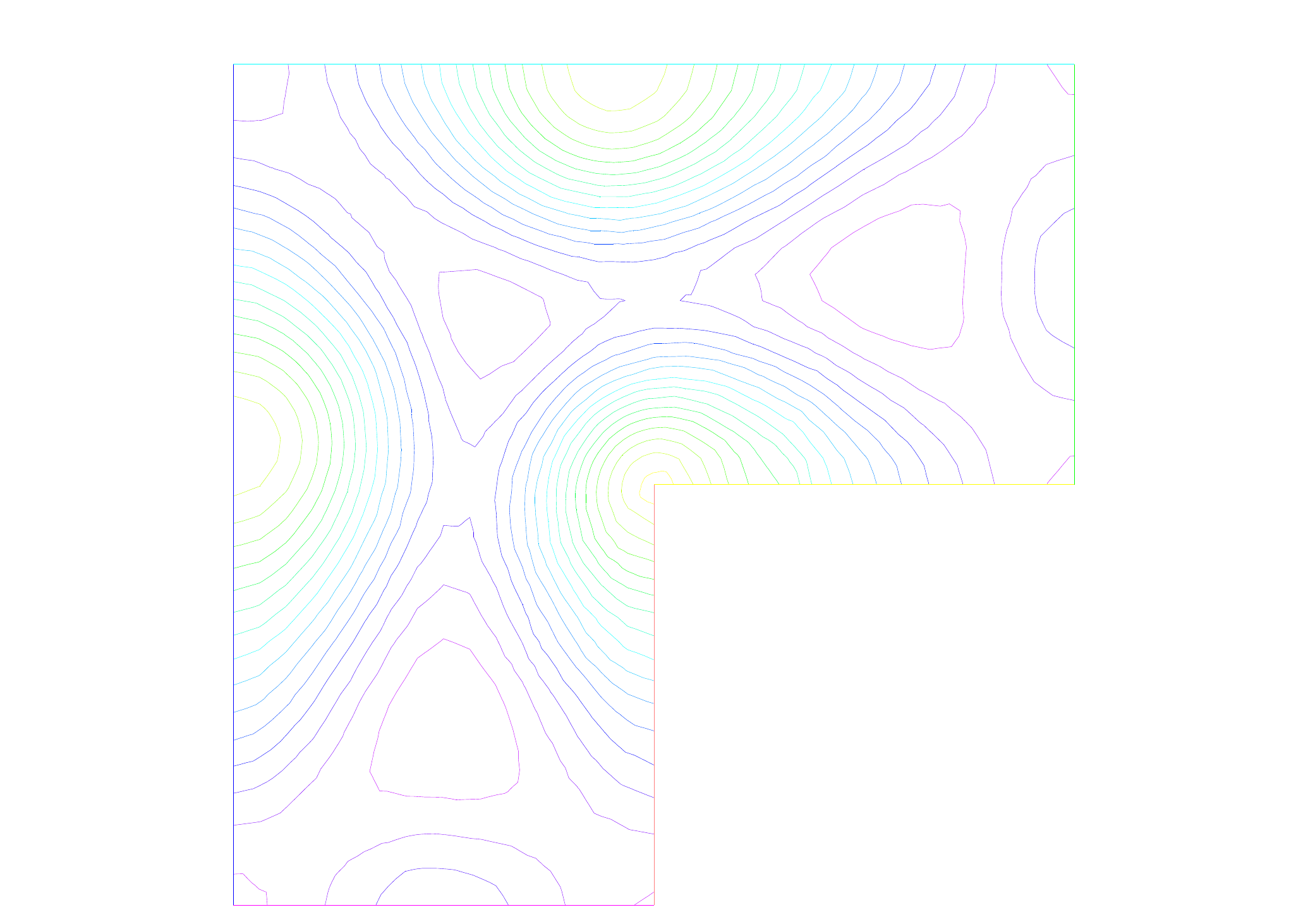}}
\subfigure[$t=20$] 
{\includegraphics[height=1.15in,width=1.6in]{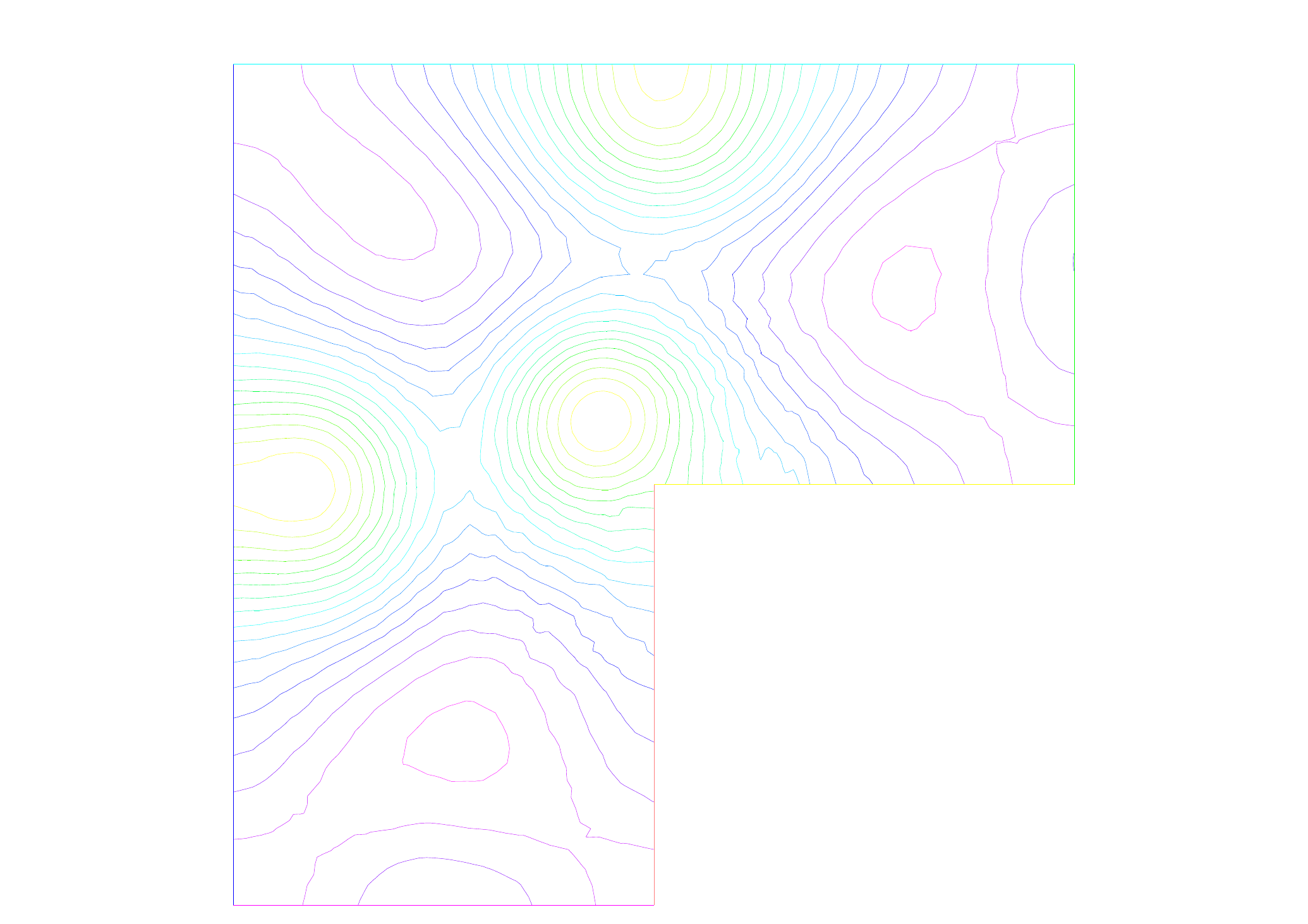}}
\subfigure[$t=40$]
{\includegraphics[height=1.15in,width=1.6in]{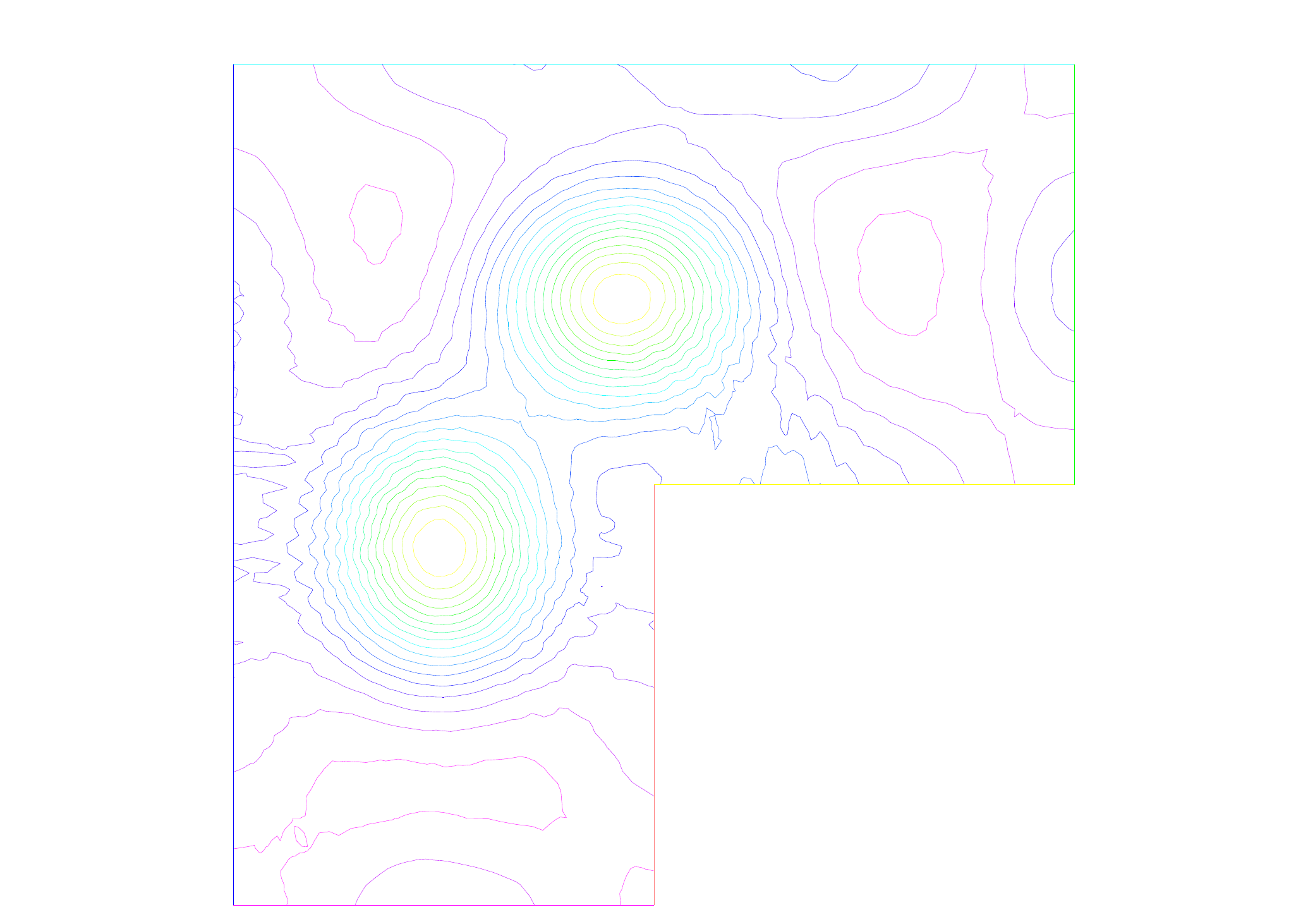}}
\caption{Contour of $|\psi|^2$ by solving 
the TDGL under the temporal gauge
with $h=1/32$.}
\label{LShp4}\vspace{10pt}
\subfigure[$t=5$] 
{\includegraphics[height=1.15in,width=1.6in]{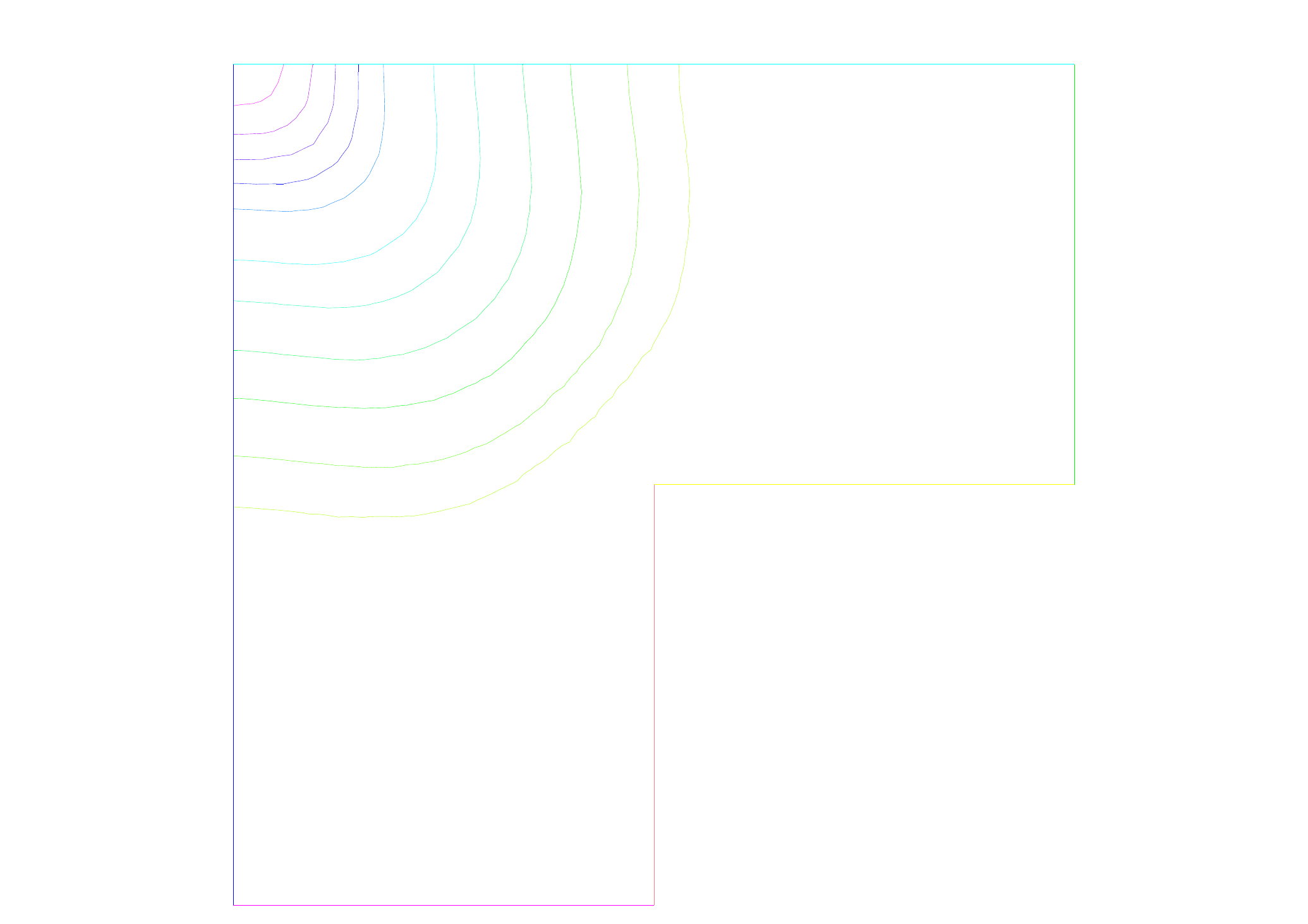}}
\subfigure[$t=20$] 
{\includegraphics[height=1.15in,width=1.6in]{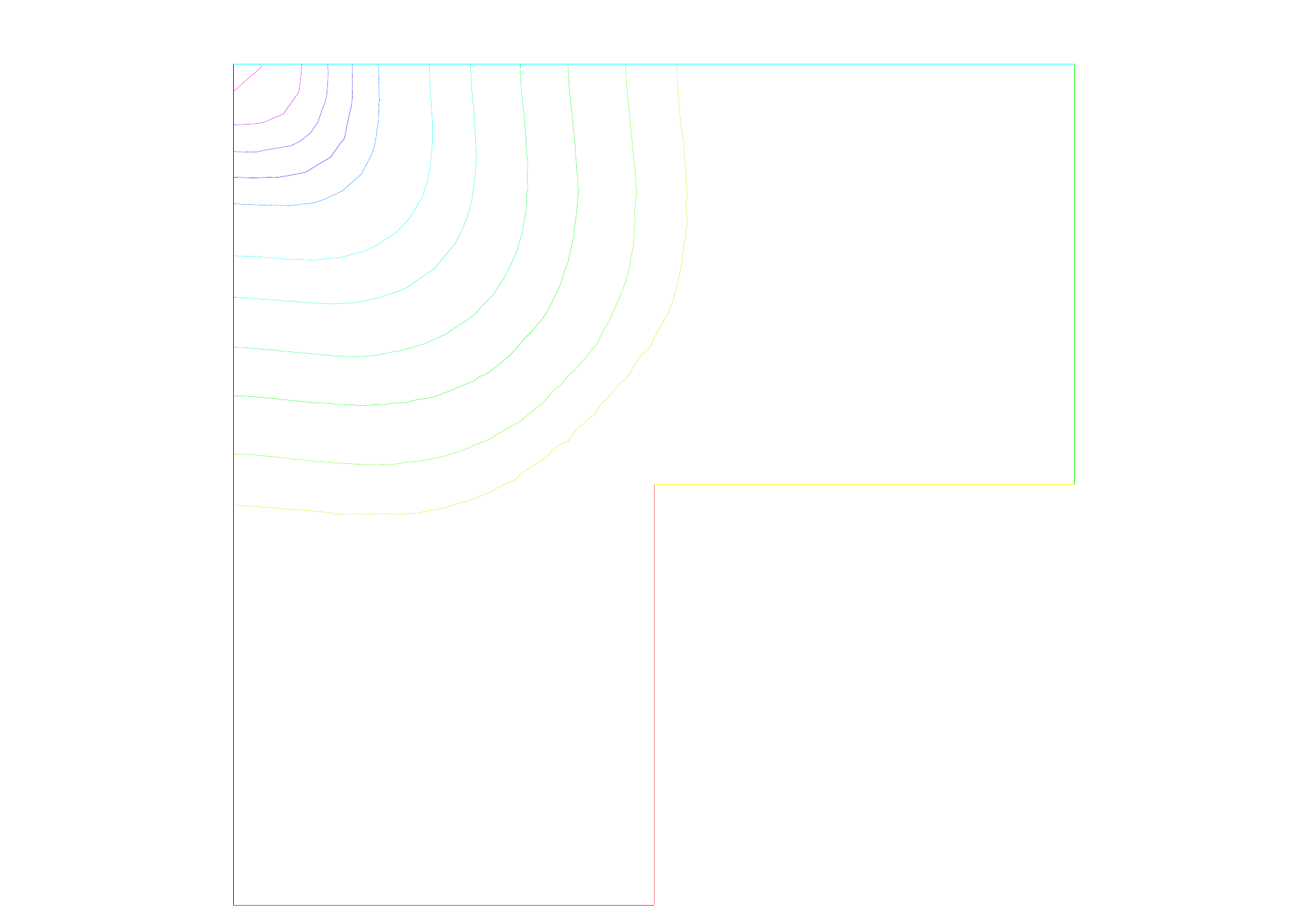}}
\subfigure[$t=40$] 
{\includegraphics[height=1.15in,width=1.6in]{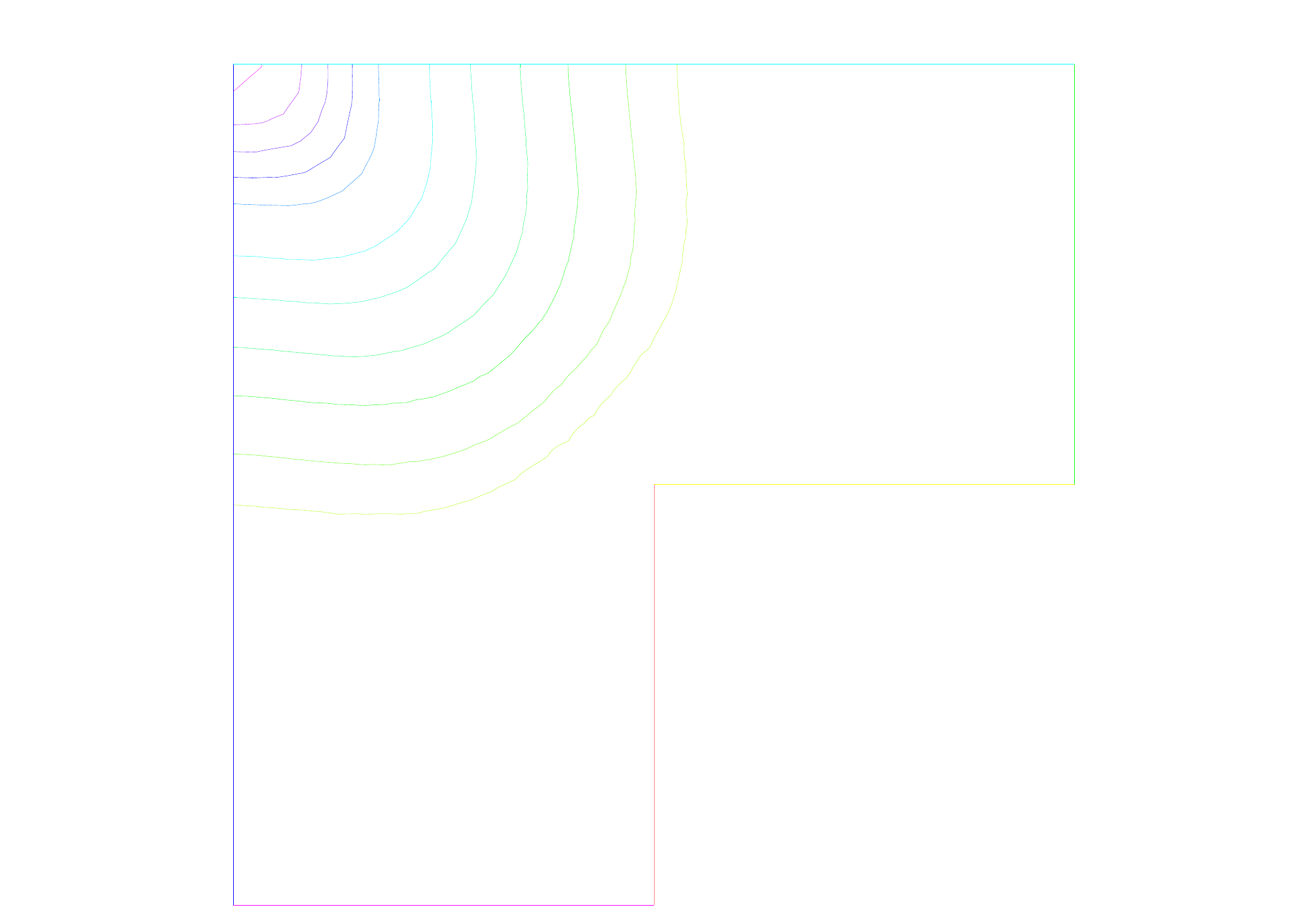}}
\caption{Contour of $|\psi|^2$ by solving 
the TDGL under the Lorentz gauge
with $h=1/32$.}
\label{LShp5}\vspace{10pt}
\subfigure[$t=5$] 
{\includegraphics[height=1.15in,width=1.6in]{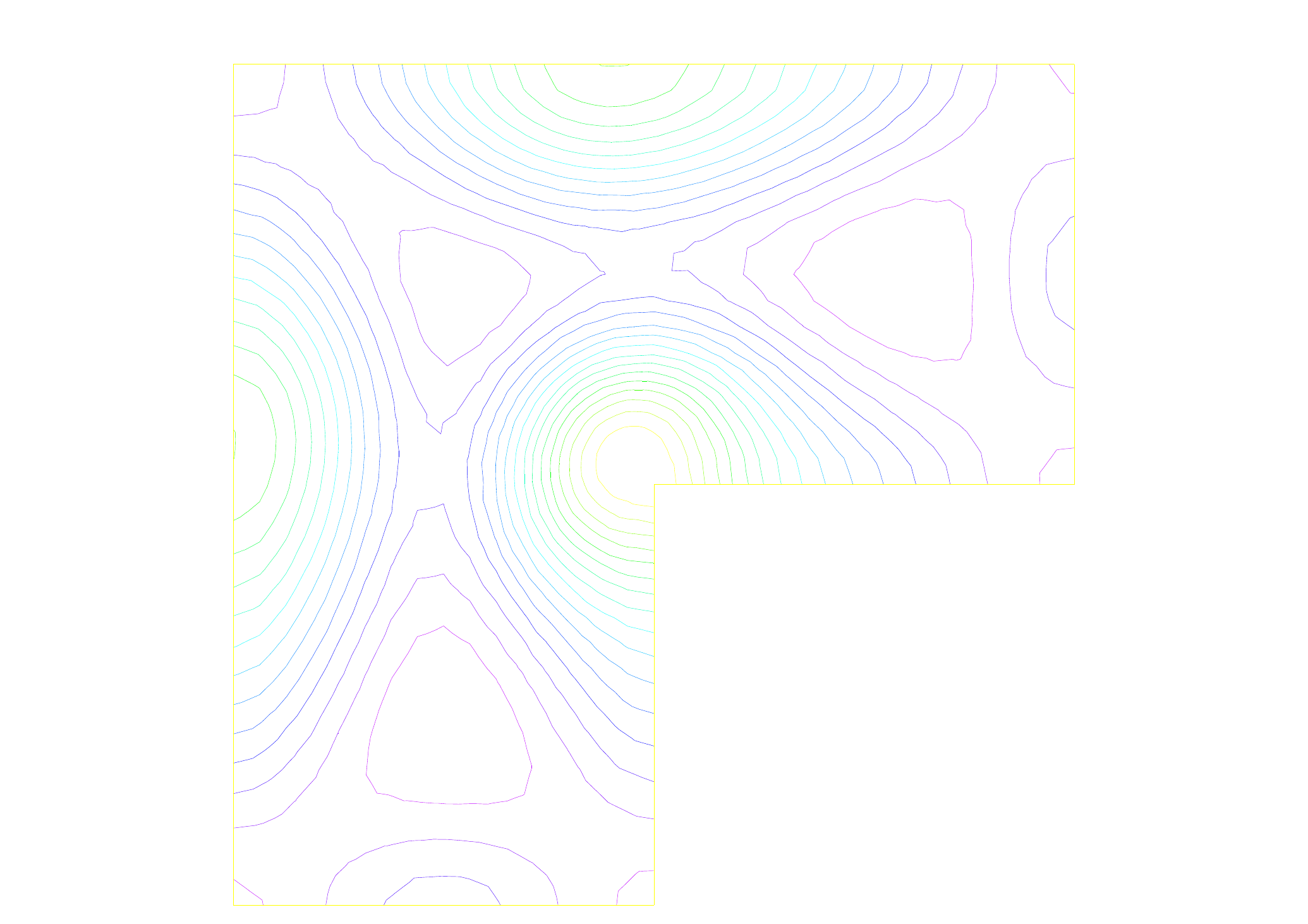}}
\subfigure[$t=20$] 
{\includegraphics[height=1.15in,width=1.6in]{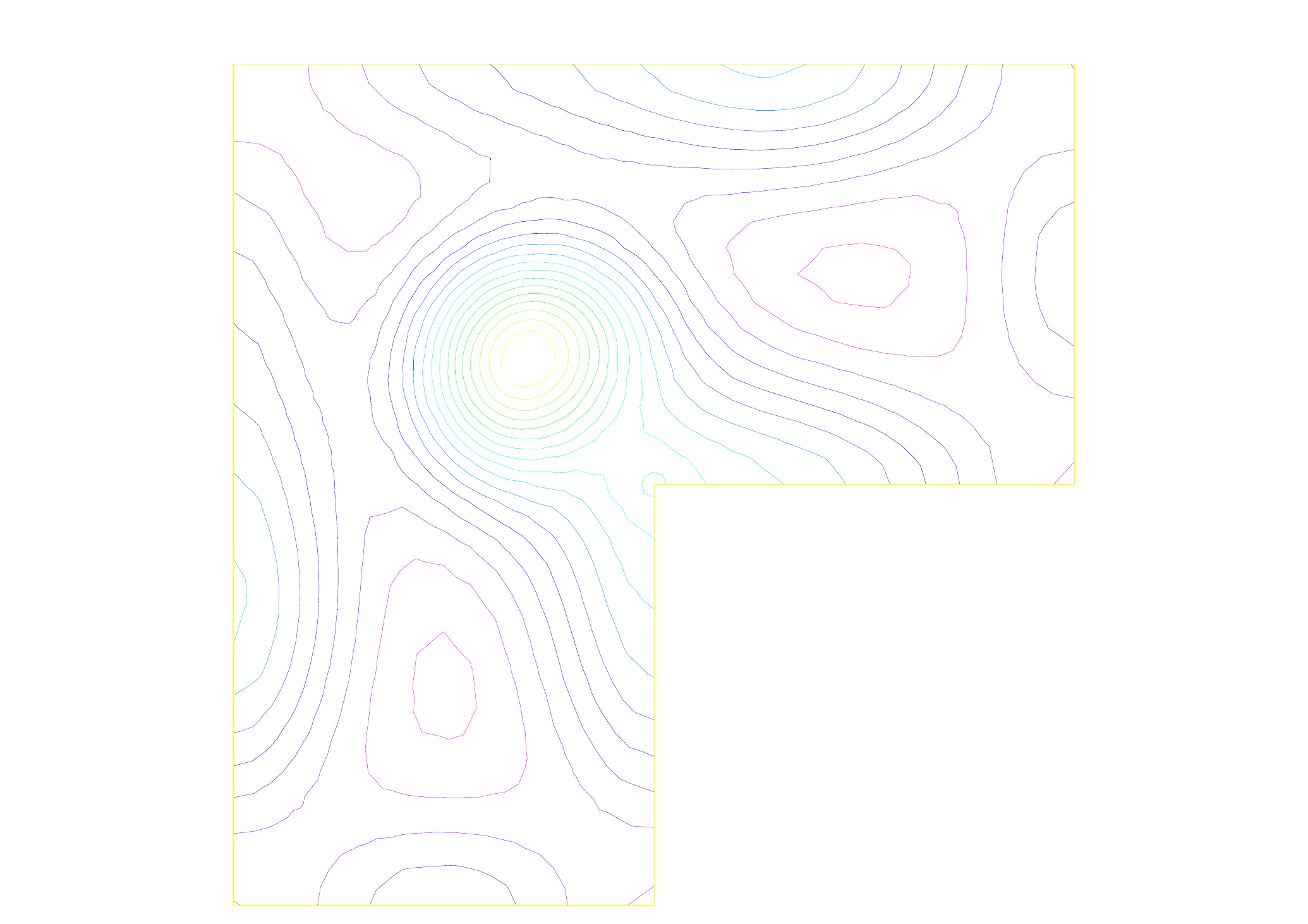}}
\subfigure[$t=40$] 
{\includegraphics[height=1.15in,width=1.6in]{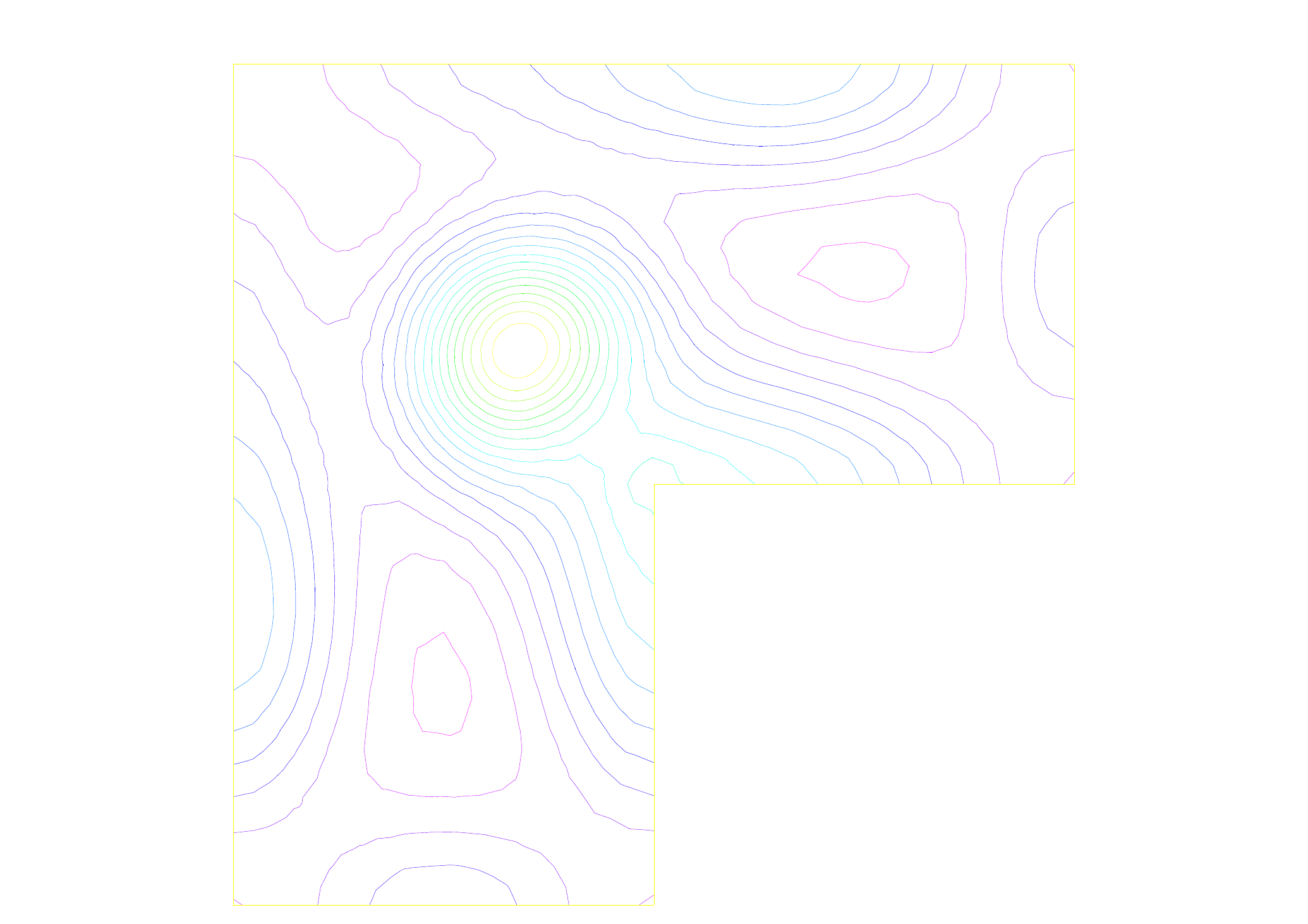}}
\caption{Contour of $|\psi|^2$ computed by the new approach
with $h=1/32$.}
\label{LShp6}
\end{figure}

\begin{figure}[htp]
\centering
\subfigure[$t=5$] 
{\includegraphics[height=1.15in,width=1.6in]{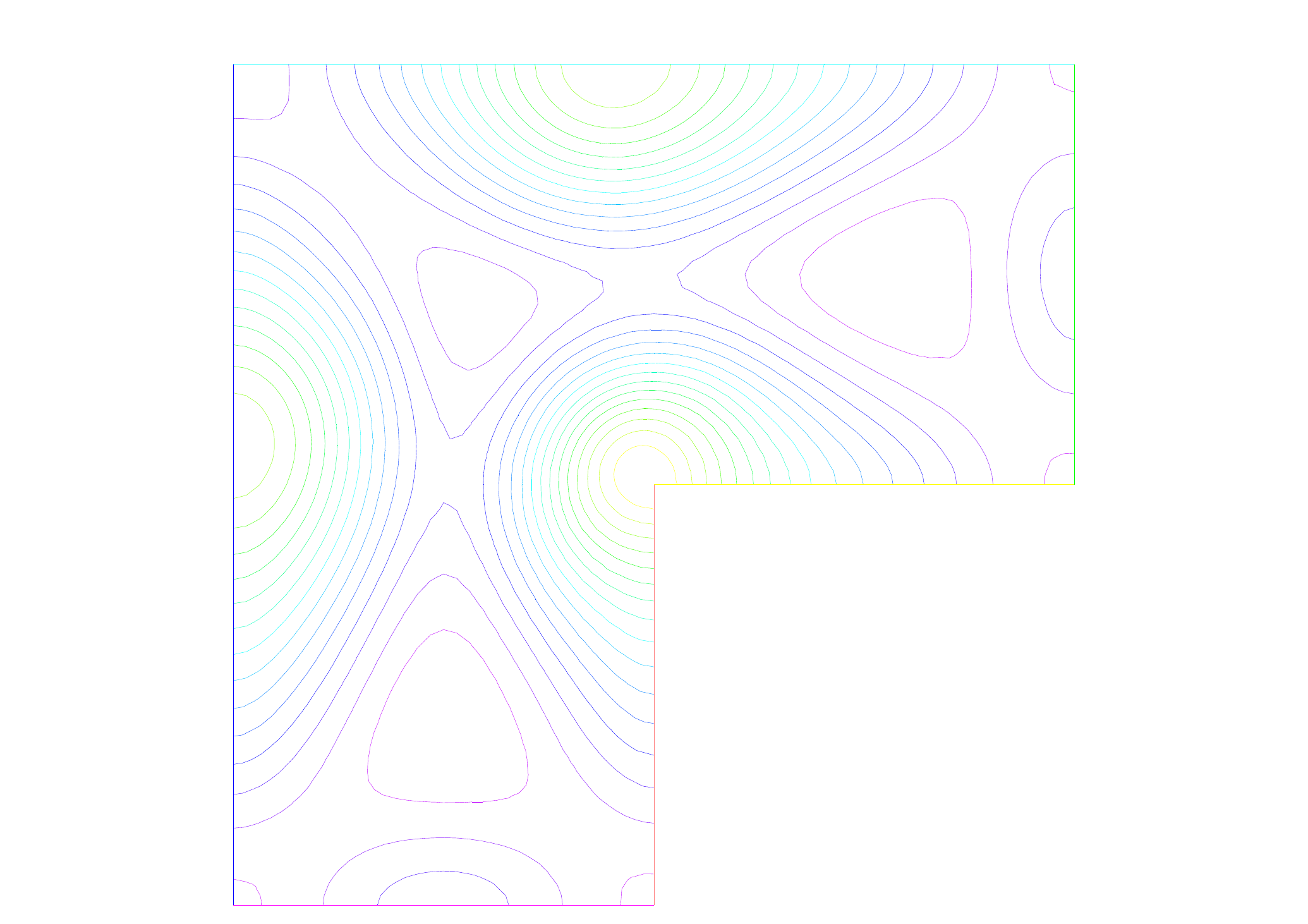}}
\subfigure[$t=20$] 
{\includegraphics[height=1.15in,width=1.6in]{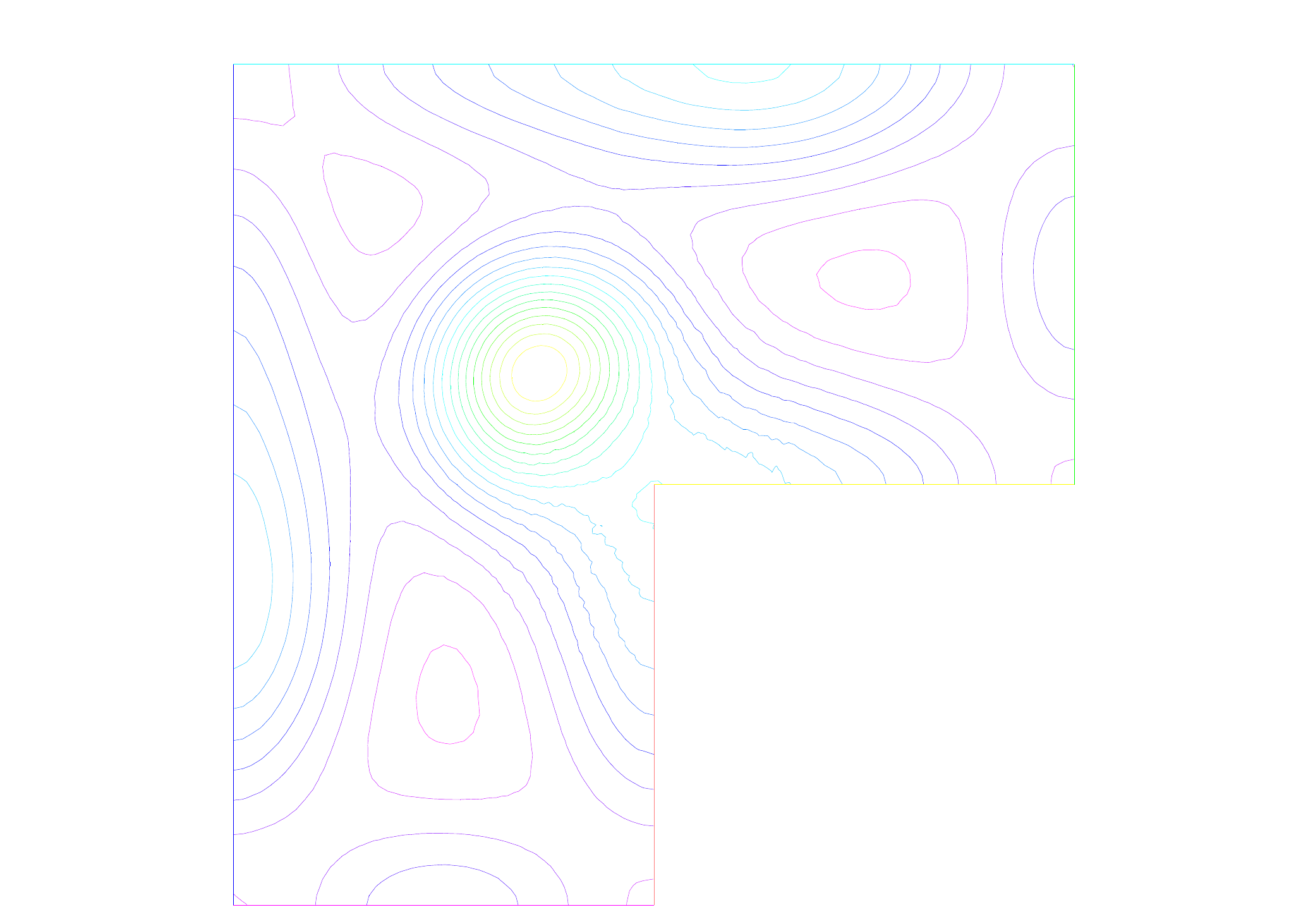}}
\subfigure[$t=40$]
{\includegraphics[height=1.15in,width=1.6in]{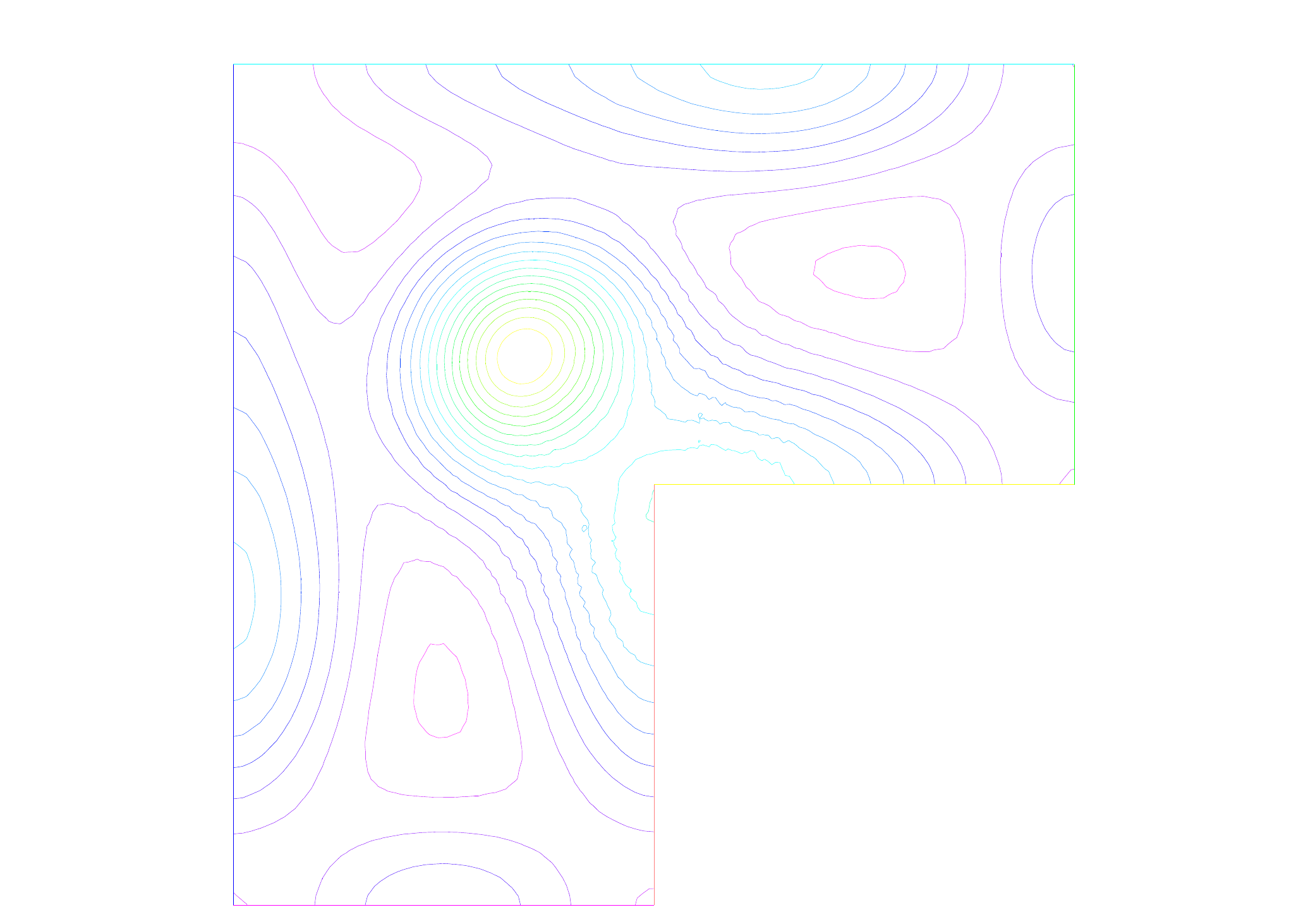}}
\caption{Contour of $|\psi|^2$ by solving 
the TDGL under the temporal gauge
with $h=1/64$.}
\label{LShp7}\vspace{10pt}
\subfigure[$t=5$] 
{\includegraphics[height=1.15in,width=1.6in]{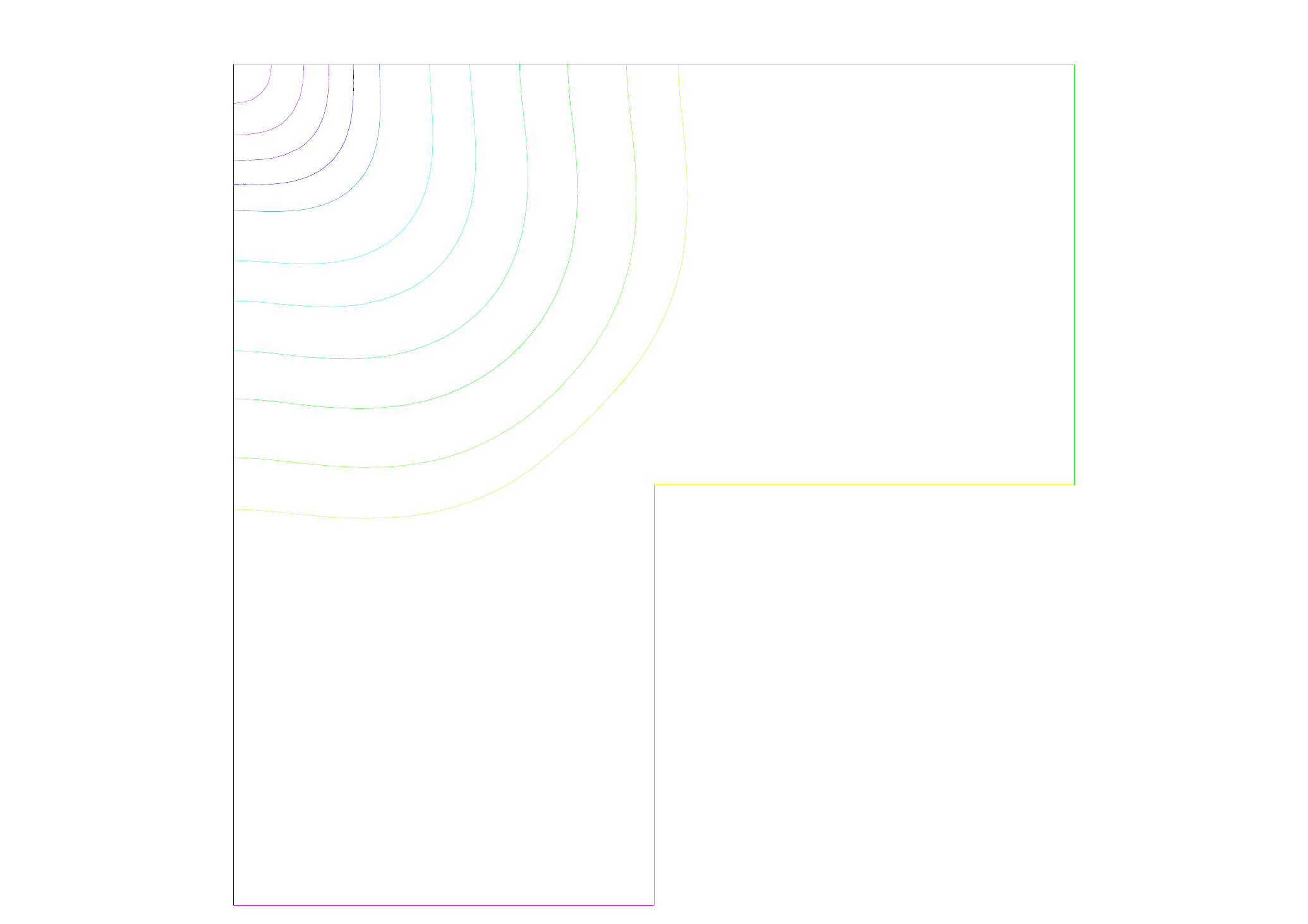}}
\subfigure[$t=20$] 
{\includegraphics[height=1.15in,width=1.6in]{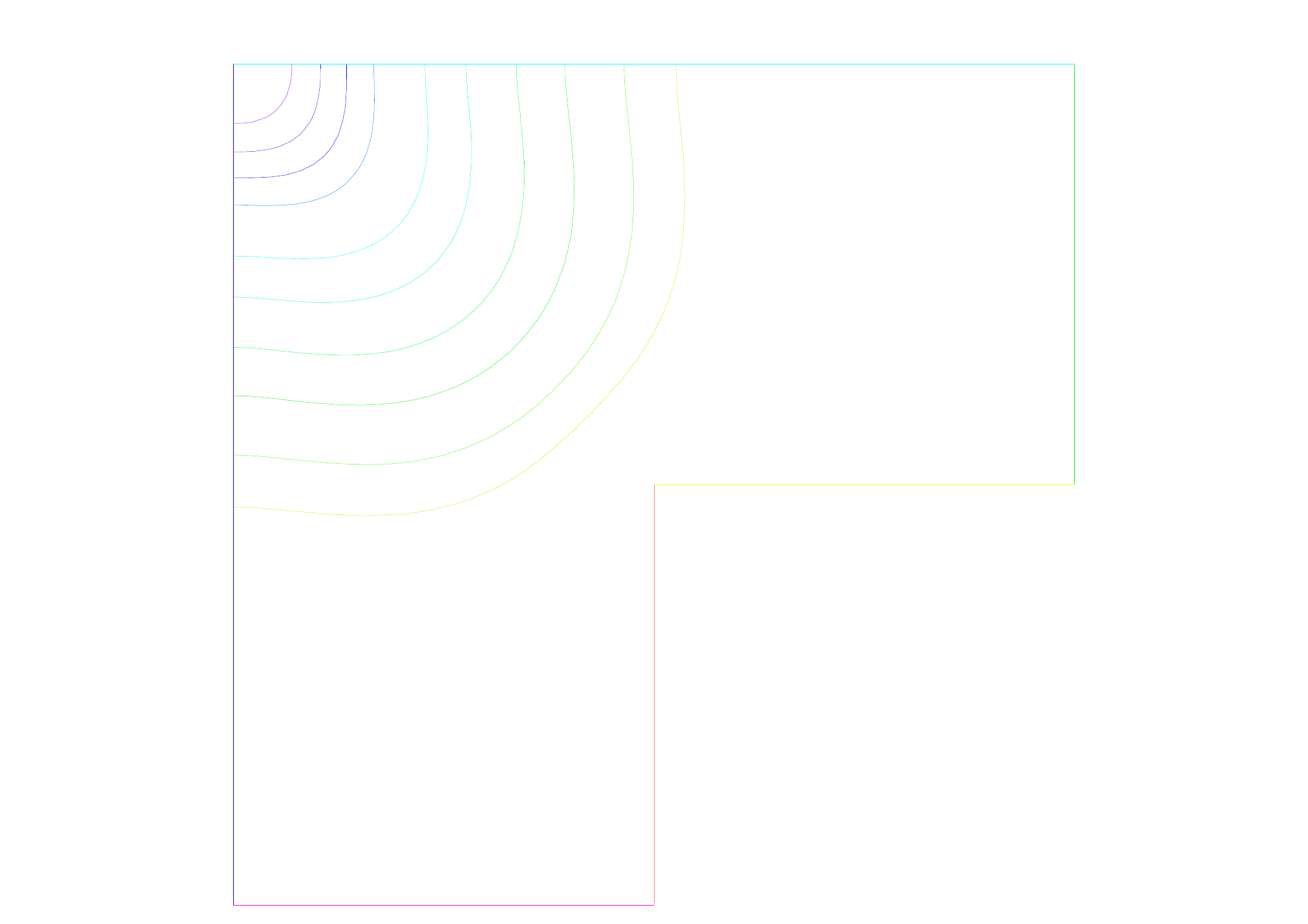}}
\subfigure[$t=40$] 
{\includegraphics[height=1.15in,width=1.6in]{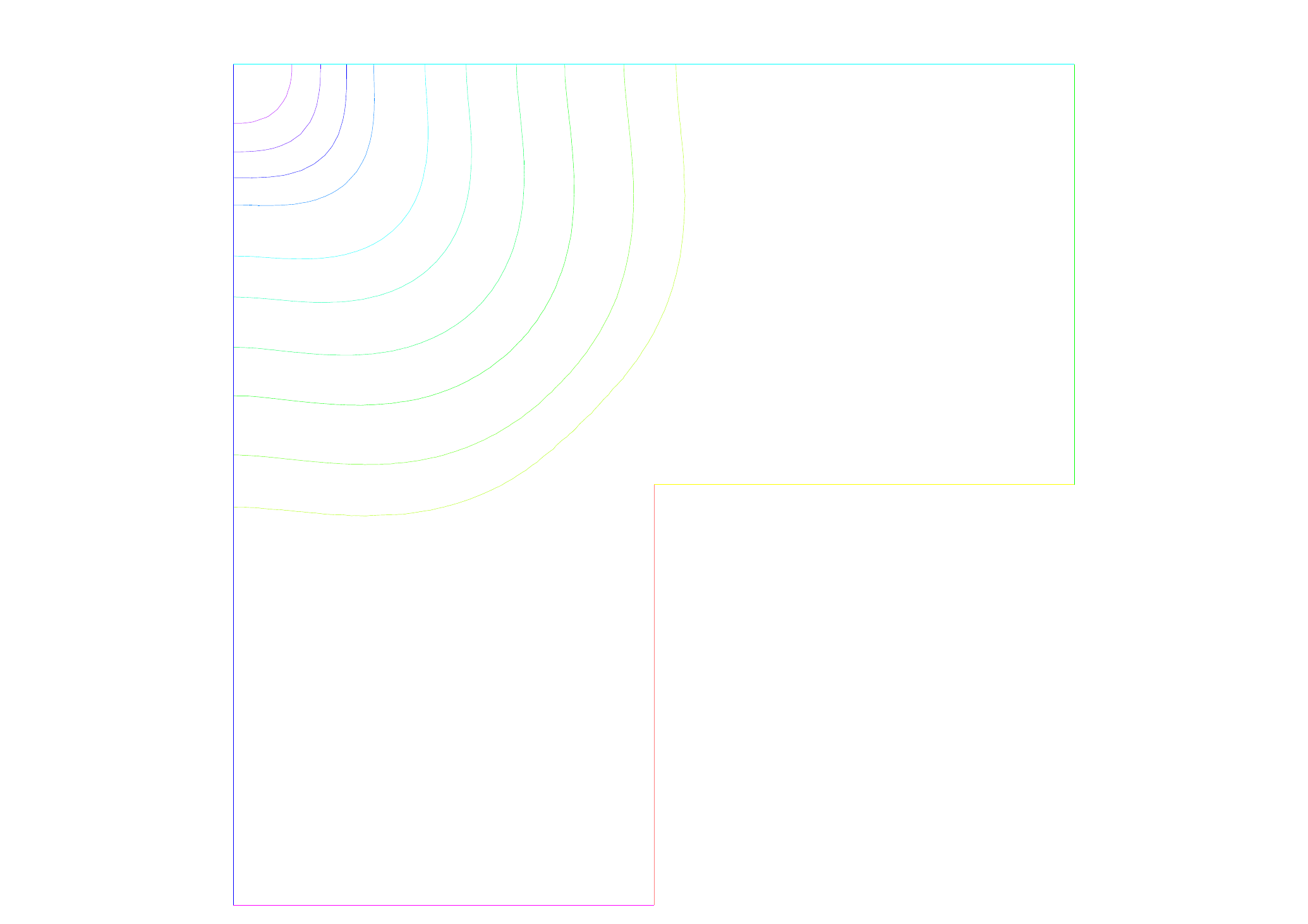}}
\caption{Contour of $|\psi|^2$ by solving 
the TDGL under the Lorentz gauge
with $h=1/64$.}
\label{LShp8}\vspace{10pt}
\subfigure[$t=5$] 
{\includegraphics[height=1.15in,width=1.6in]{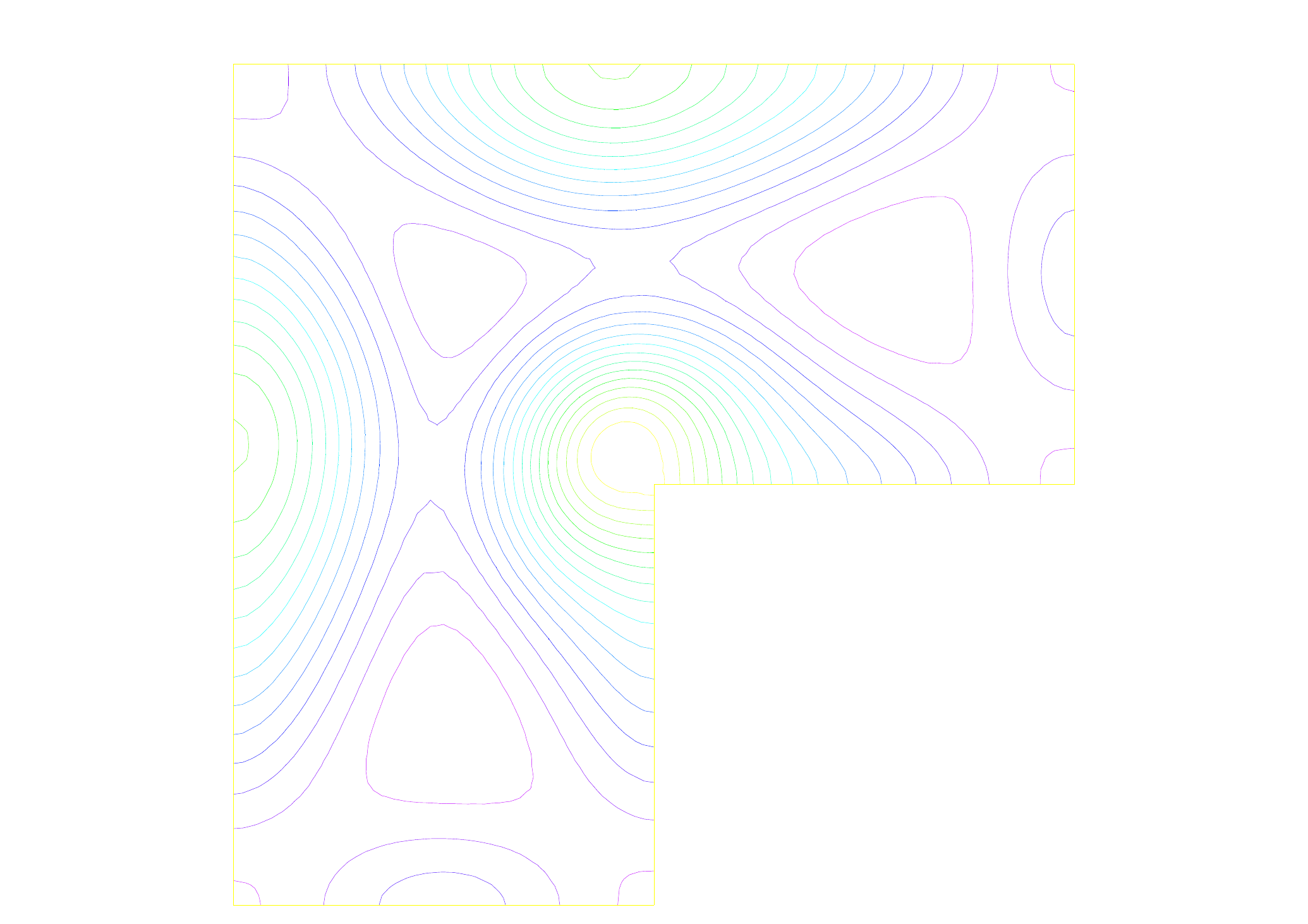}}
\subfigure[$t=20$] 
{\includegraphics[height=1.15in,width=1.6in]{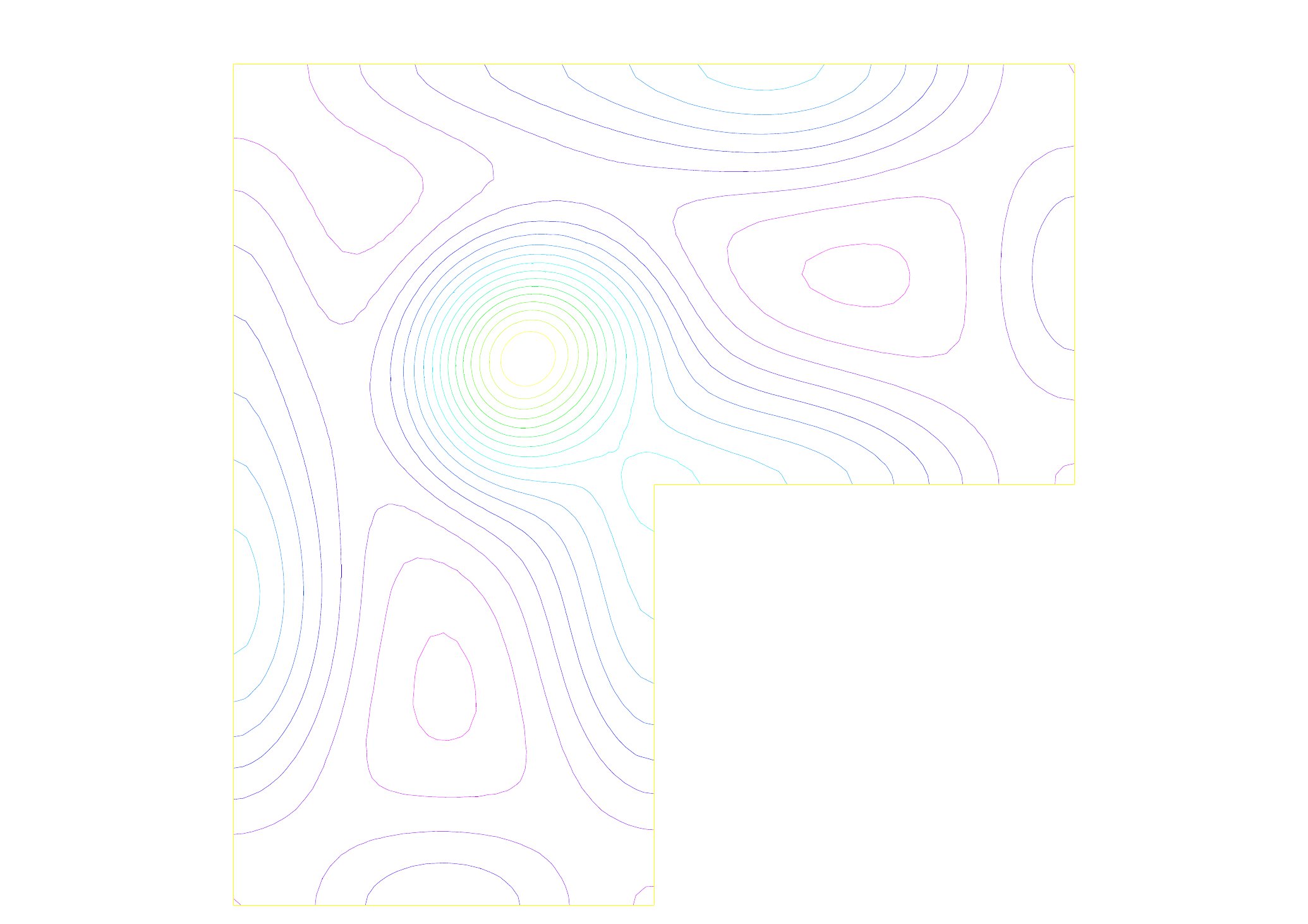}}
\subfigure[$t=40$] 
{\includegraphics[height=1.15in,width=1.6in]{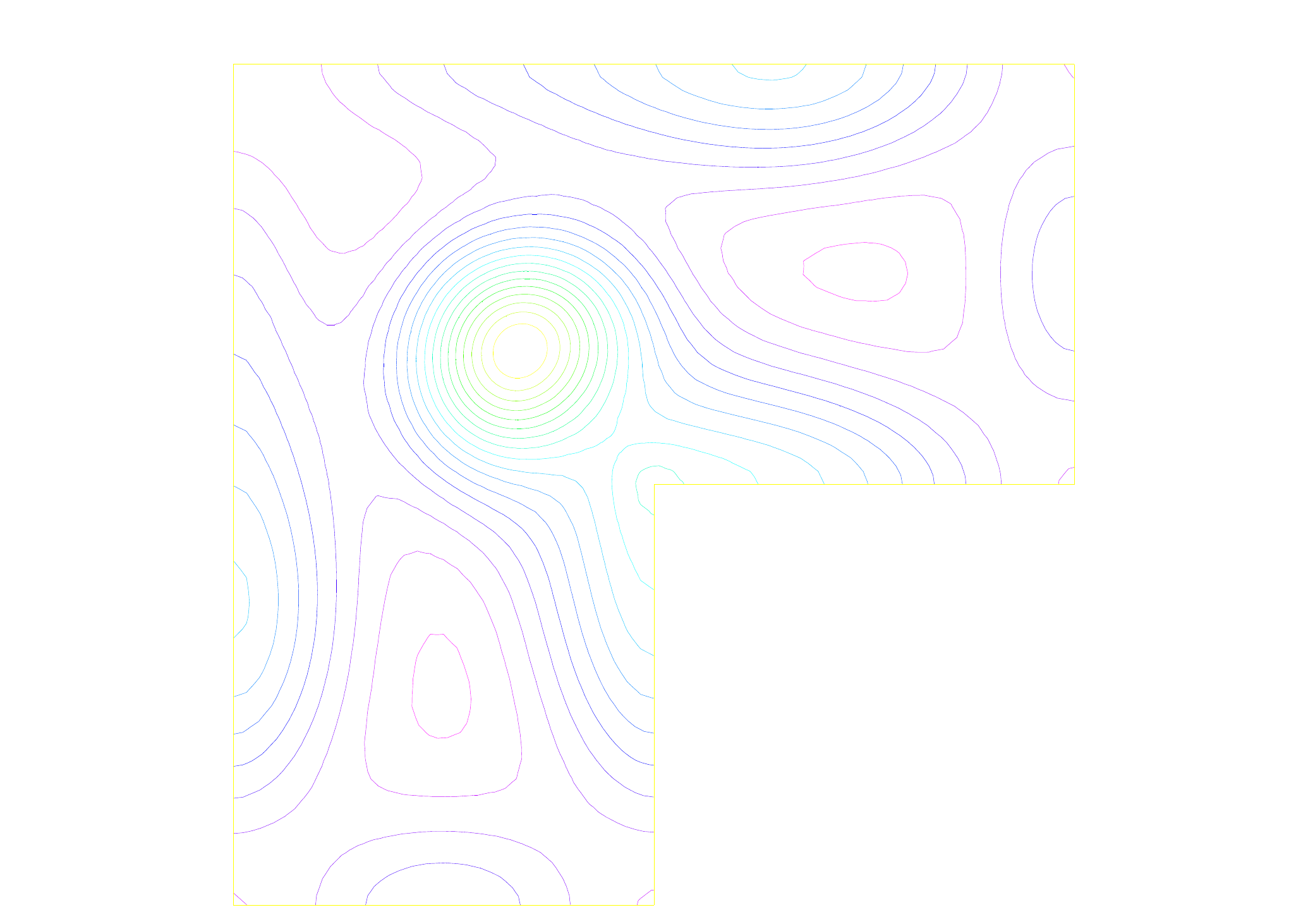}}
\caption{Contour of $|\psi|^2$ computed by the new approach
with $h=1/64$.}
\label{LShp9}\vspace{10pt}
\subfigure[quasi-uniform mesh] 
{\includegraphics[height=1.15in,width=1.6in]{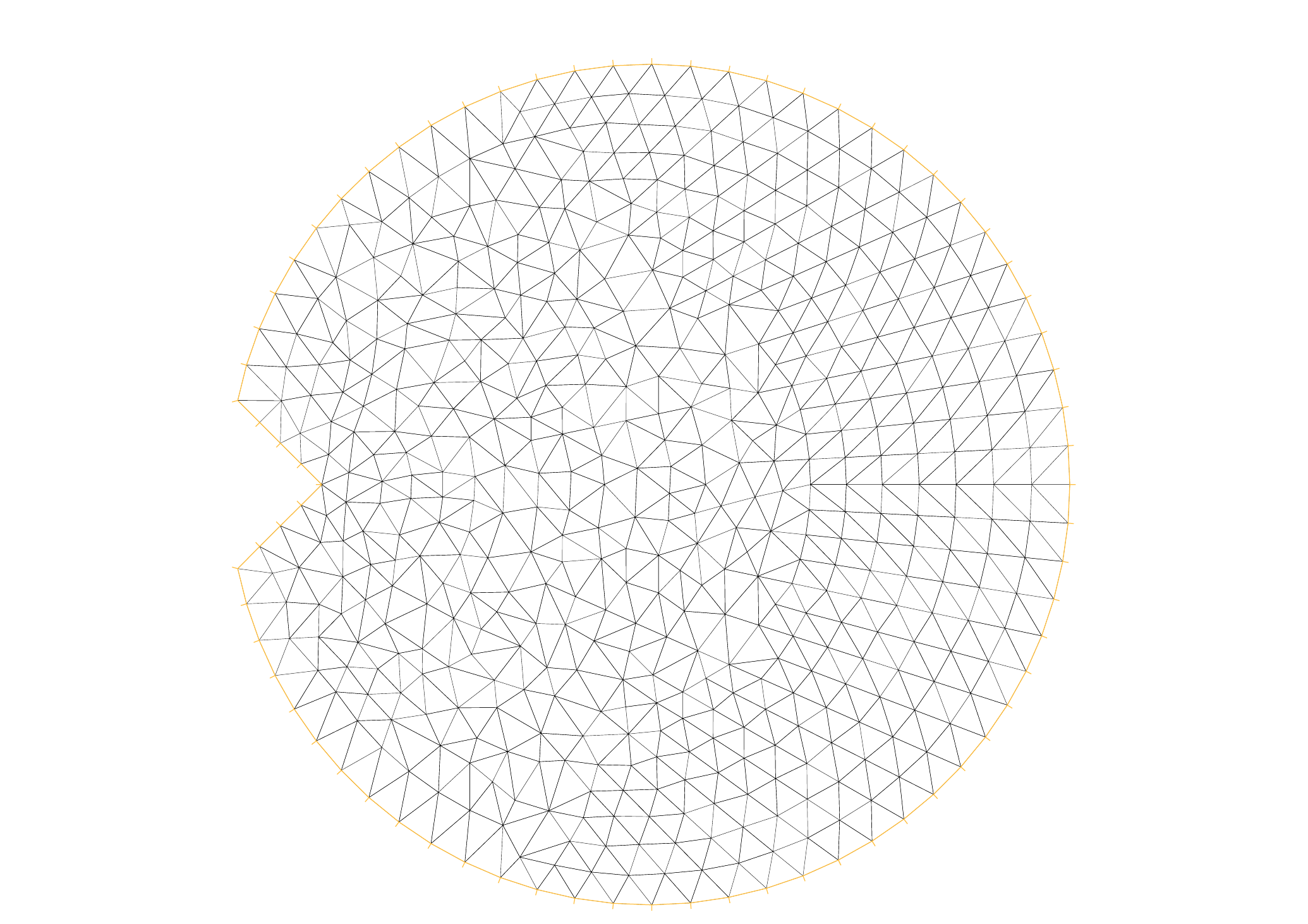}}
\subfigure[locally refined mesh] 
{\includegraphics[height=1.15in,width=1.6in]{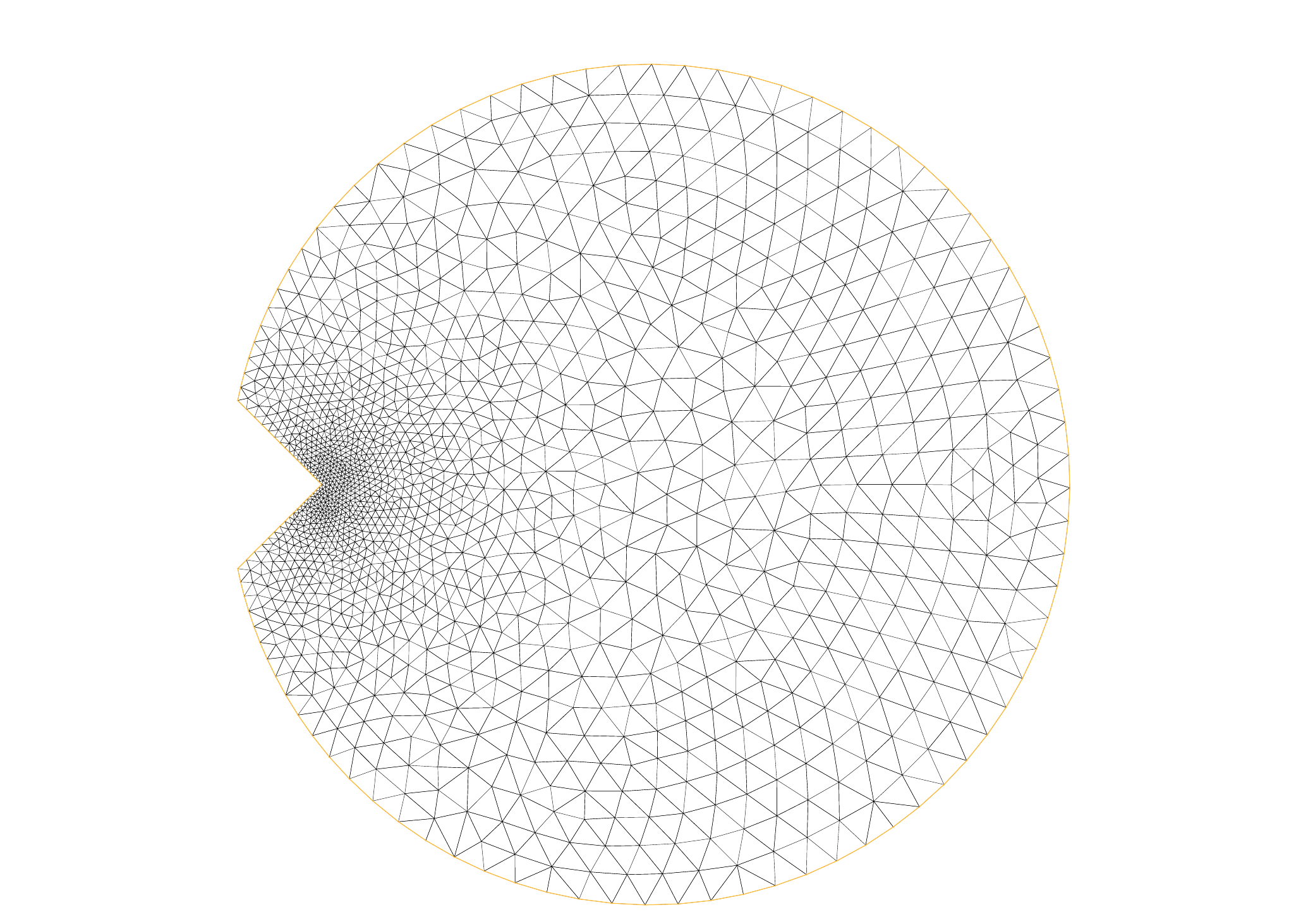}}
\caption{Triangulation of the circular domain with a 
triangular defect on the boundary.}
\label{CircularD}
\end{figure}

We solve 
\refe{RFPDE1}-\refe{RFinit} with the proposed numerical scheme
with piecewise linear finite elements and $\tau=0.1$, 
and compare the numerical results of the numerical solutions of 
\refe{TPDEq1}-\refe{TBDC} and \refe{PDE1}-\refe{TBDC2}, respectively,
by using the same finite element mesh and time-step size.
The contours of the numerical solutions
of $|\psi|^2$ are presented in Figure \ref{LShp1}--\ref{LShp9}
with $h=1/16$, $1/32$ and $1/64$.
One can see that the numerical solution of \refe{TPDEq1}-\refe{TBDC}
changes much as the mesh is refined from $h=1/16$ to $h=1/64$;
the numerical solution of \refe{PDE1}-\refe{TBDC2} and 
\refe{RFPDE1}-\refe{RFinit} are relatively stable as the mesh is refined.
When $h=1/16$ or $h=1/32$ the contours of the three numerical solutions
are very different, while when $h=1/64$ the 
the numerical solutions of \refe{TPDEq1}-\refe{TBDC}
and \refe{RFPDE1}-\refe{RFinit} agree.
Based on these numerical results, we 
see the following interesting phenomenons:

(1) the numerical solution of \refe{TPDEq1}-\refe{TBDC} 
is unstable with respect to the the mesh size;

(2) although the numerical solution of \refe{PDE1}-\refe{TBDC2} is stable
as the mesh refines, it converges to an incorrect solution;

(3) the numerical solution of \refe{RFPDE1}-\refe{RFinit} 
is stable and correct.

Athough the system \refe{RFPDE1}-\refe{RFinit} is derived from 
\refe{PDE1}-\refe{TBDC2} and the two systems are
equivalent theoretically, they are not equivalent computationally.
Clearly, the reformulated system can be solved easily by
the FEMs, while the original system requires extra work to
overcome its computational difficulty in a domain with reentrant corners.
\bigskip

{\it Example 3.2}$\quad$
Secondly, we present simulations 
of vortex dynamics of a type II superconductor in a 
circular disk with a triangular defect on the boundary.
This example has been tested in \cite{ASPM11} 
with 
$$
\eta=1,\quad \kappa=4,\quad
 \psi_0 = 1.0, \quad {\bf A}_0 = (0,0) ,
$$
and with several different values of $H$,
by solving the TDGL under the temporal gauge.
Details of the geometry of the domain can be found in the 
reference \cite{ASPM11}.
Figure \ref{CircularD} contains a quasi-uniform mesh and 
a locally refined mesh on the circular domain with 64
points on the boundary. Our computations below
use similar triangulations with 256 points on the boundary.

\begin{figure}[htp]
\centering
\subfigure[$t=20$] 
{\includegraphics[height=1.15in,width=1.6in]{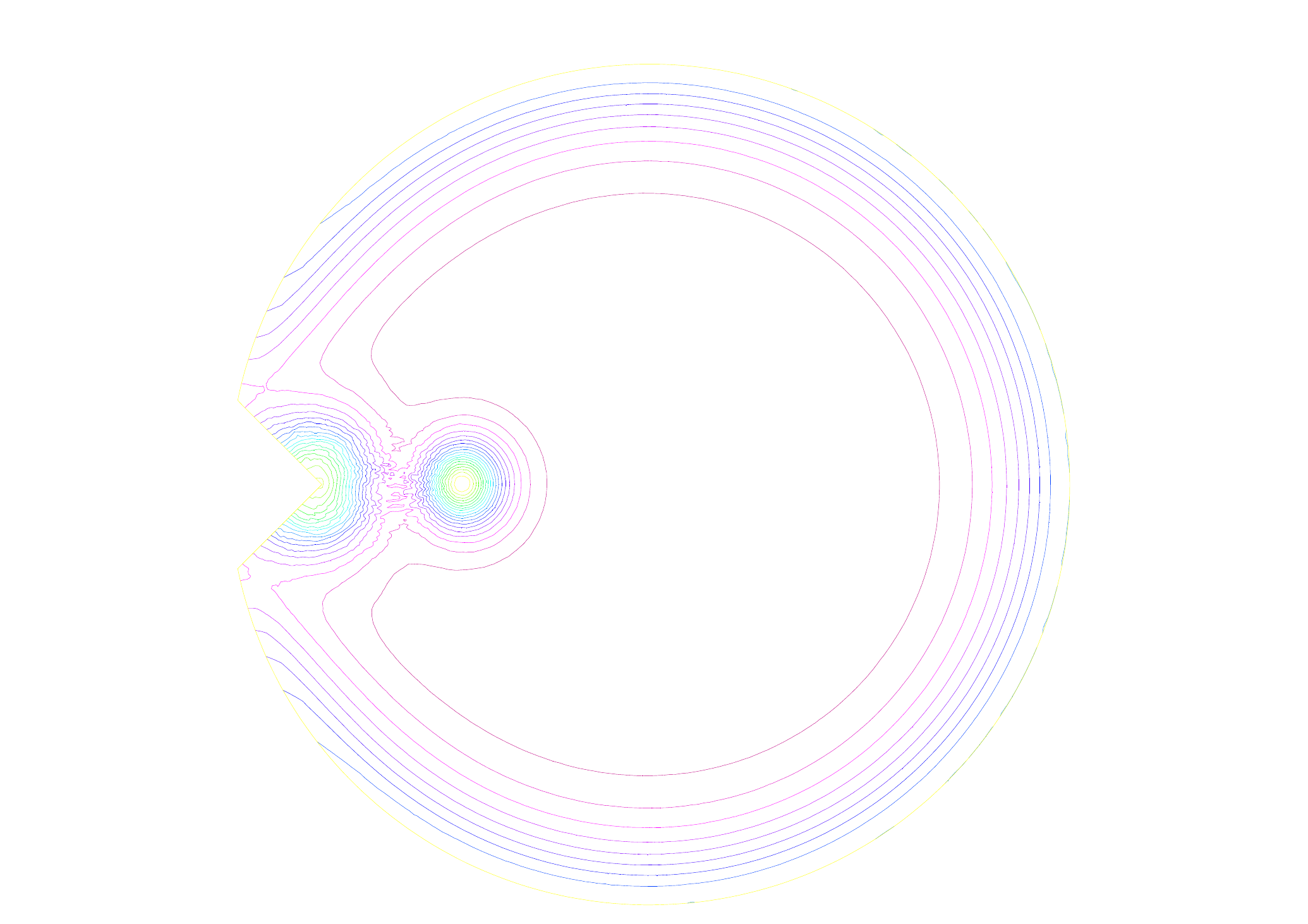}}
\subfigure[$t=100$] 
{\includegraphics[height=1.15in,width=1.6in]{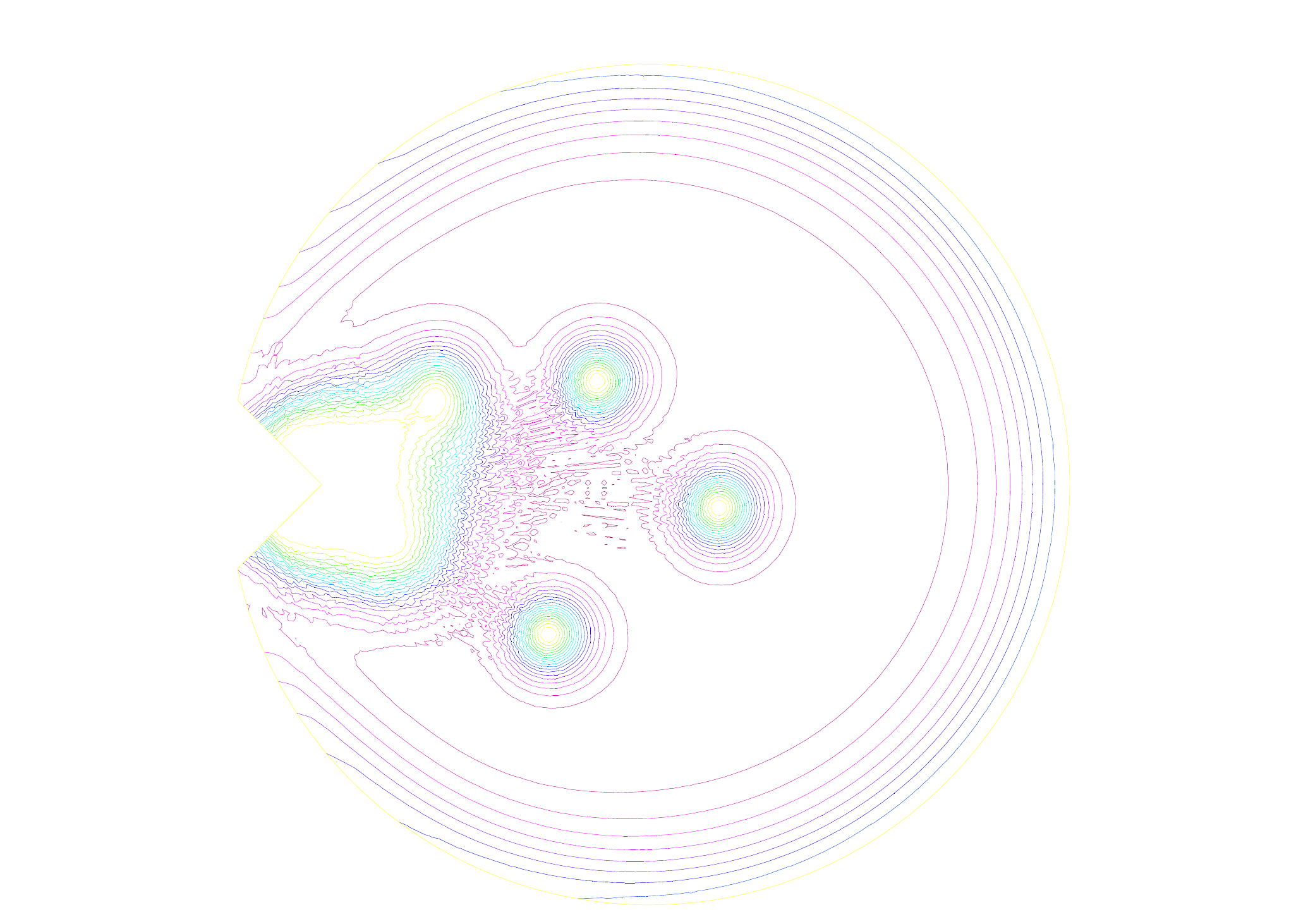}}
\subfigure[$t=15000$]
{\includegraphics[height=1.15in,width=1.6in]{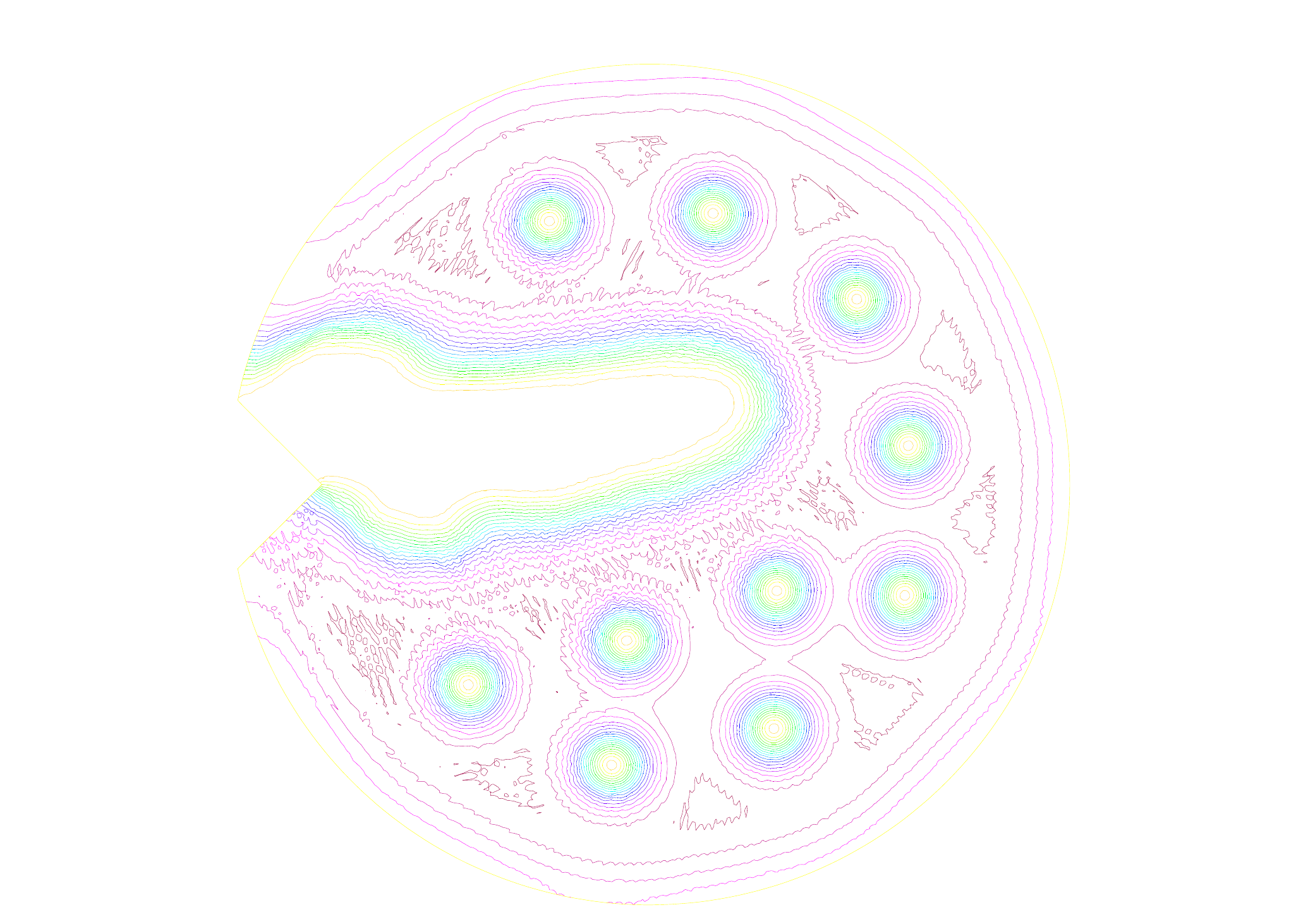}}
\caption{Contour of $|\psi|^2$ with $H=0.8$ by solving 
the TDGL under the temporal gauge.}
\label{FigS1}\vspace{10pt}
\subfigure[$t=20$] 
{\includegraphics[height=1.15in,width=1.6in]{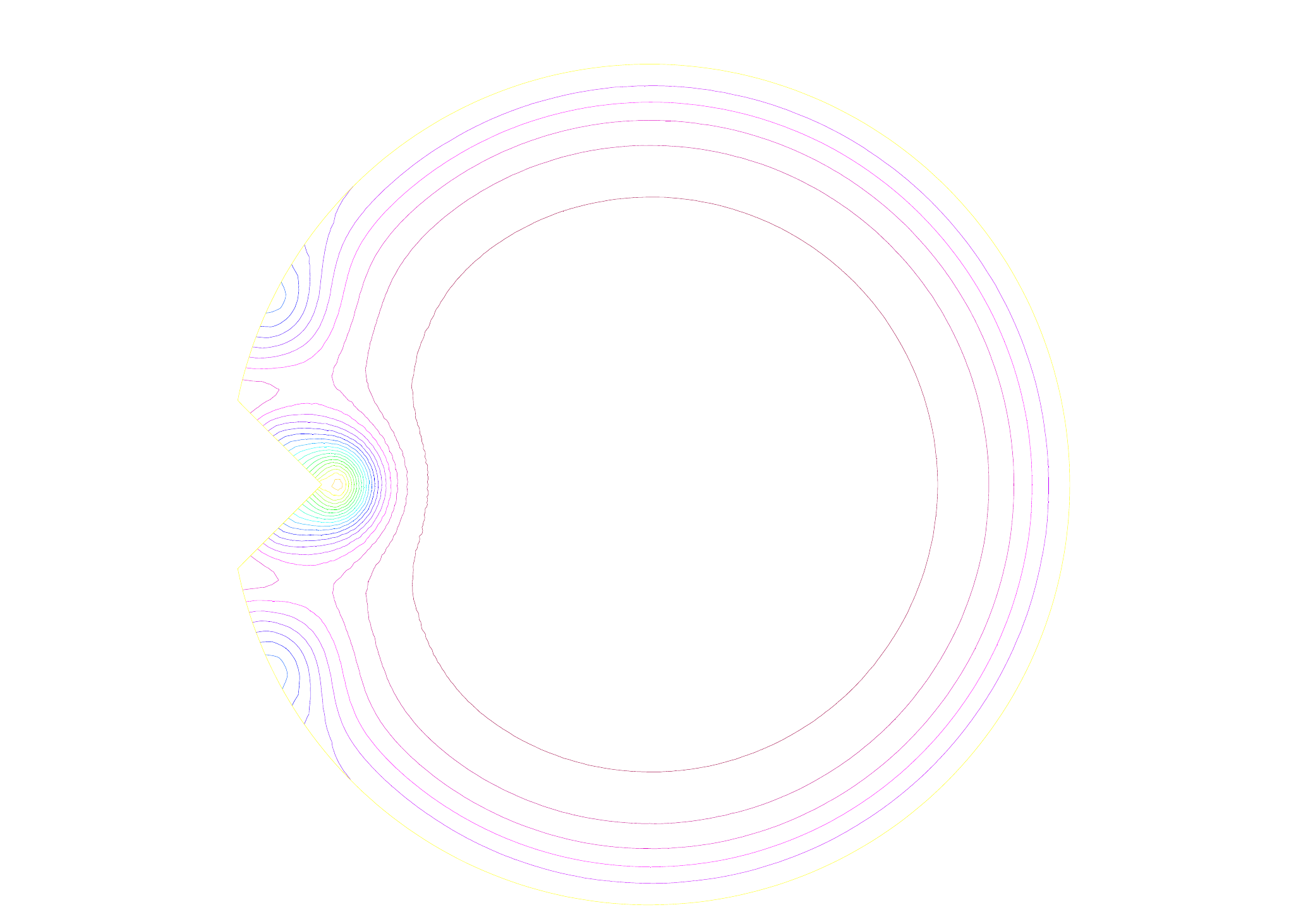}}
\subfigure[$t=100$] 
{\includegraphics[height=1.15in,width=1.6in]{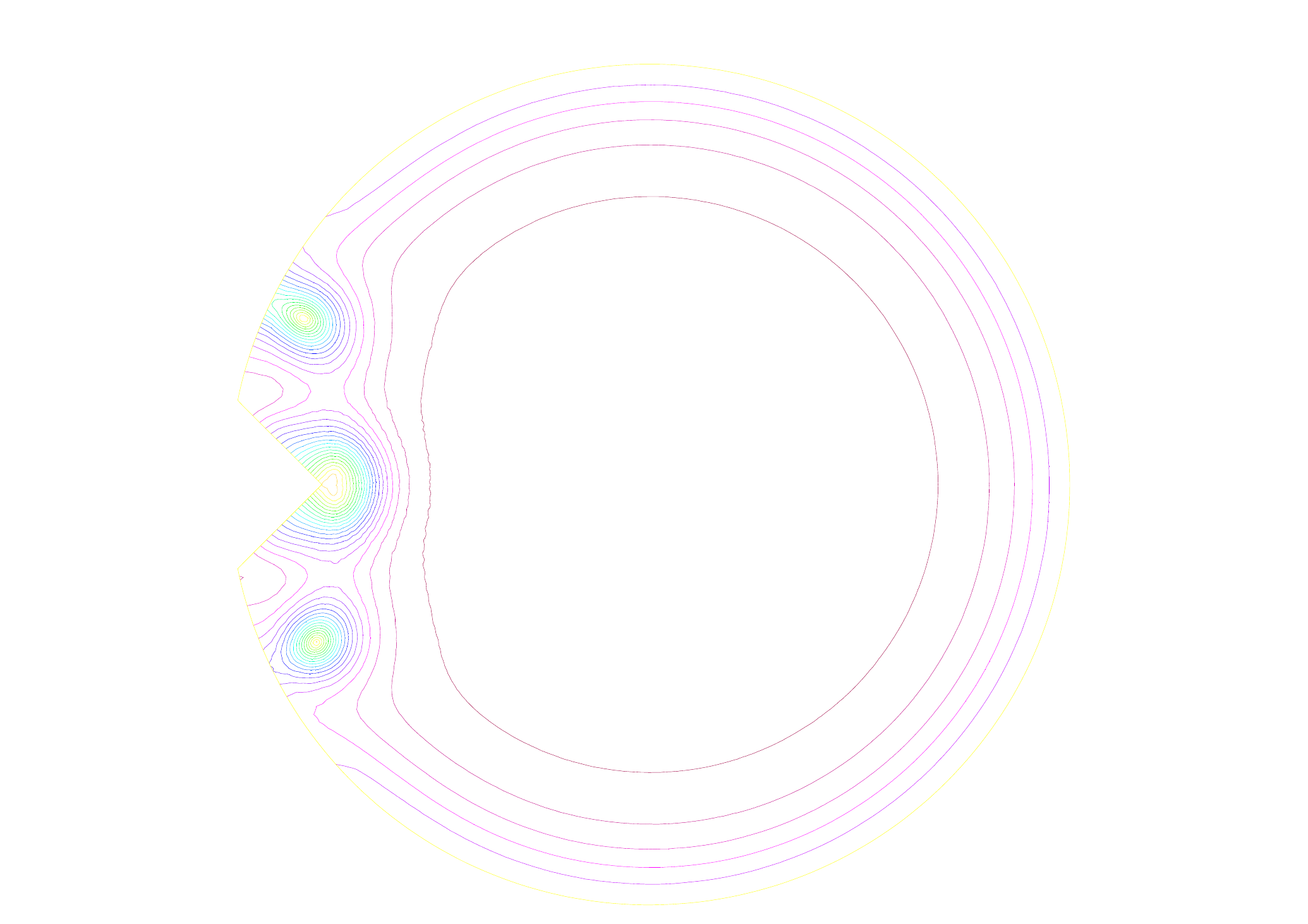}}
\subfigure[$t=15000$] 
{\includegraphics[height=1.15in,width=1.6in]{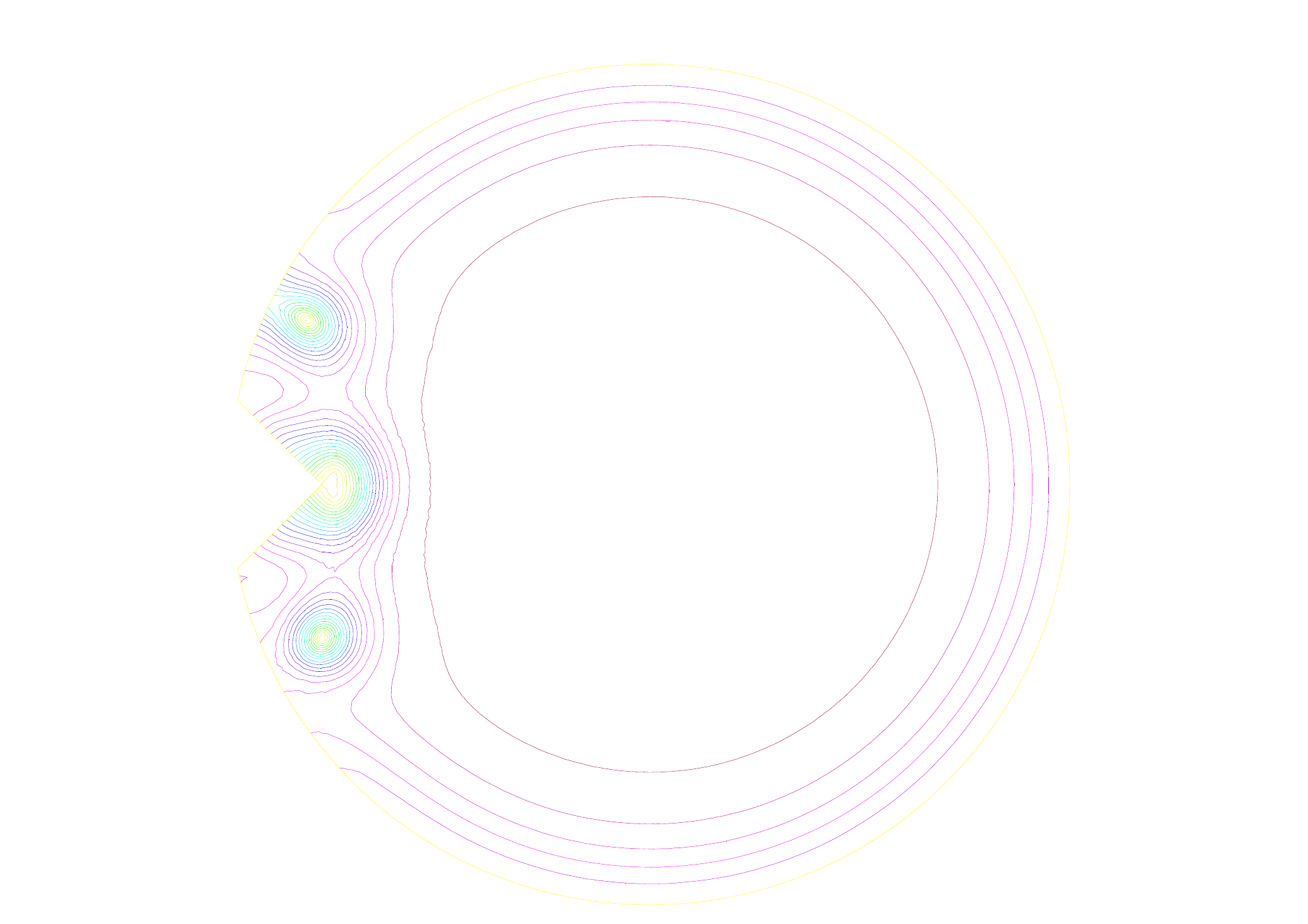}}
\caption{Contour of $|\psi|^2$ with $H=0.8$ by solving 
the TDGL under the Lorentz gauge.}
\label{FigS2}\vspace{10pt}
\subfigure[$t=20$] 
{\includegraphics[height=1.15in,width=1.6in]{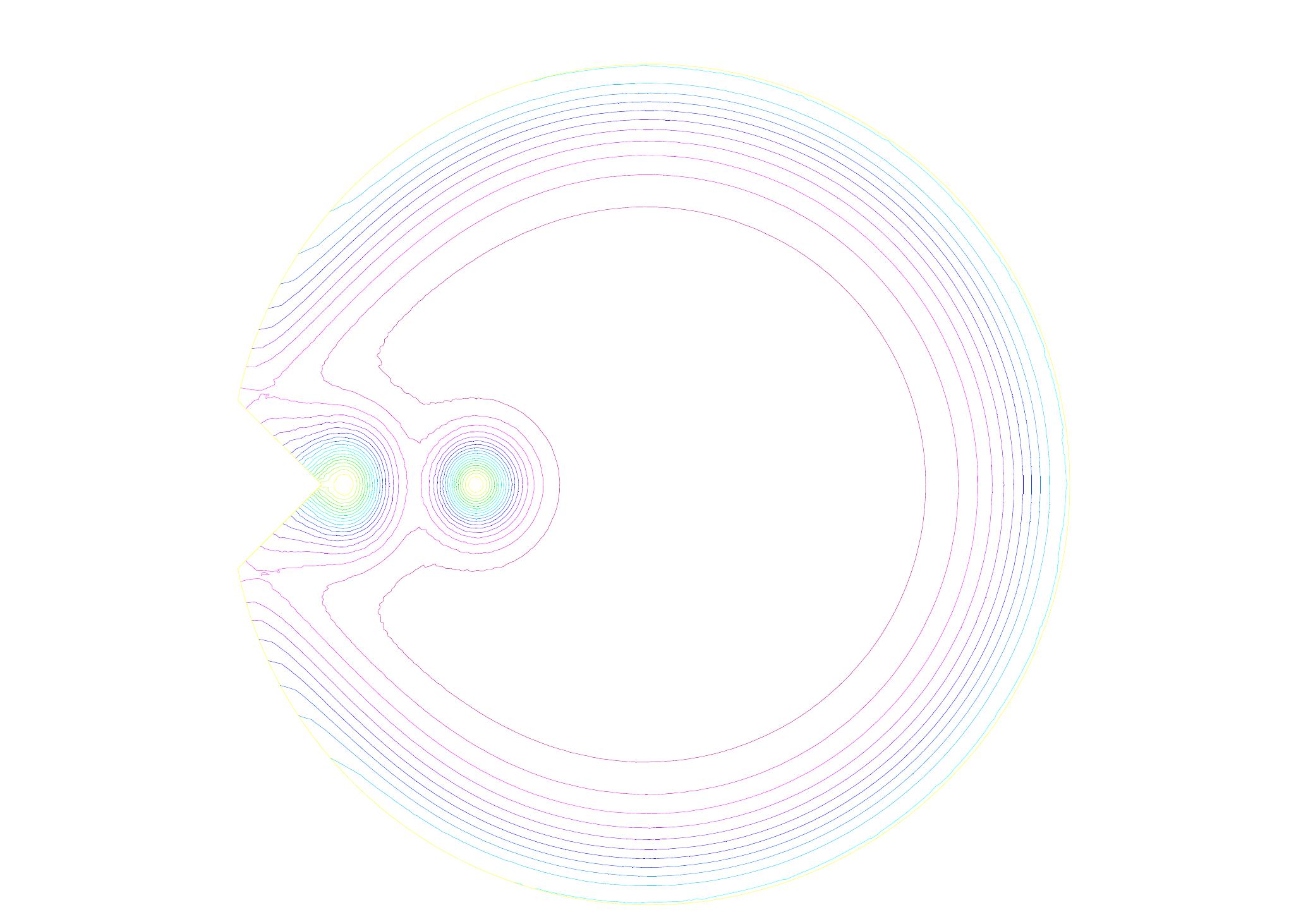}}
\subfigure[$t=100$] 
{\includegraphics[height=1.15in,width=1.6in]{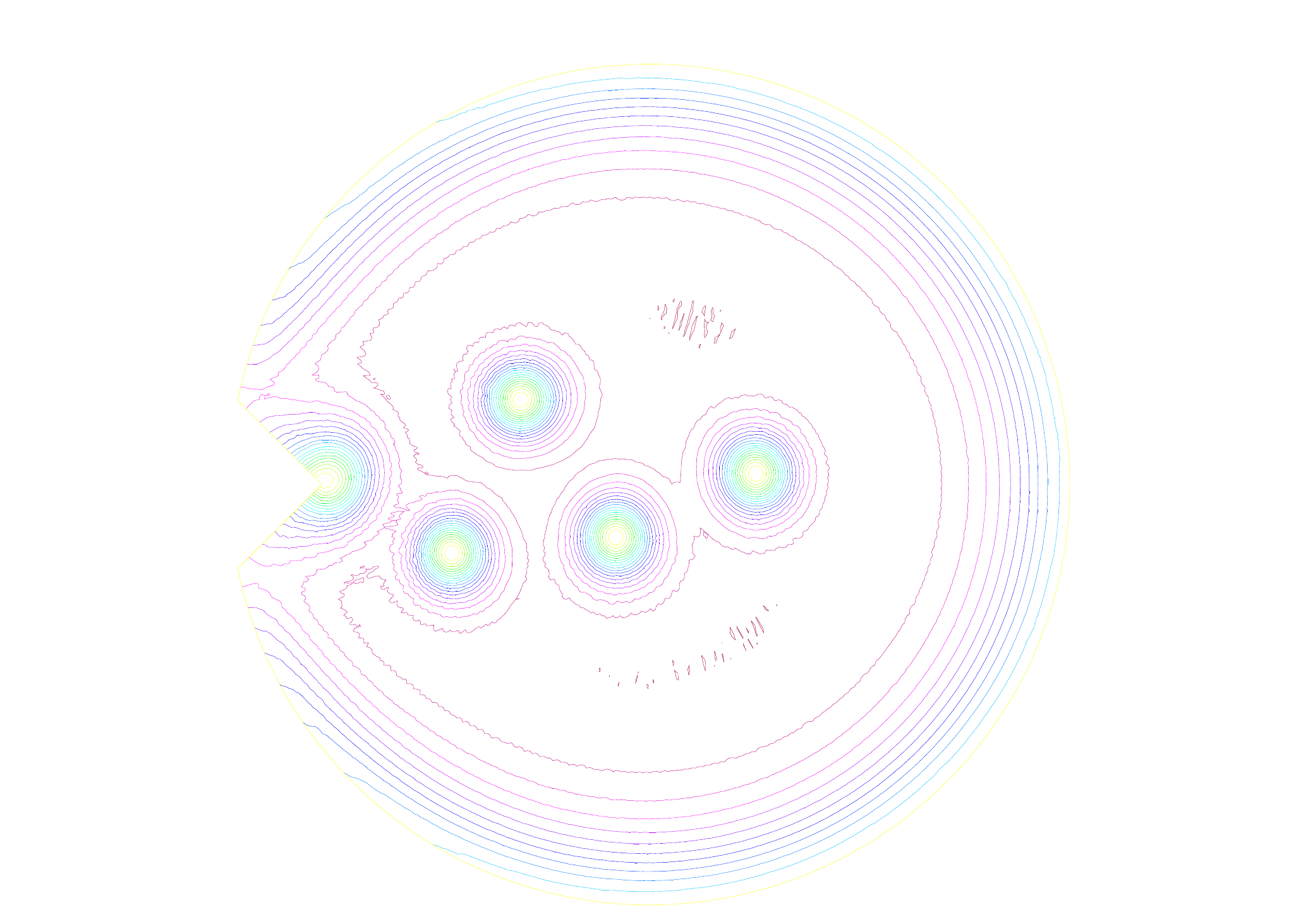}}
\subfigure[$t=15000$] 
{\includegraphics[height=1.15in,width=1.6in]{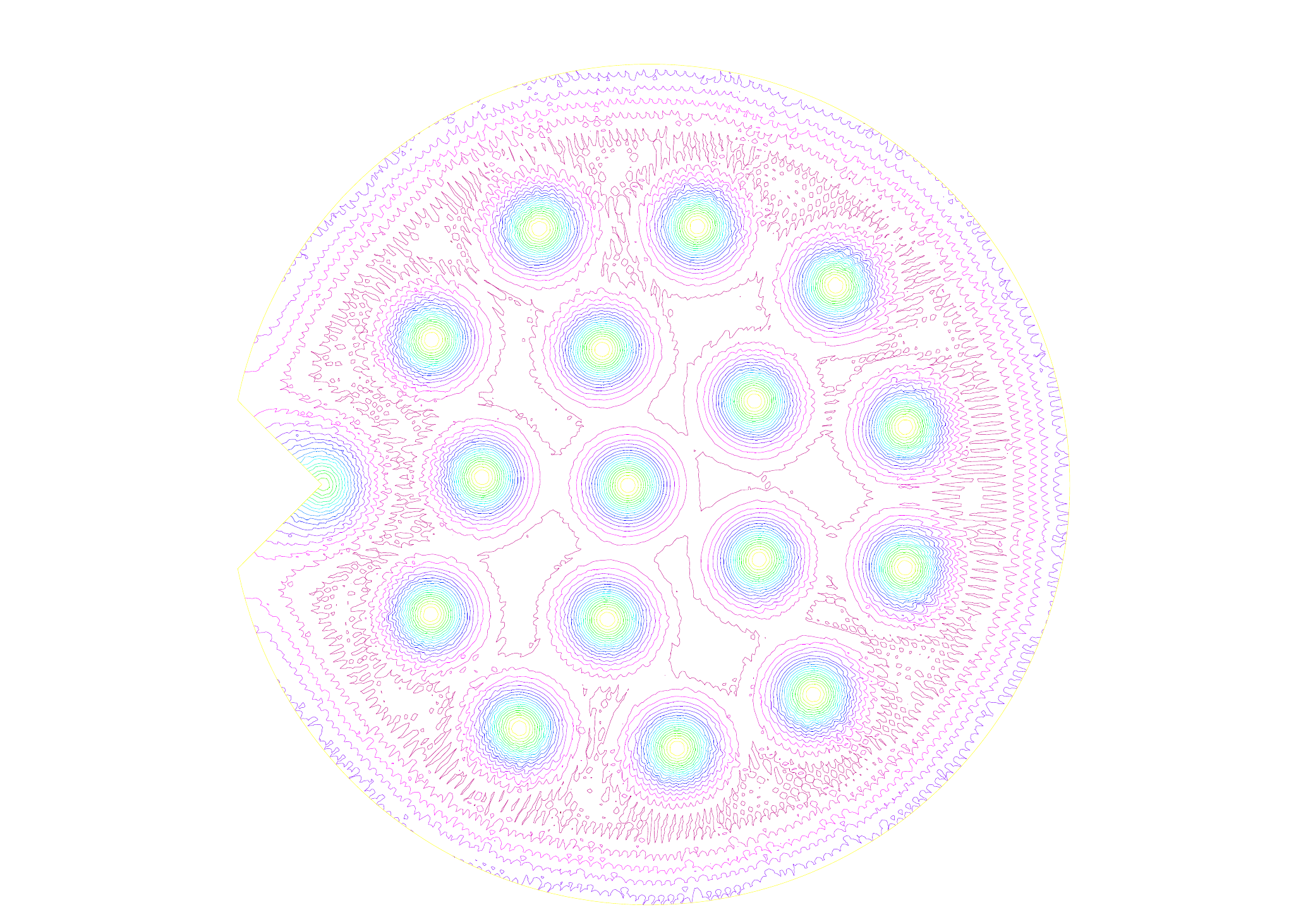}}
\caption{Contour of $|\psi|^2$ with $H=0.8$ computed by the new approach.}
\label{FigS3}\vspace{10pt}
\subfigure[$t=20$] 
{\includegraphics[height=1.15in,width=1.6in]{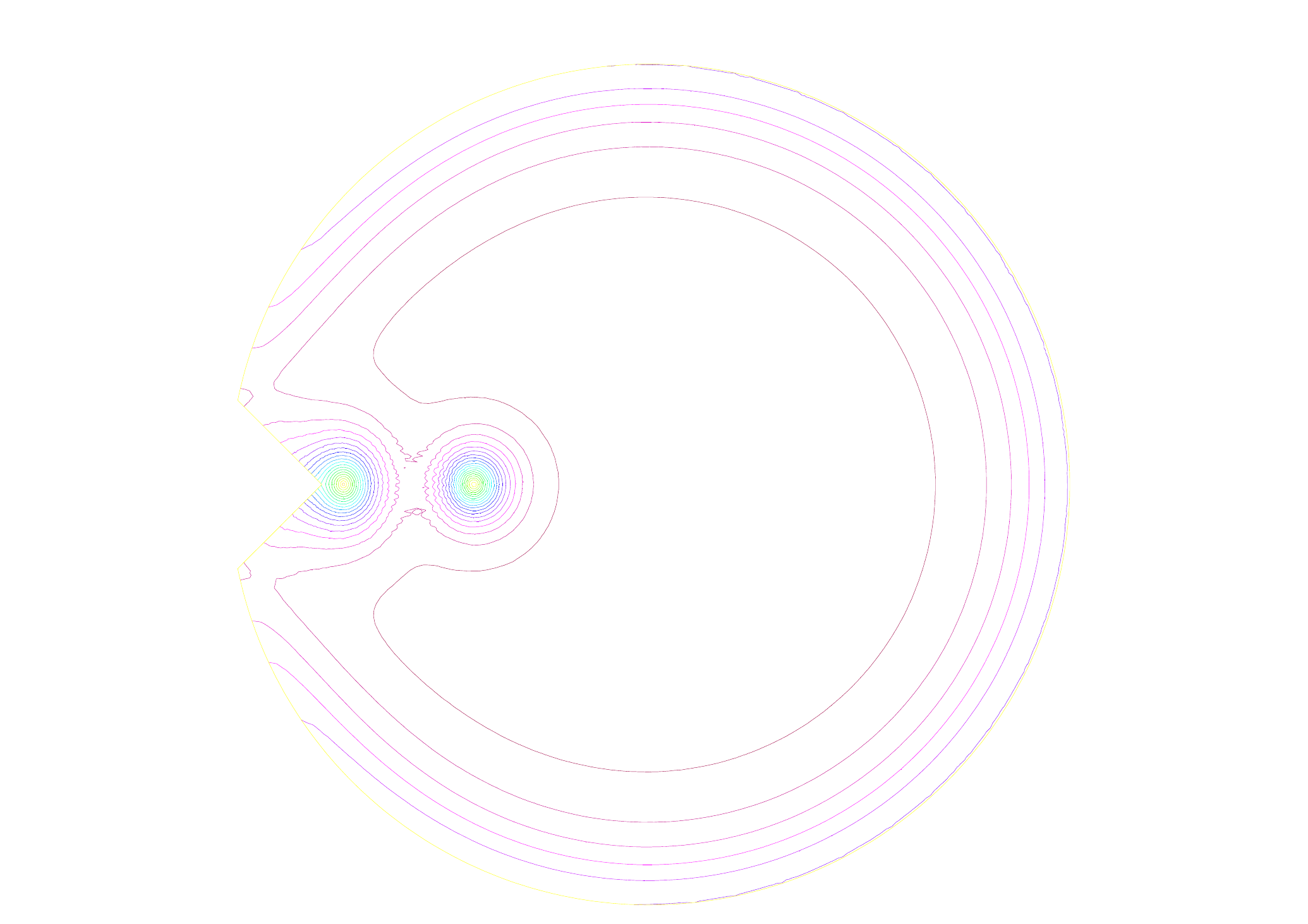}}
\subfigure[$t=100$] 
{\includegraphics[height=1.15in,width=1.6in]{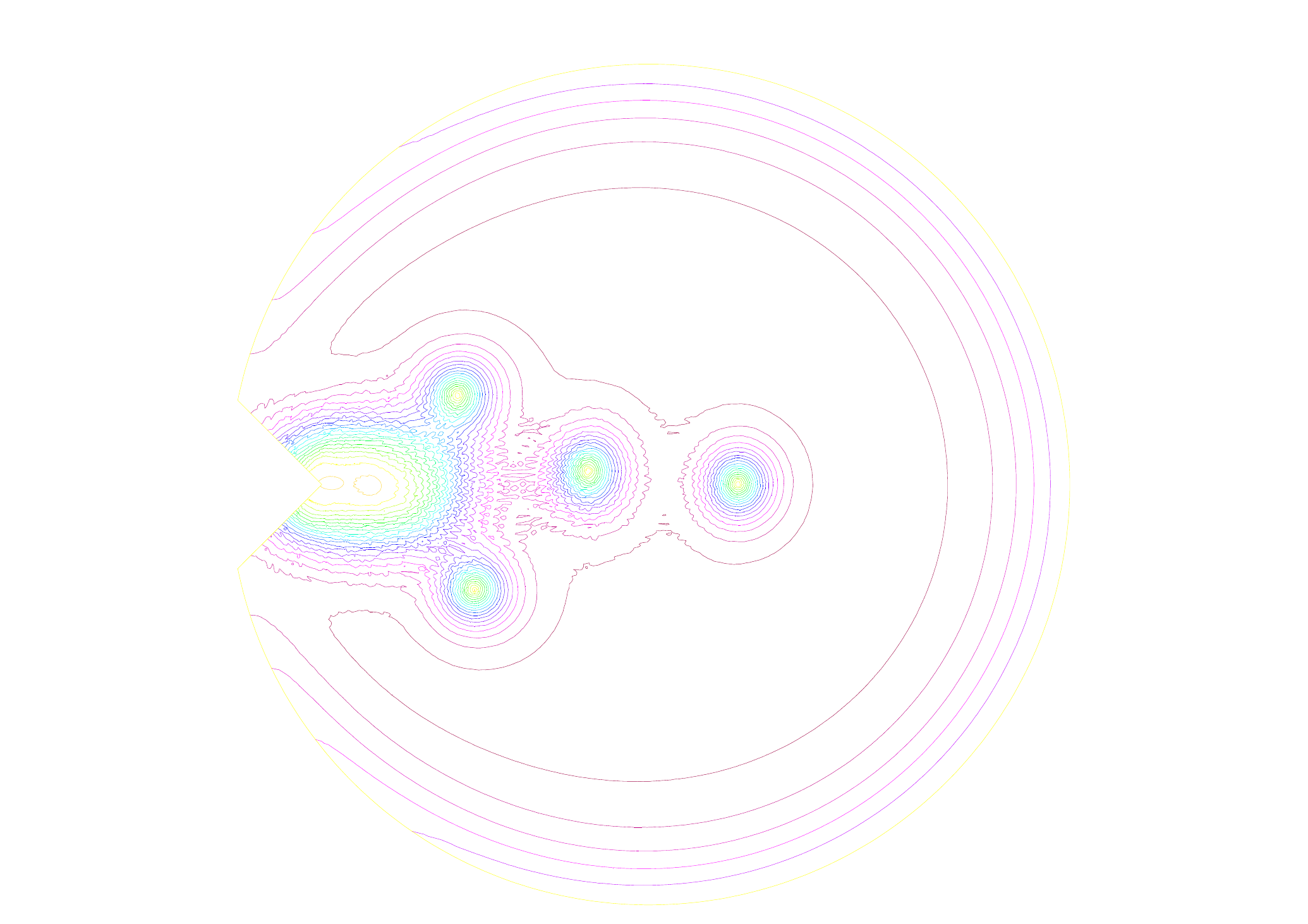}}
\subfigure[$t=15000$] 
{\includegraphics[height=1.15in,width=1.6in]{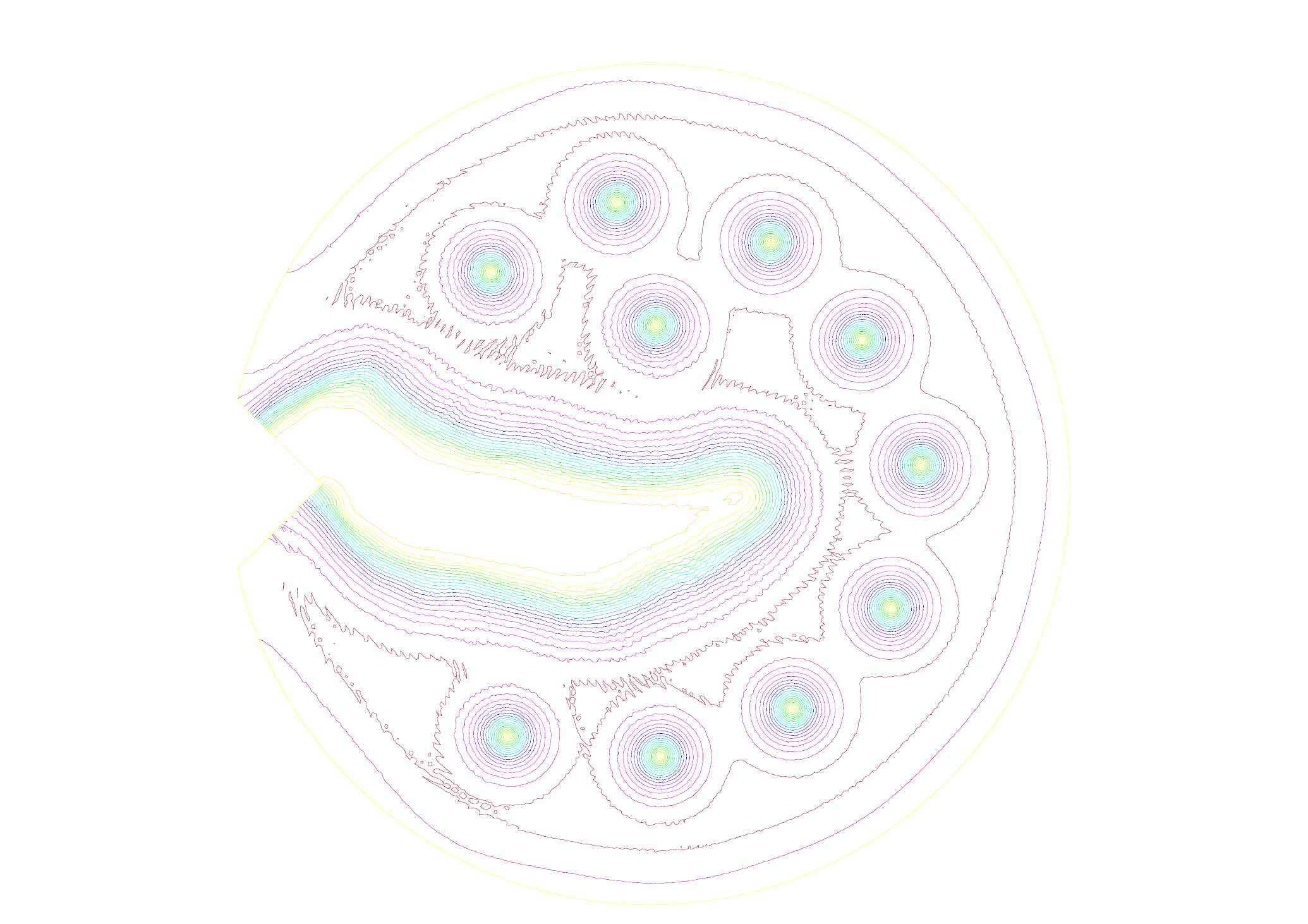}}
\caption{Contour of $|\psi|^2$ with $H=0.8$ by 
solving the TDGL under the temporal gauge 
with a locally refined mesh.}
\label{FigS4}
\end{figure}

For $H=0.8$, we solve \refe{TPDEq1}-\refe{TBDC}, 
\refe{PDE1}-\refe{TBDC2} and \refe{RFPDE1}-\refe{RFinit}, respectively,
with $\tau=0.1$ and piecewise linear finite elements subject to a
quasi-uniform triangulation of the domain with 256 points on the boundary;
Figure \ref{CircularD}-(a) contains such a triangulation 
with only 64 points on the boundary.
The contours of $|\psi|^2$ are presented in Figure \ref{FigS1}--\ref{FigS3}.
From Figure \ref{FigS1} one can see that a vortex at the concave corner
grows larger and larger as time grows, penetrating into the superconductor,
while this giant vortex is not reflected in Figure \ref{FigS3},
which we believe is a correct approximation of the exact solution. 
Excluding the giant vortex at the corner, the rest part of Figure \ref{FigS1}
looks similar as Figure \ref{FigS3} when $t$ is very large.
We believe that the giant vortex is a numerical pollution, whose 
shape will change if mesh changes (see Figure \ref{FigS4}). 
This indicates that the numerical solution of \refe{TPDEq1}-\refe{TBDC} 
is unstable compared with the numerical solution of \refe{RFPDE1}-\refe{RFinit}.
Clearly, Figure \ref{FigS2} is different from both Figure 
\ref{FigS1} and Figure \ref{FigS3},
and this implies that the finite element solution 
of \refe{PDE1}-\refe{TBDC2} may converge to an incorrect solution.

\begin{figure}[htp]
\centering
\subfigure[$t=25$] 
{\includegraphics[height=1.15in,width=1.6in]{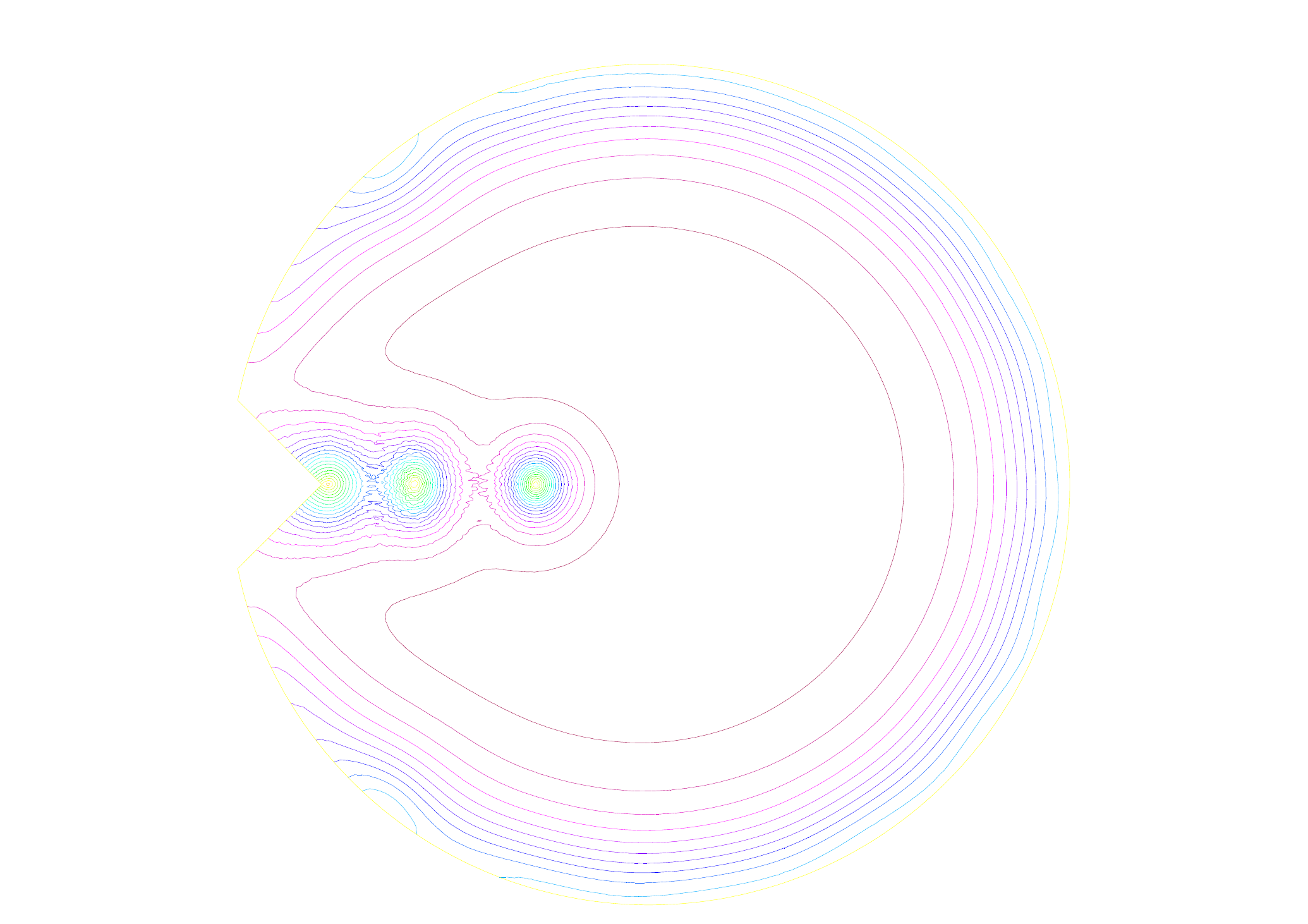}}
\subfigure[$t=30$] 
{\includegraphics[height=1.15in,width=1.6in]{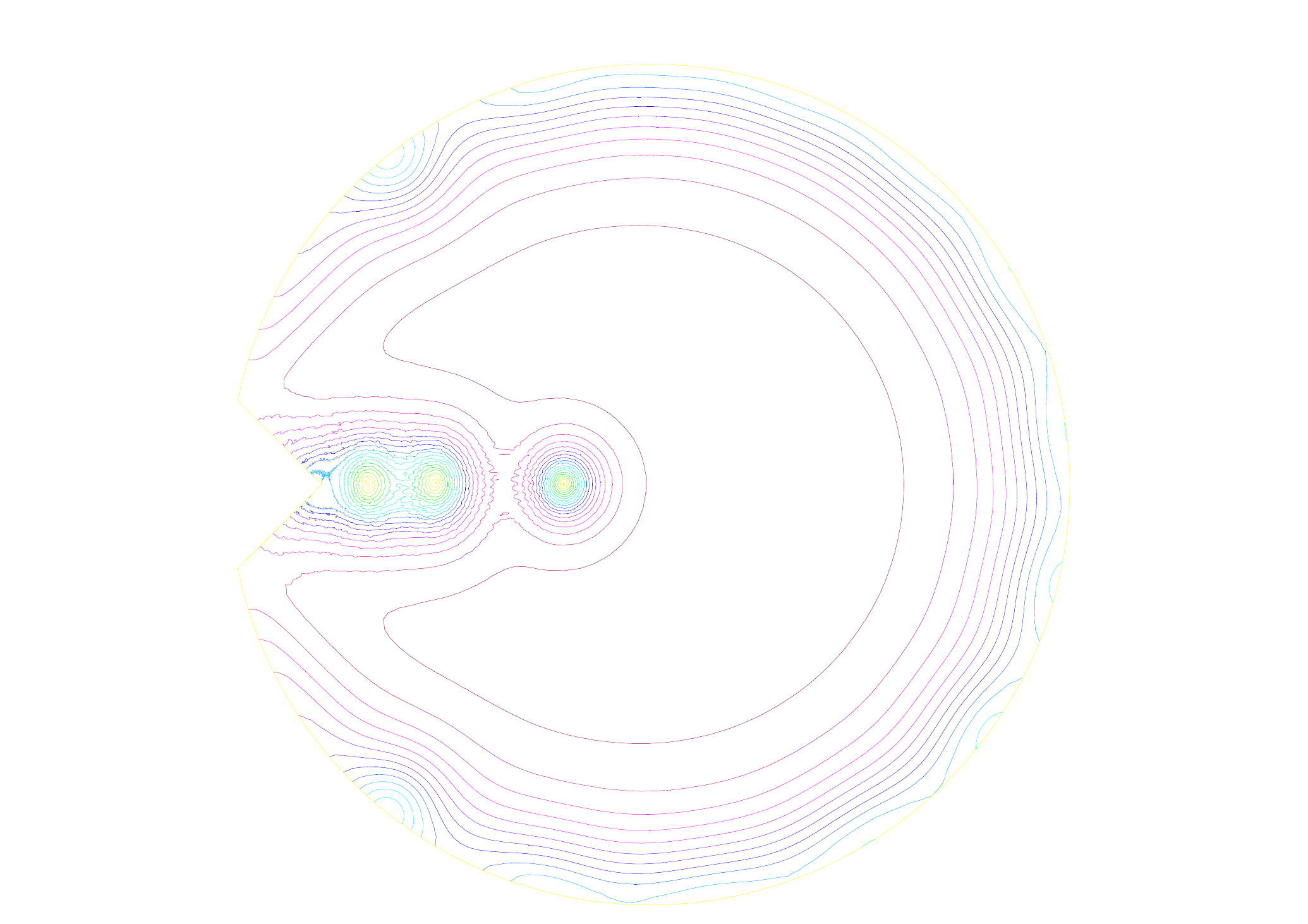}}
\subfigure[$t=5000$]
{\includegraphics[height=1.15in,width=1.6in]{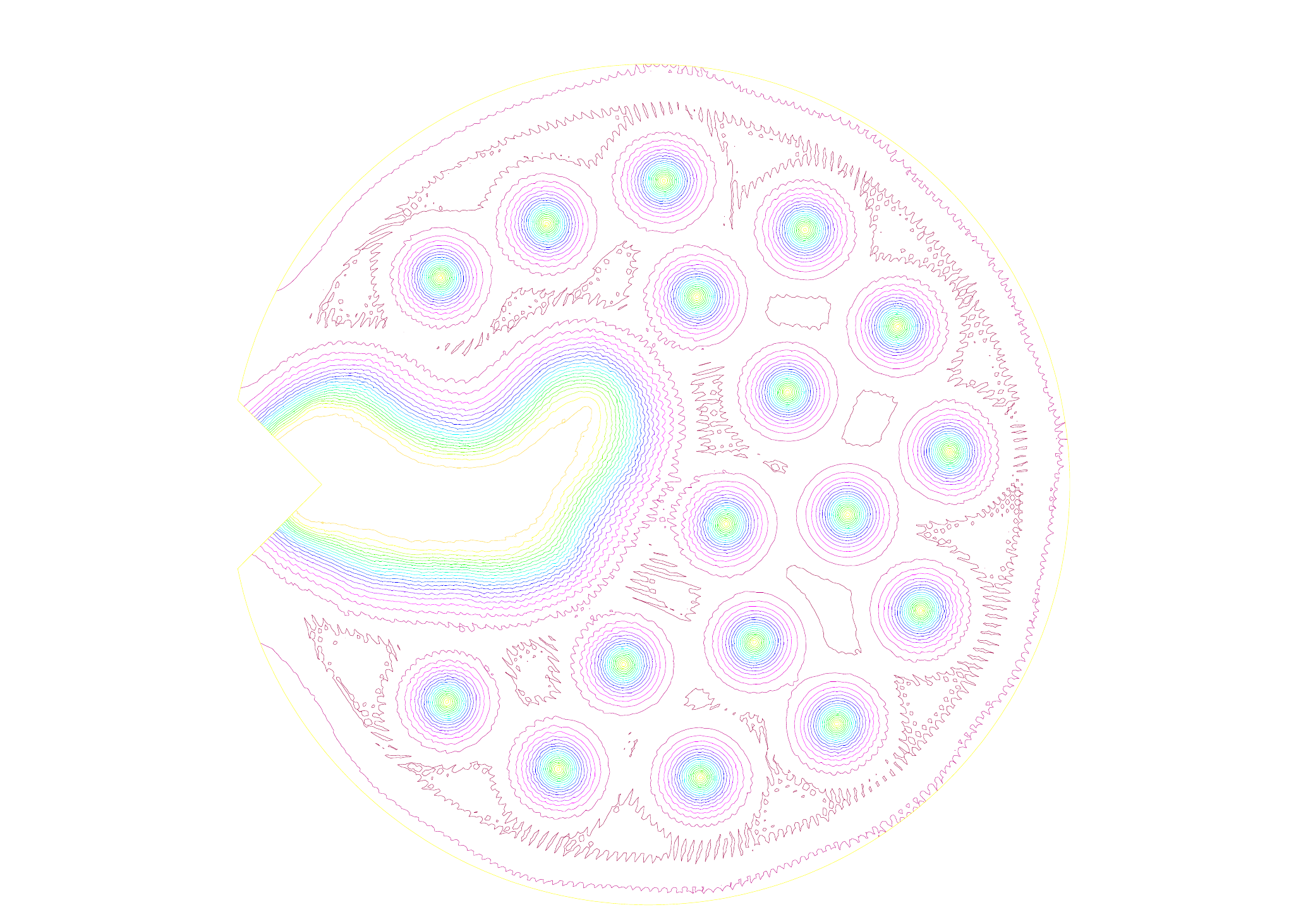}}
\caption{Contour of $|\psi|^2$ with $H=0.9$ by solving 
the TDGL under the temporal gauge.}
\label{FigS5}\vspace{10pt}
\subfigure[$t=25$] 
{\includegraphics[height=1.15in,width=1.6in]{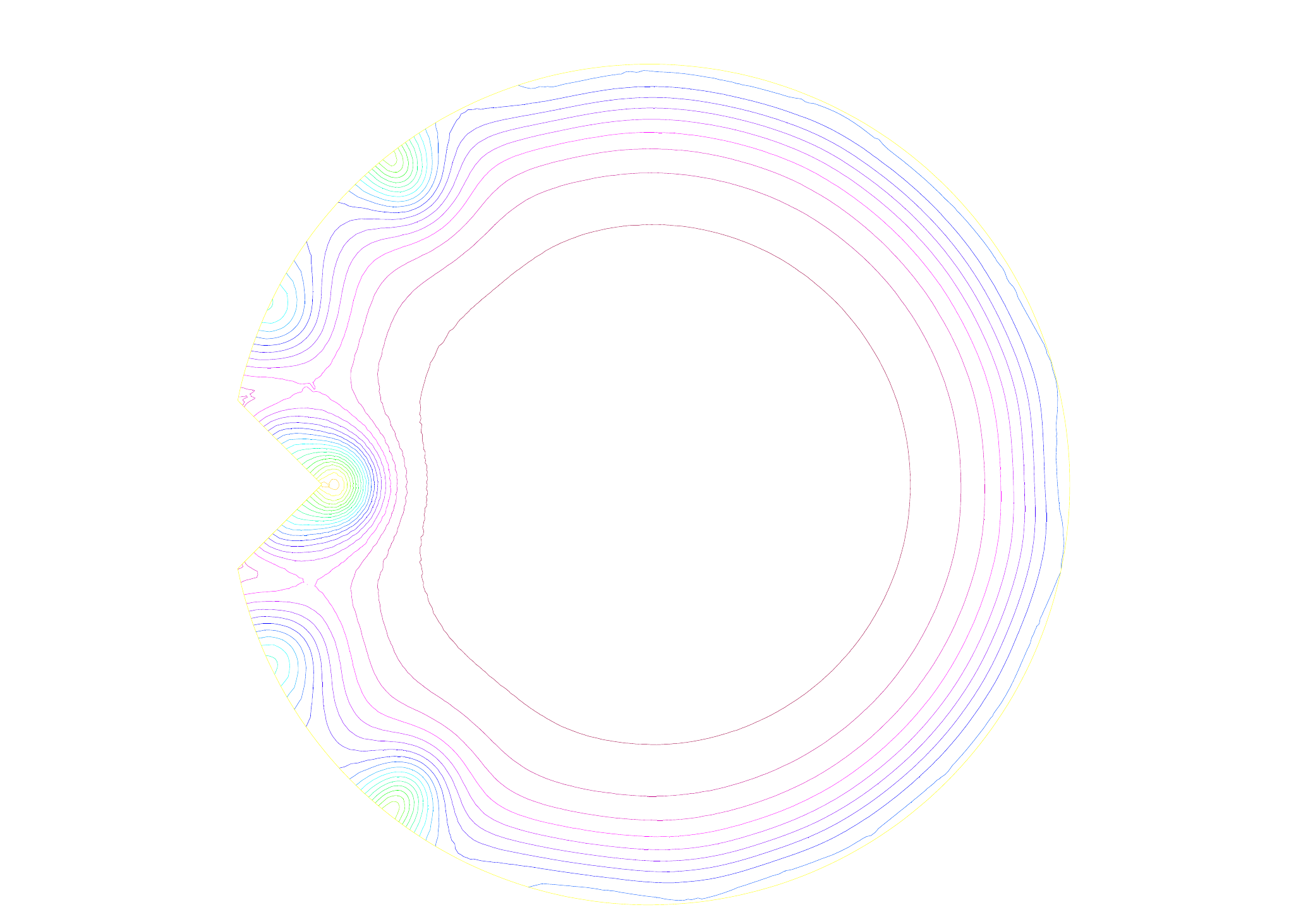}}
\subfigure[$t=30$] 
{\includegraphics[height=1.15in,width=1.6in]{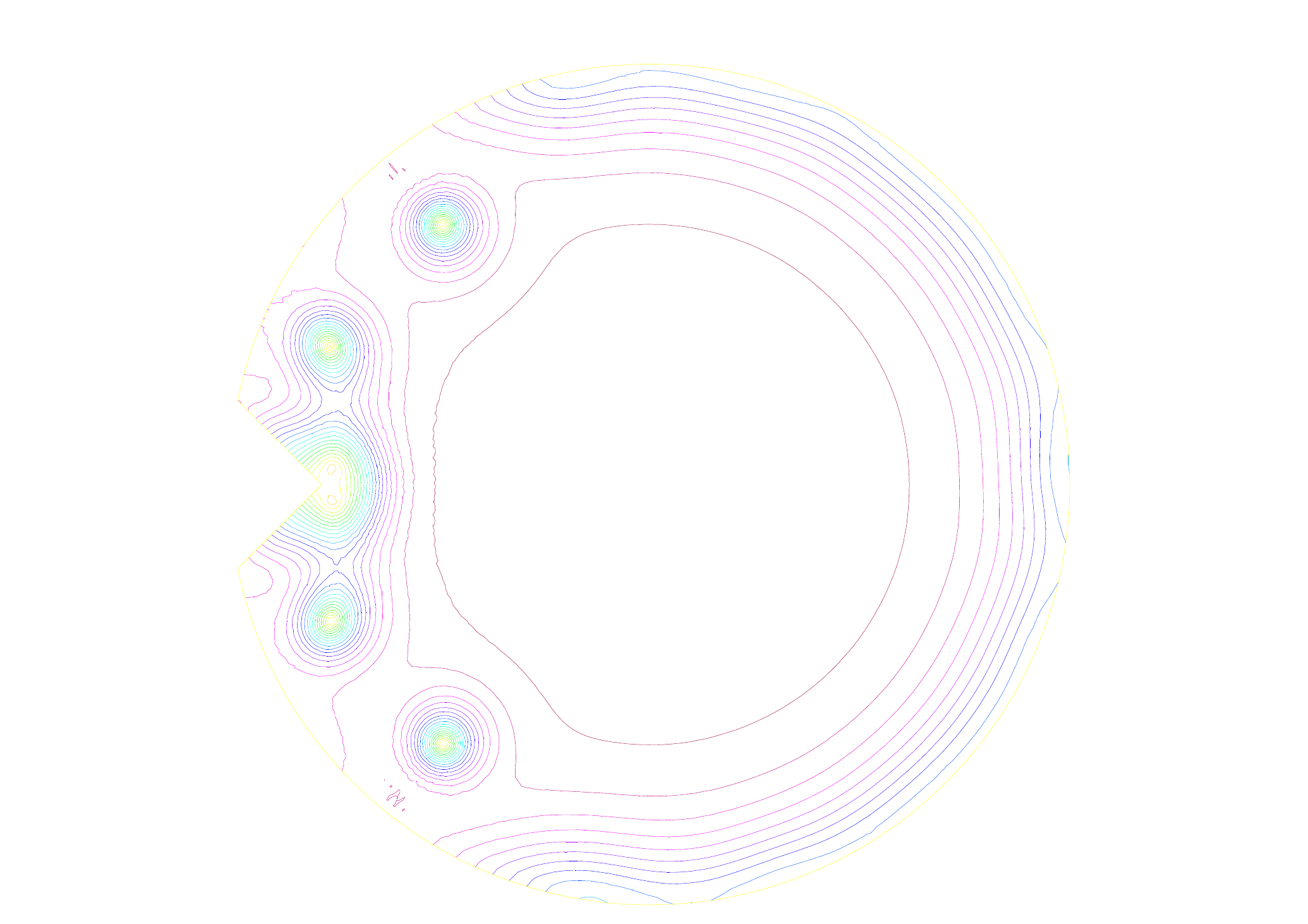}}
\subfigure[$t=5000$] 
{\includegraphics[height=1.15in,width=1.6in]{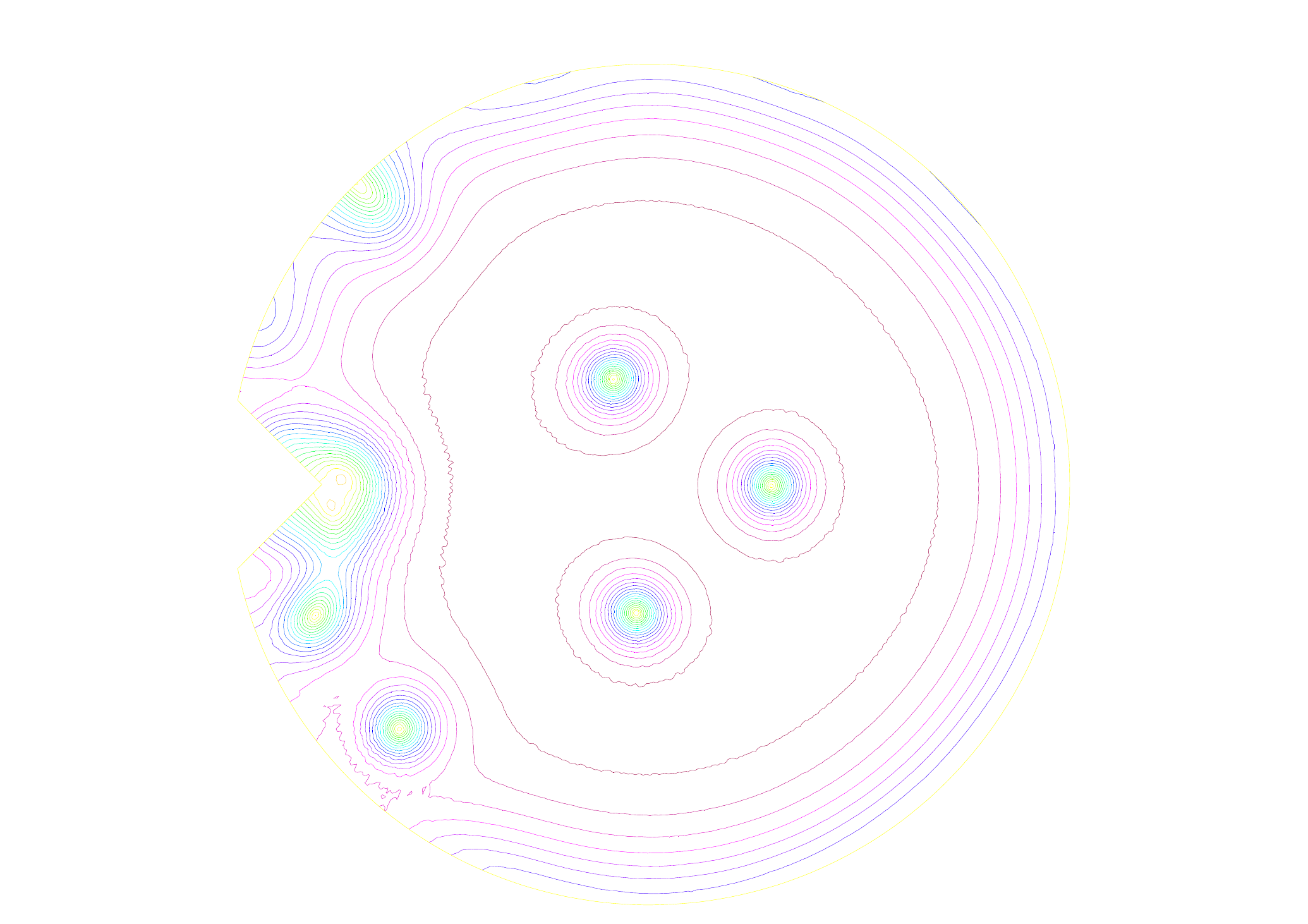}}
\caption{Contour of $|\psi|^2$ with $H=0.9$ by solving 
the TDGL under the Lorentz gauge.}
\label{FigS6}\vspace{10pt}
\subfigure[$t=25$] 
{\includegraphics[height=1.15in,width=1.6in]{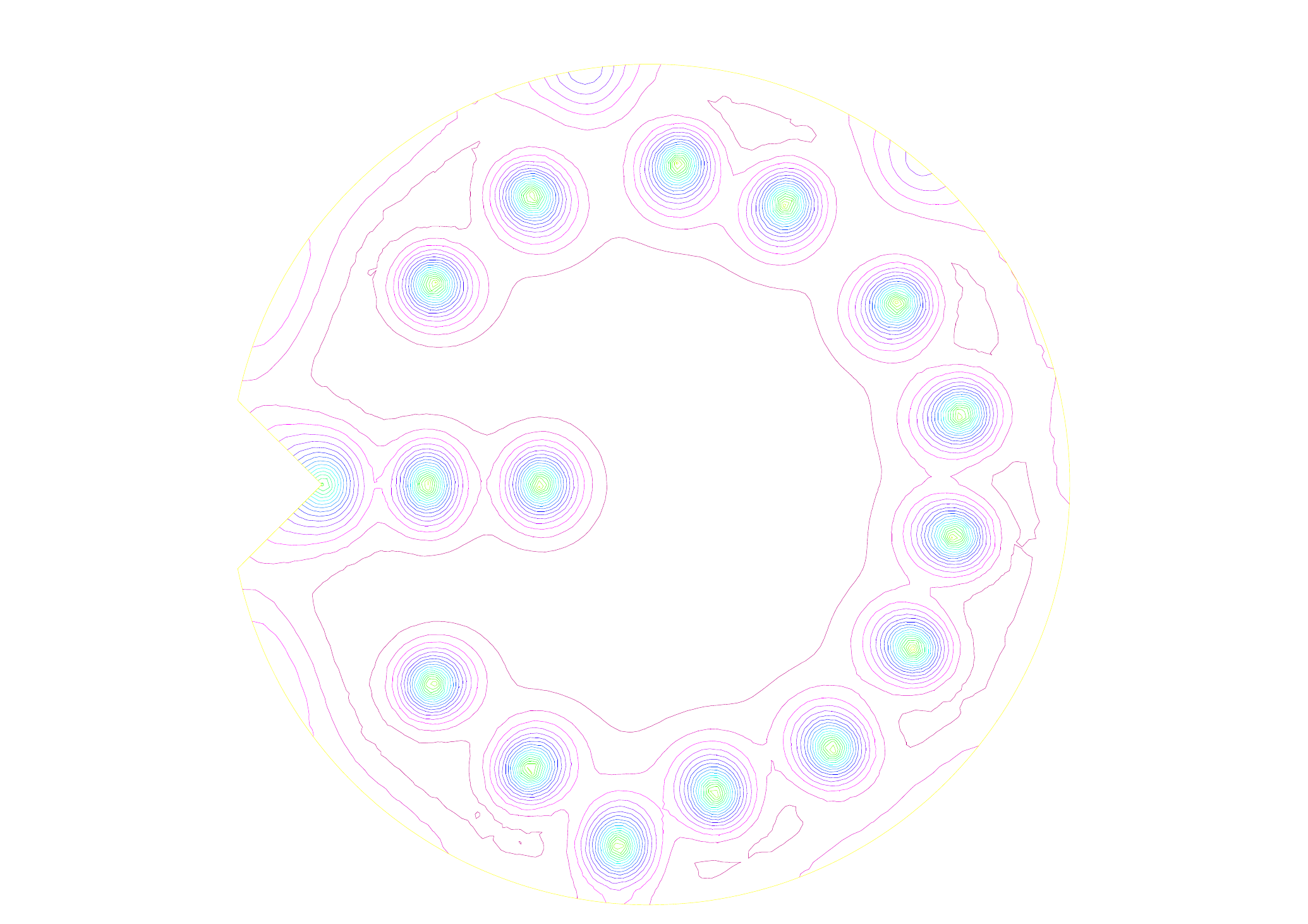}}
\subfigure[$t=30$] 
{\includegraphics[height=1.15in,width=1.6in]{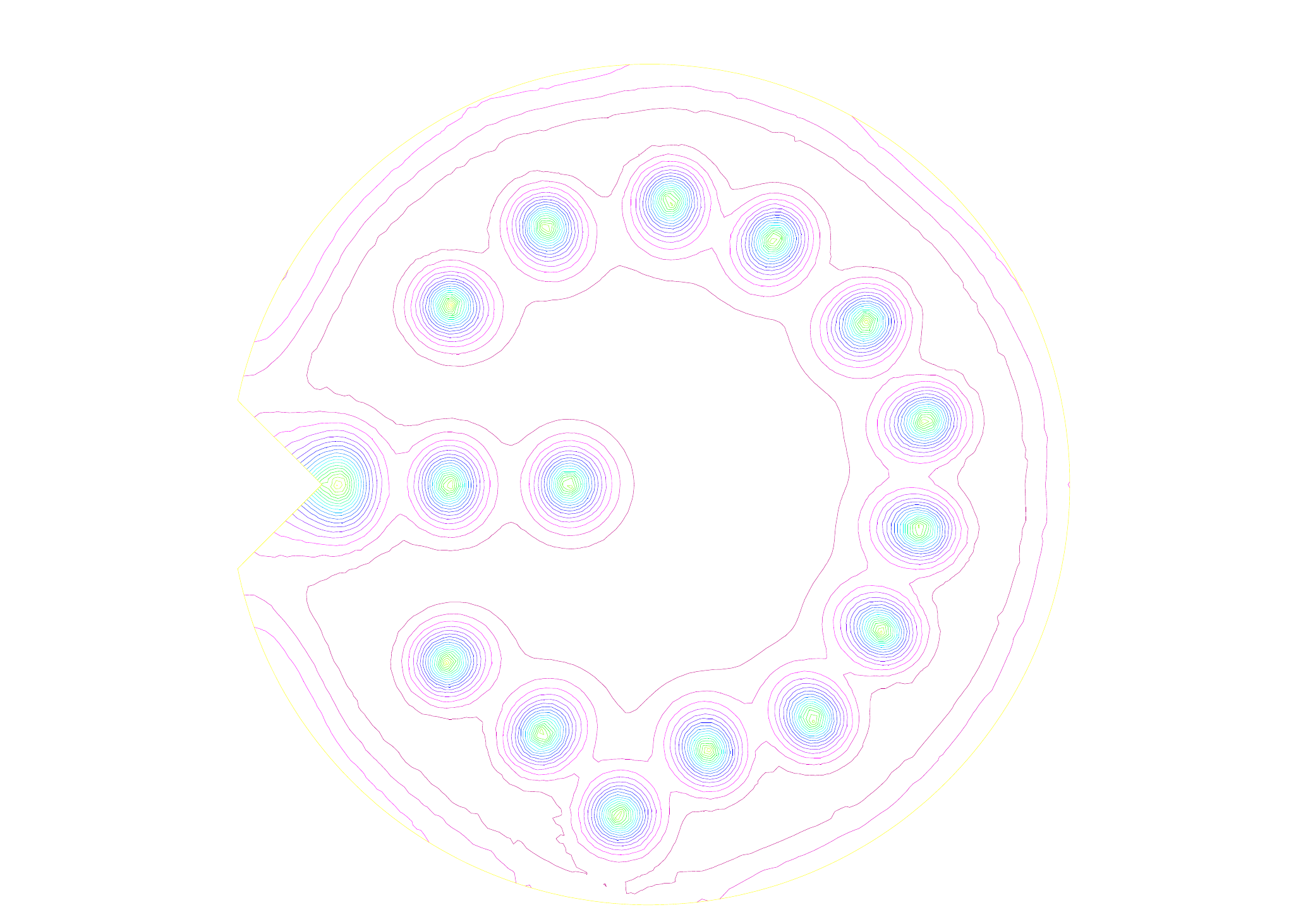}}
\subfigure[$t=5000$] 
{\includegraphics[height=1.15in,width=1.6in]{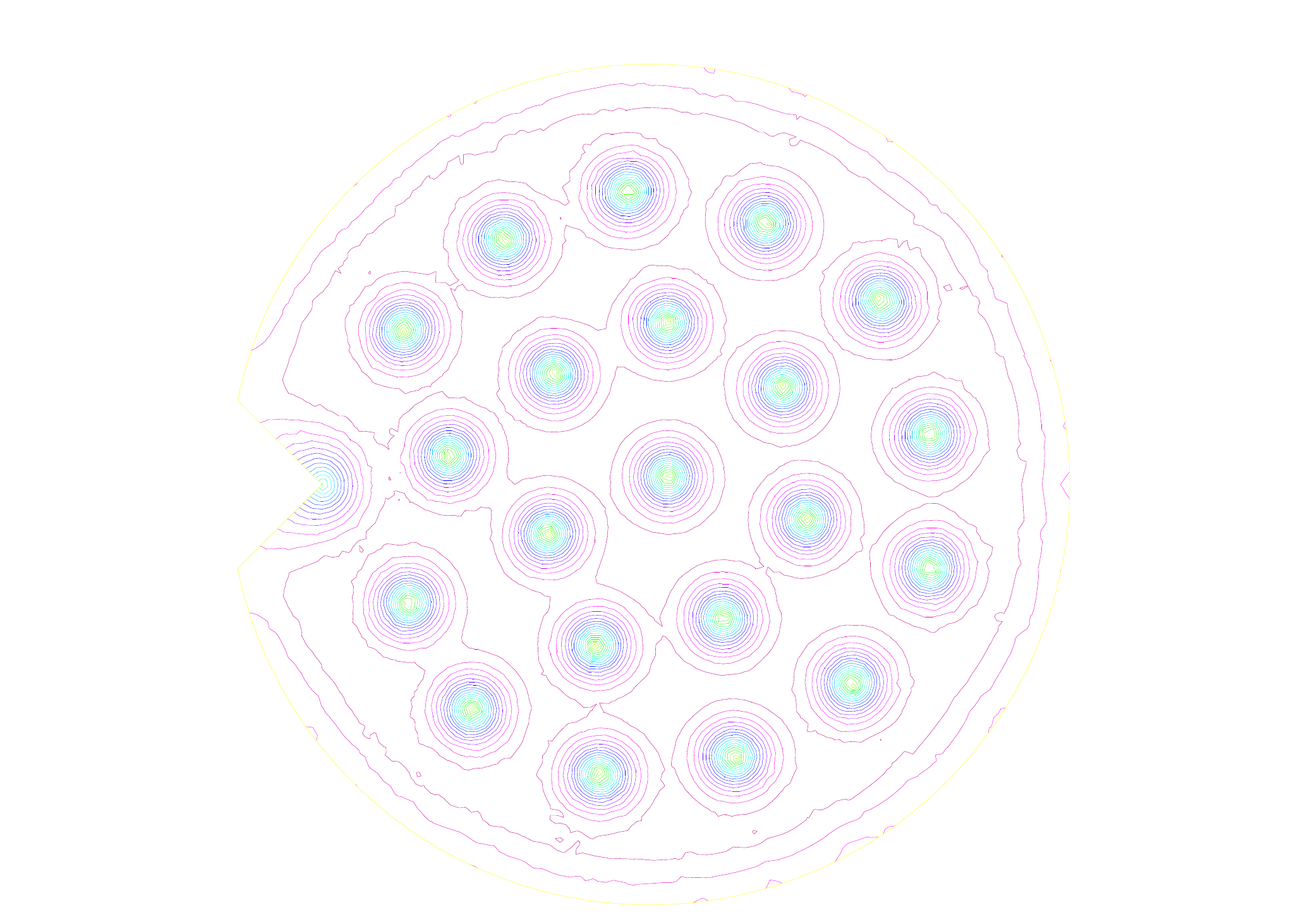}}
\caption{Contour of $|\psi|^2$ with $H=0.9$ computed by the new approach.}
\label{FigS7}\vspace{10pt}
\subfigure[$t=25$] 
{\includegraphics[height=1.15in,width=1.6in]{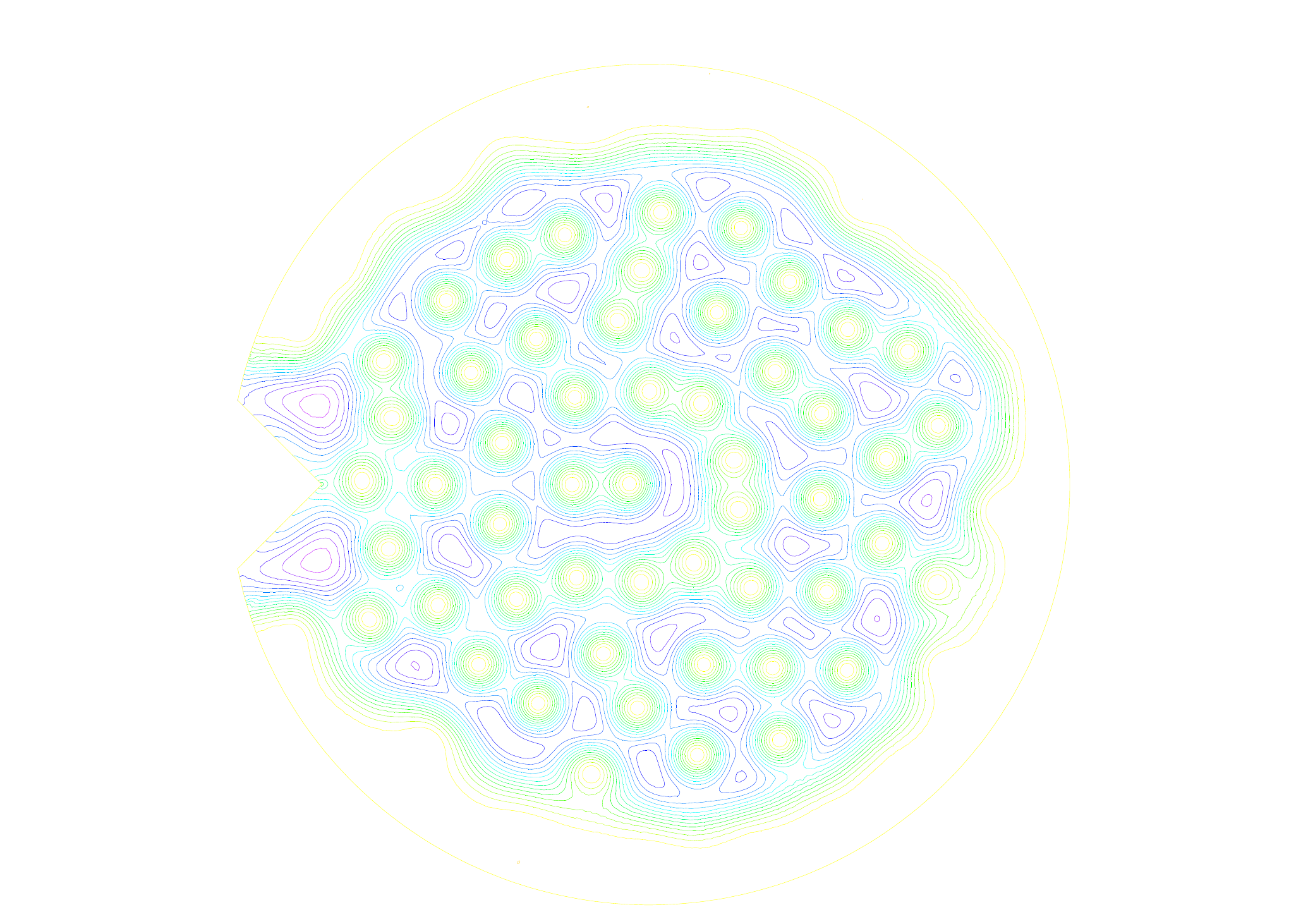}}
\subfigure[$t=50$] 
{\includegraphics[height=1.15in,width=1.6in]{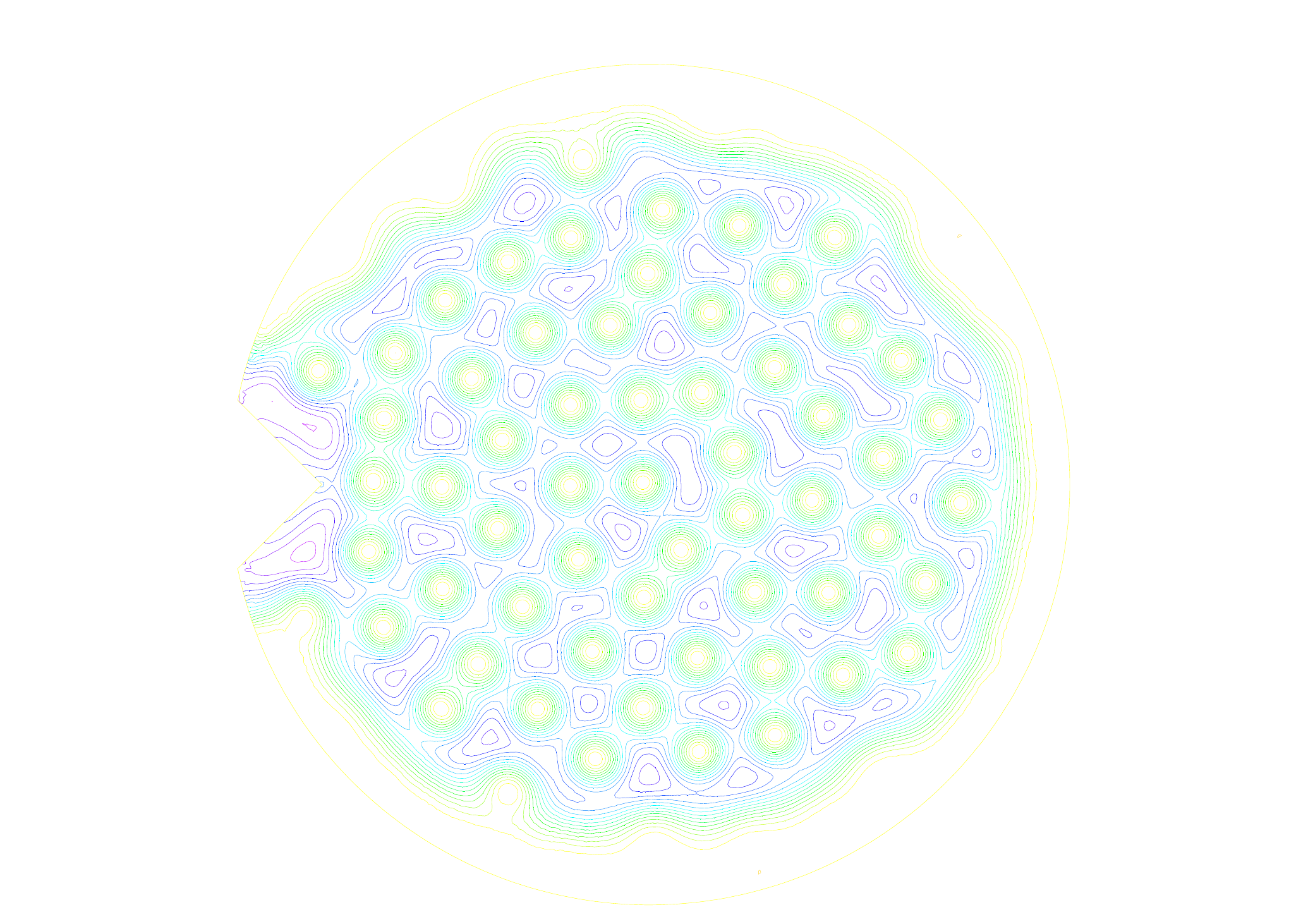}}
\subfigure[$t=100$] 
{\includegraphics[height=1.15in,width=1.6in]{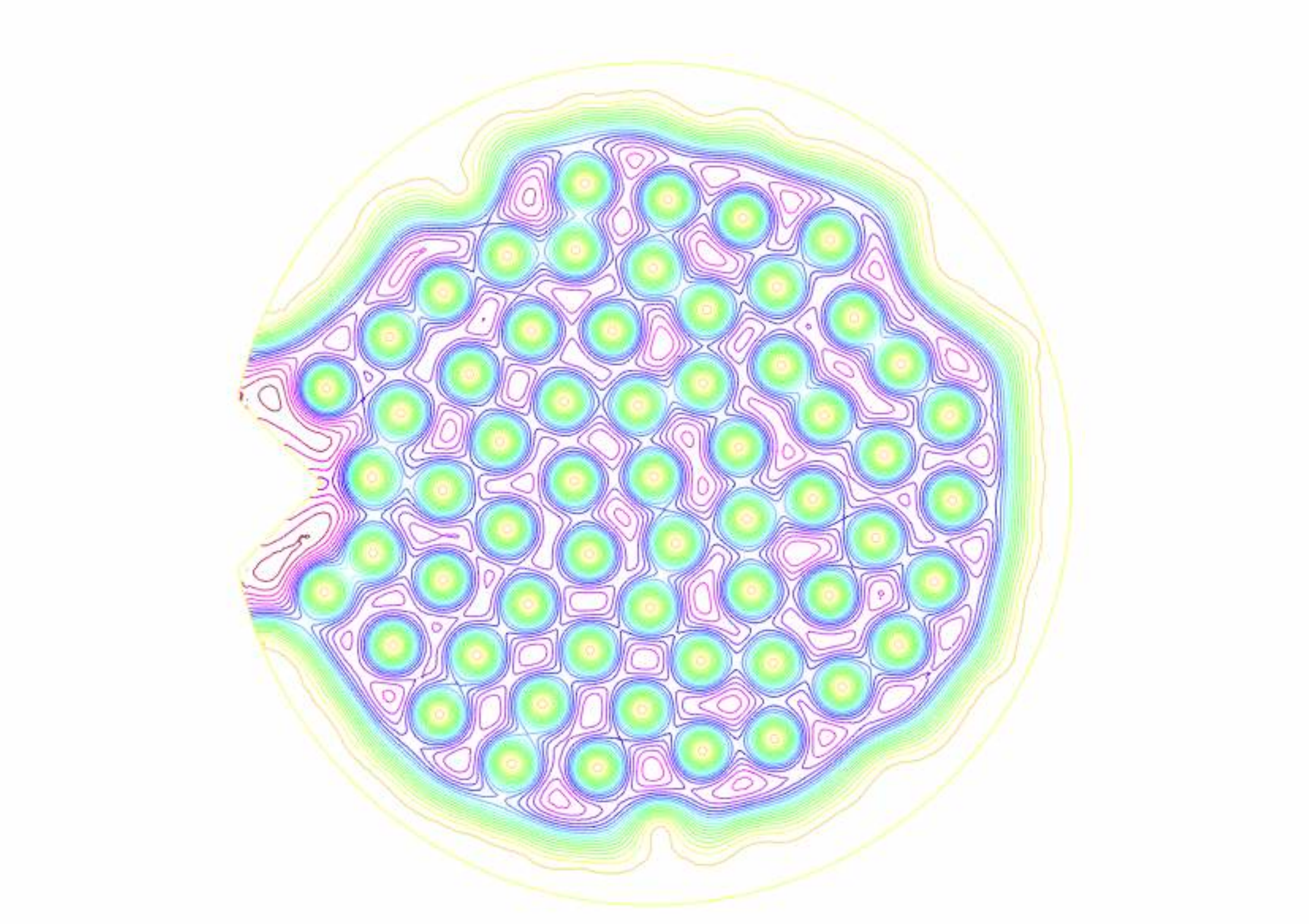}}
\caption{Contour of $|\psi|^2$ with $H=2.02$ 
by solving the TDGL under the temporal gauge.}
\label{FigS8}
\end{figure}

As the external magnetic field $H$ grows, the number of vortices 
increases and the problem becomes more difficult.
The numerical results with $H=0.9$ are present in Figure \ref{FigS5}--\ref{FigS7}.
Comparing Figure \ref{FigS5} with Figure \ref{FigS7} we see that, 
not only a wrong giant vortex may grow at the concave corner
when solving \refe{TPDEq1}-\refe{TBDC}, 
but many vortices are lost at $t=25$ and $t=30$ near the circular boundary. 

\begin{figure}[htp]
\centering
\subfigure[$t=25$] 
{\includegraphics[height=1.15in,width=1.6in]{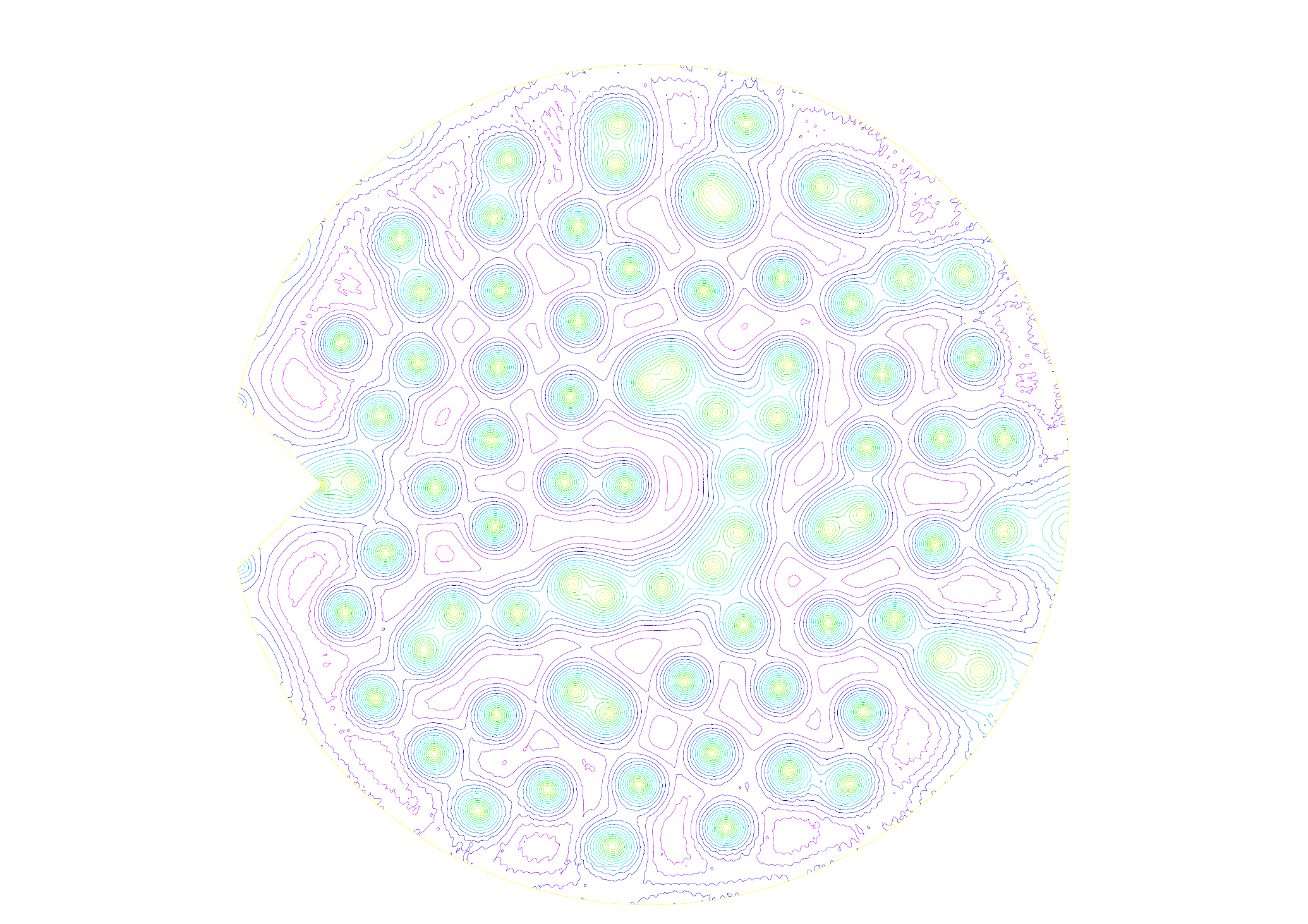}}
\subfigure[$t=50$] 
{\includegraphics[height=1.15in,width=1.6in]{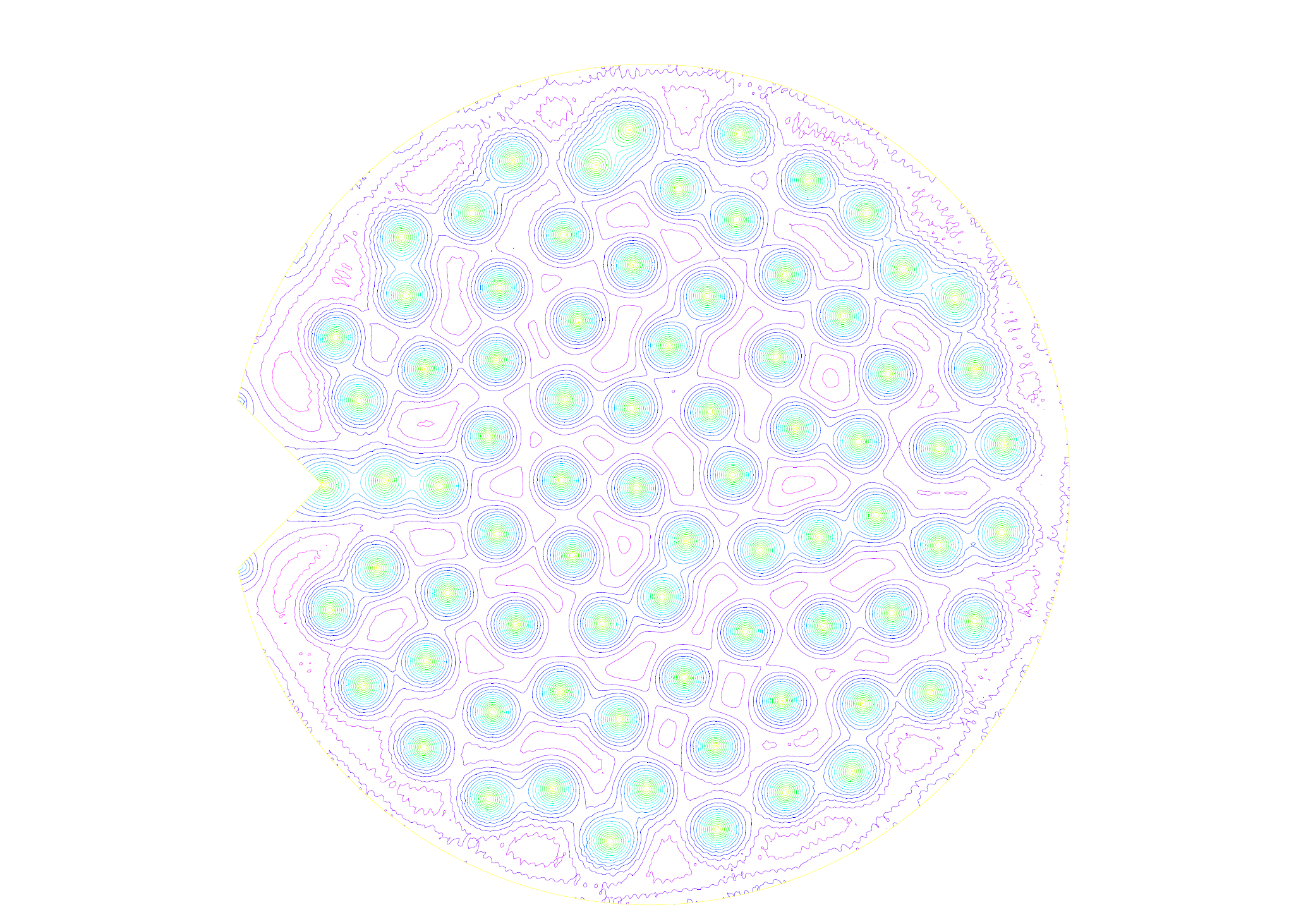}}
\subfigure[$t=100$] 
{\includegraphics[height=1.15in,width=1.6in]{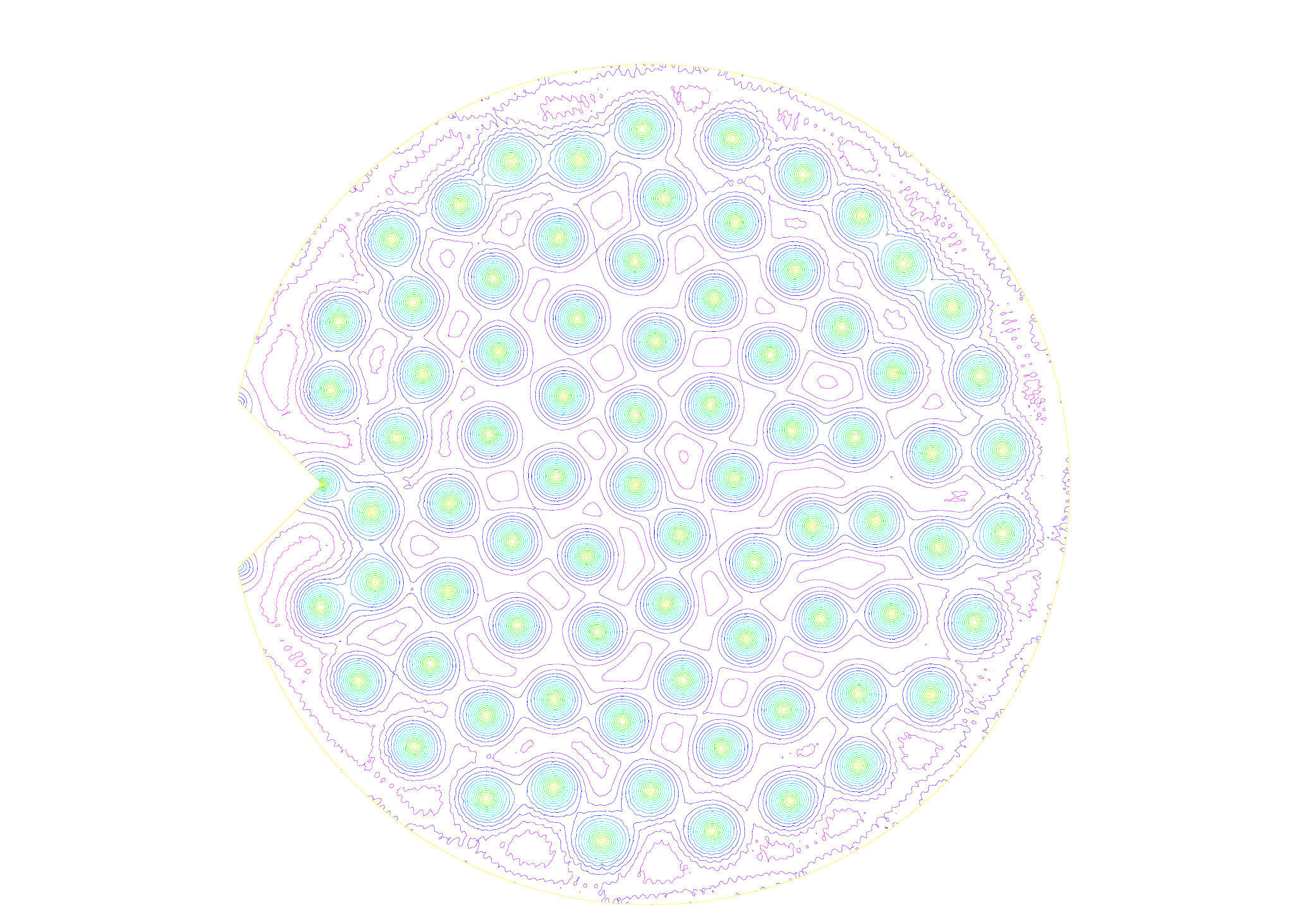}}
\caption{Contour of $|\psi|^2$ with $H=2.02$ computed by the new approach.}
\label{FigS9}\vspace{10pt}
\subfigure[$t=5$] 
{\includegraphics[height=1.15in,width=1.6in]{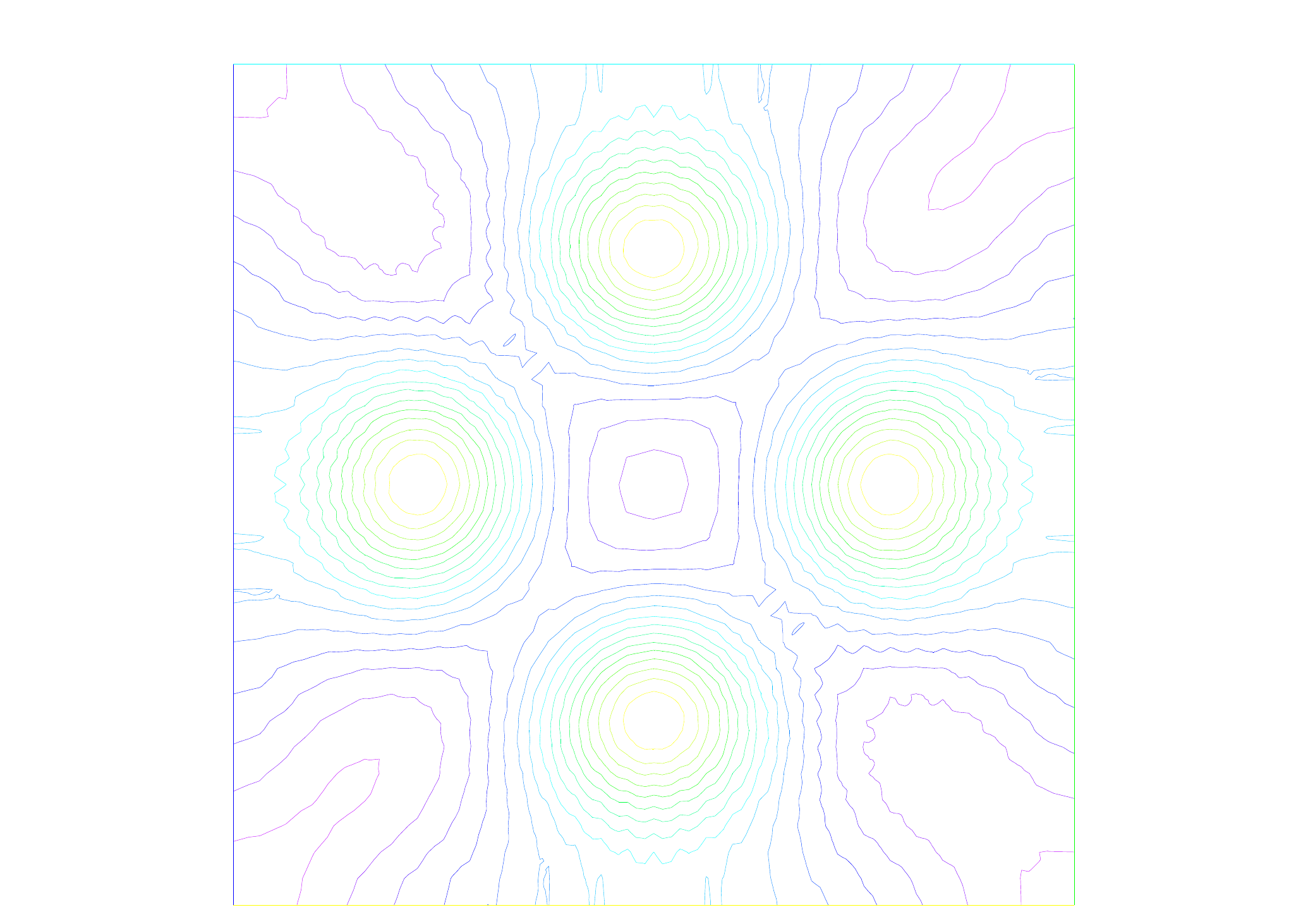}}
\subfigure[$t=20$] 
{\includegraphics[height=1.15in,width=1.6in]{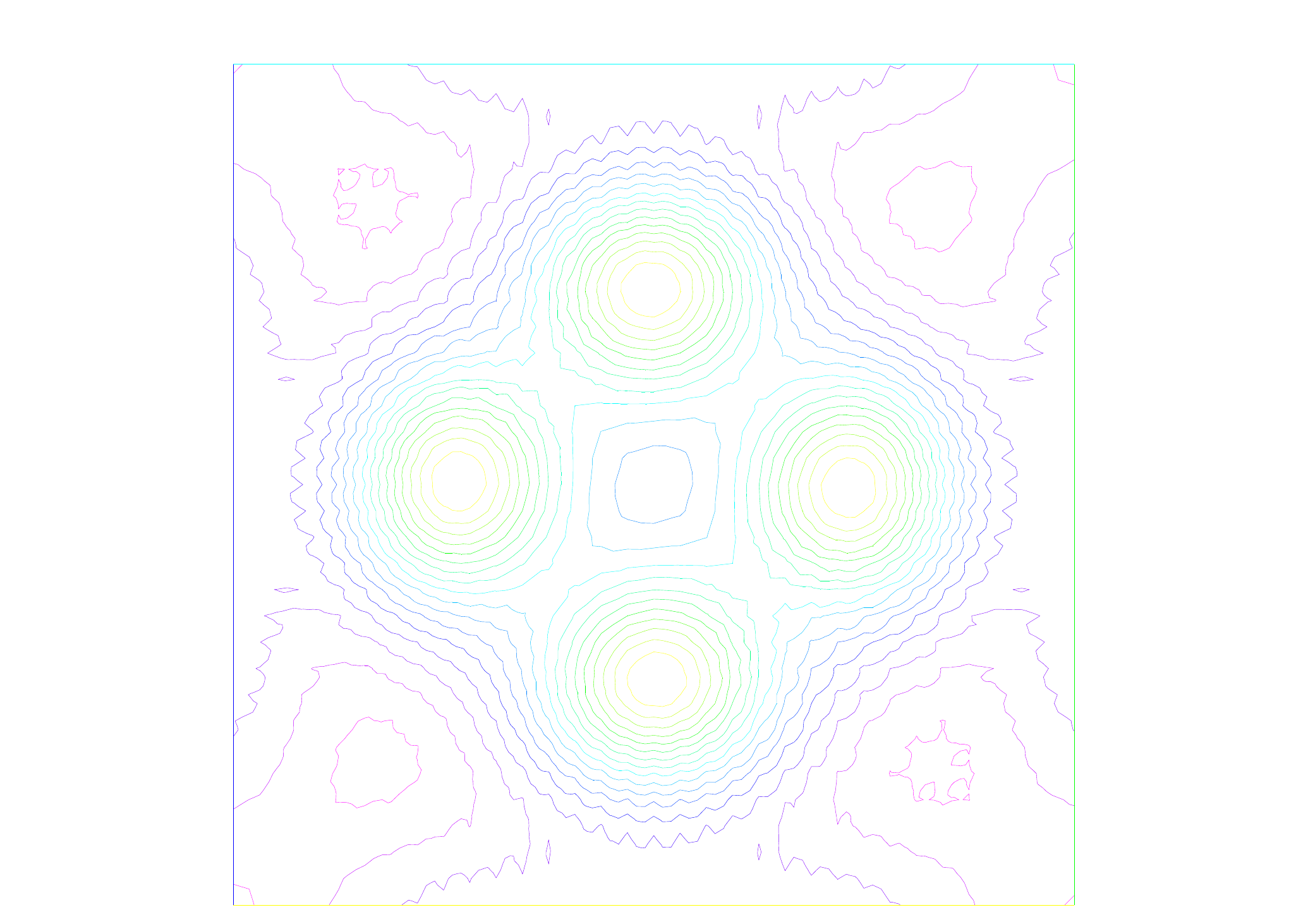}}
\subfigure[$t=40$]
{\includegraphics[height=1.15in,width=1.6in]{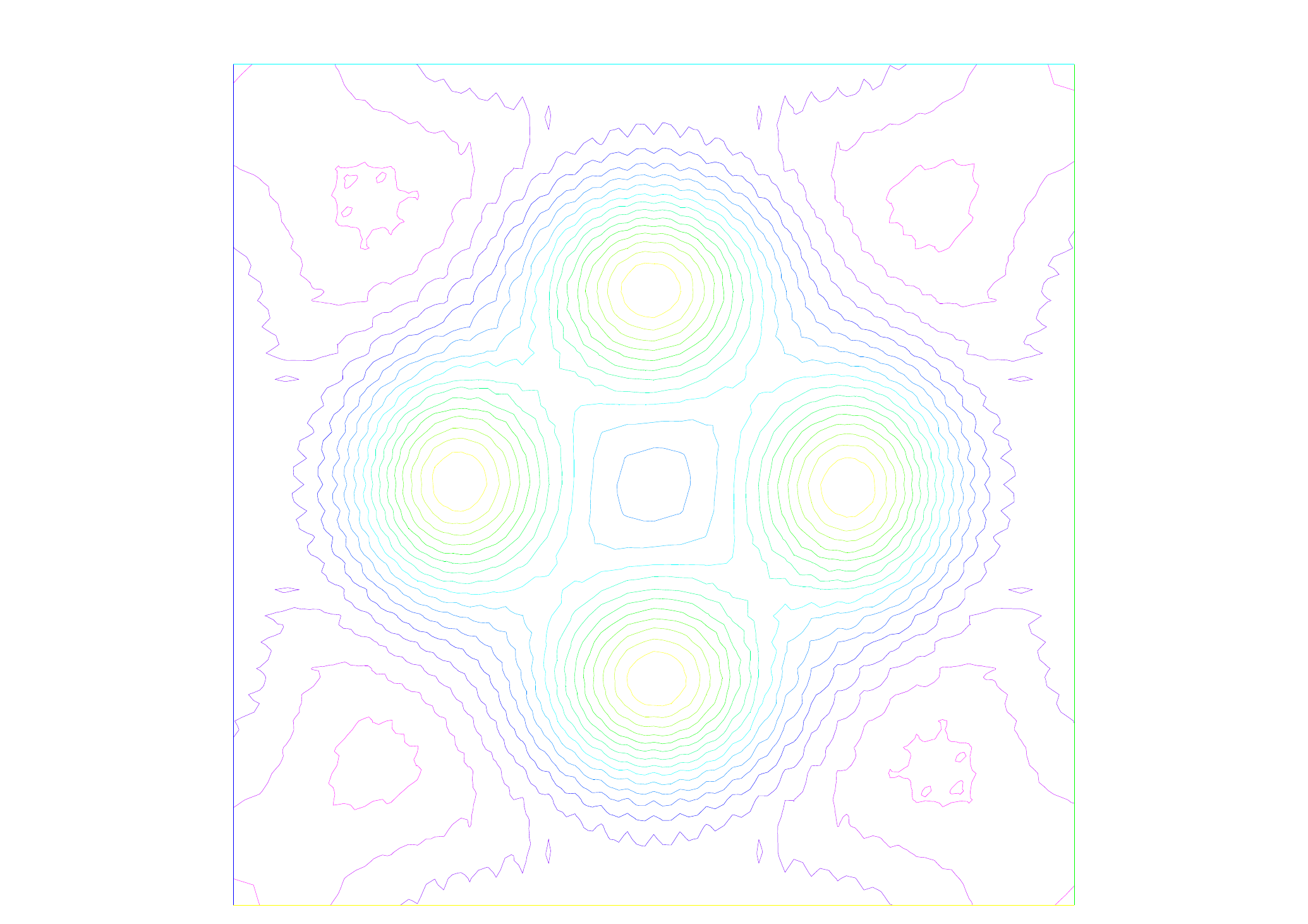}}
\caption{Contour of $|\psi|^2$ by solving 
the TDGL under the temporal gauge.}
\label{FigR1}\vspace{10pt}
\subfigure[$t=5$] 
{\includegraphics[height=1.15in,width=1.6in]{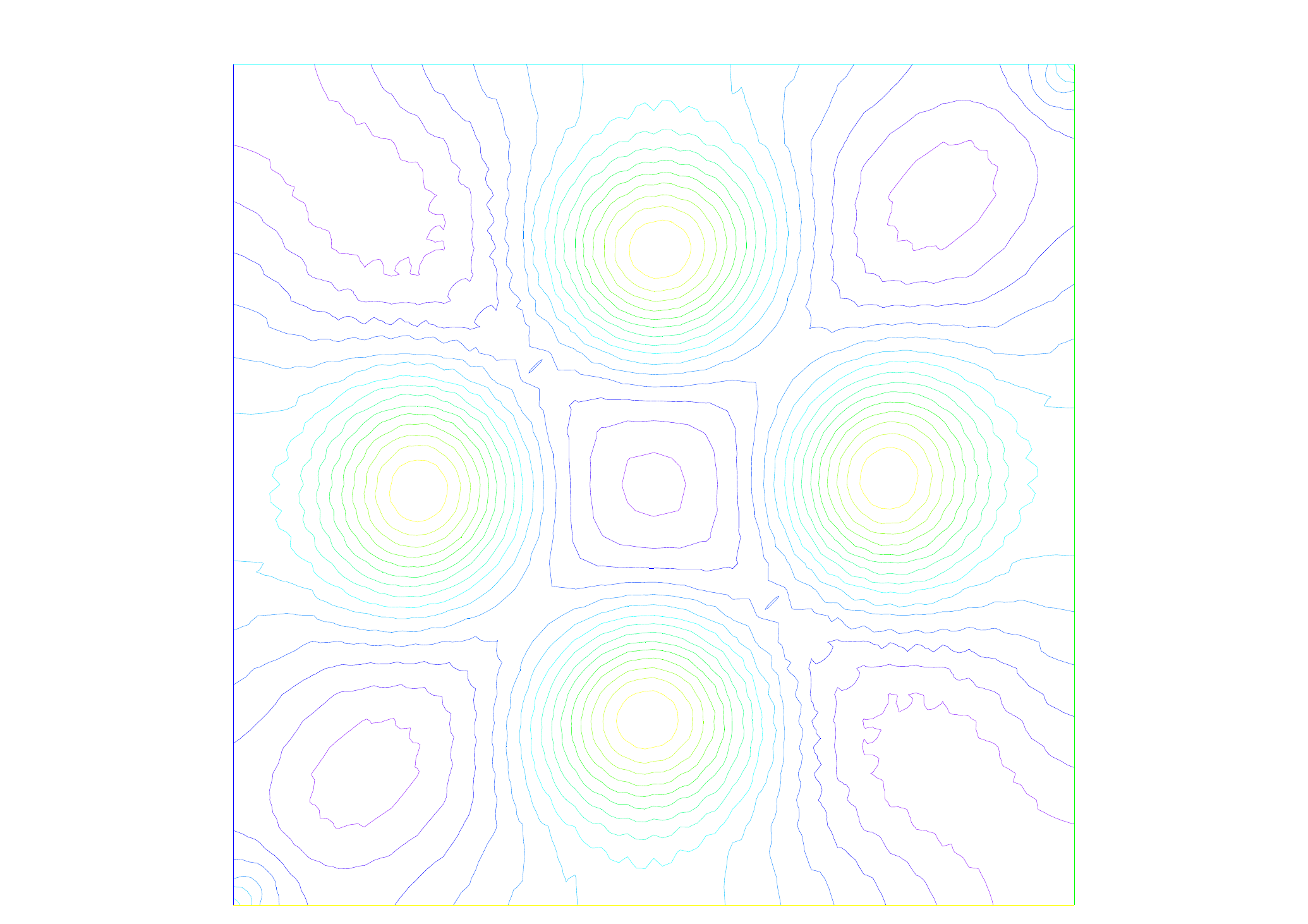}}
\subfigure[$t=20$] 
{\includegraphics[height=1.15in,width=1.6in]{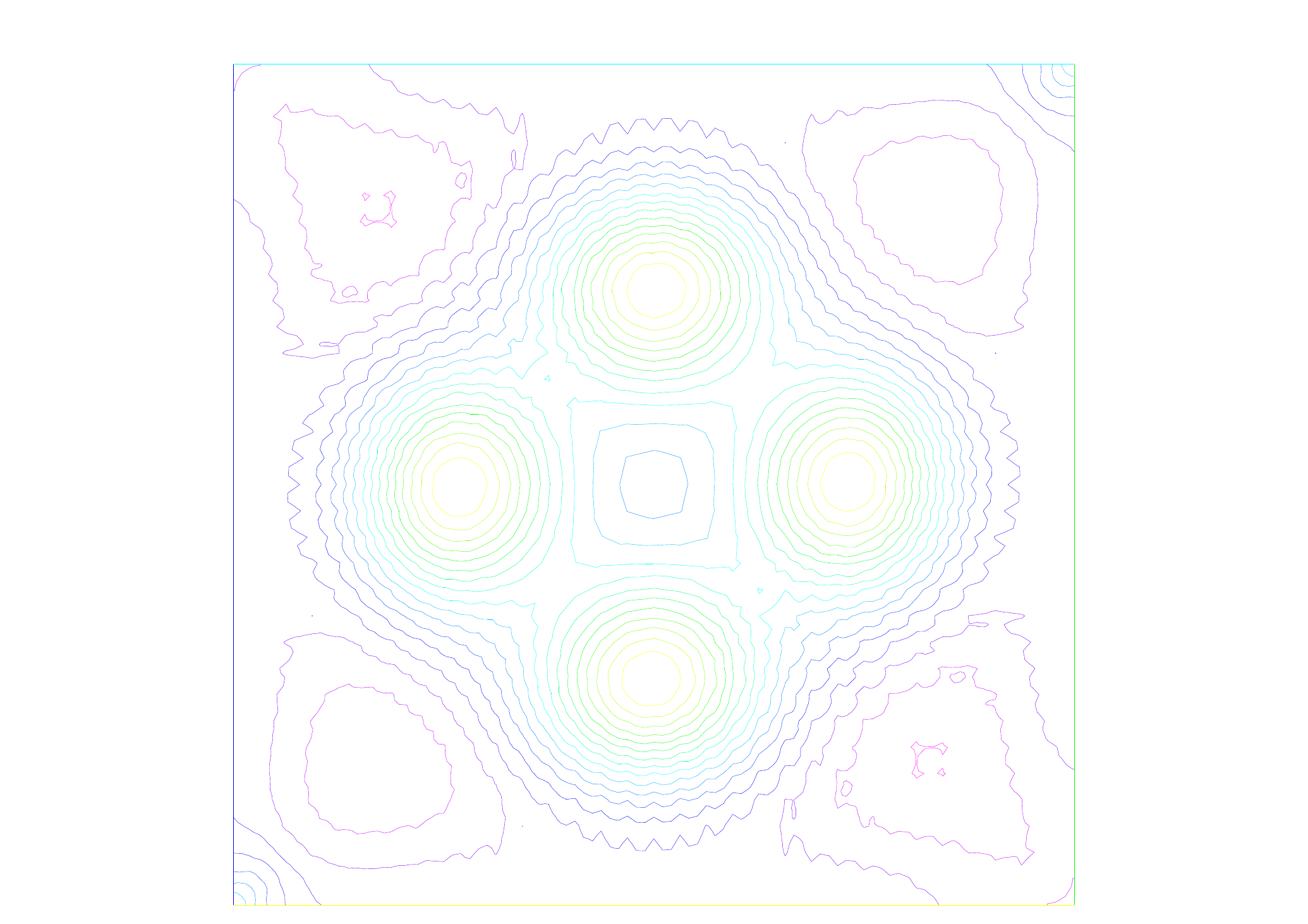}}
\subfigure[$t=40$] 
{\includegraphics[height=1.15in,width=1.6in]{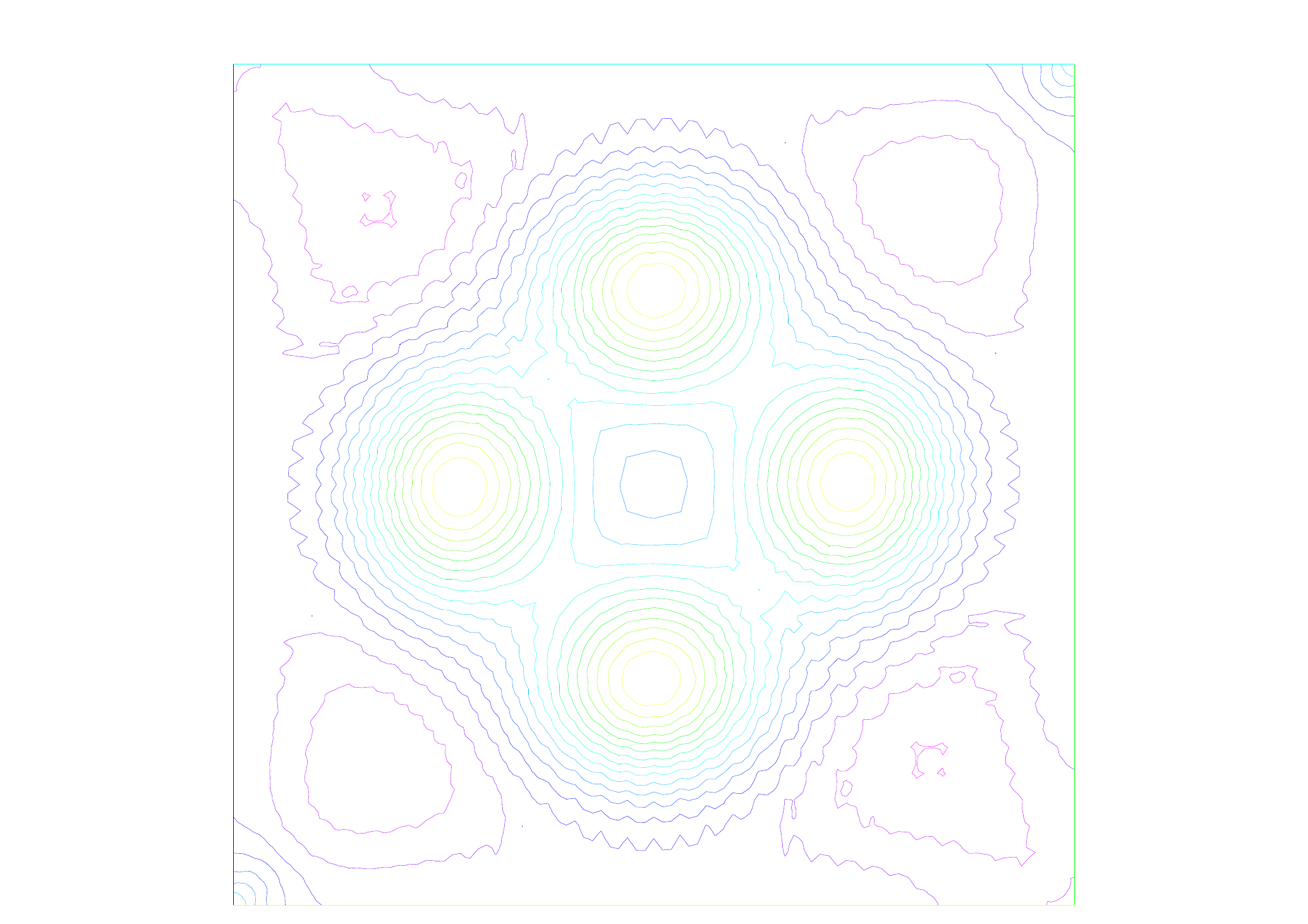}}
\caption{Contour of $|\psi|^2$ by solving 
the TDGL under the Lorentz gauge.}
\label{FigR2}\vspace{10pt}
\subfigure[$t=5$] 
{\includegraphics[height=1.15in,width=1.6in]{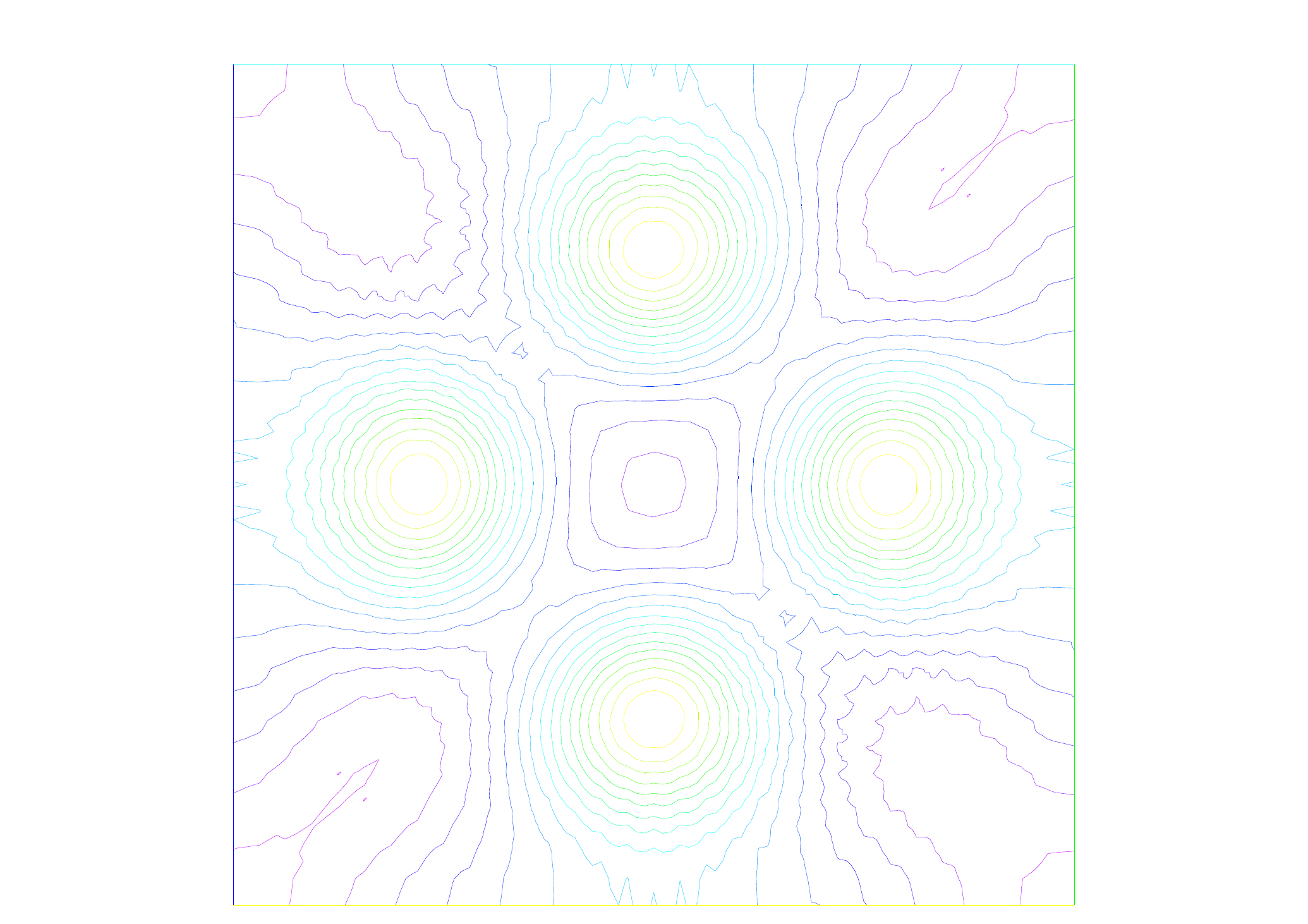}}
\subfigure[$t=20$] 
{\includegraphics[height=1.15in,width=1.6in]{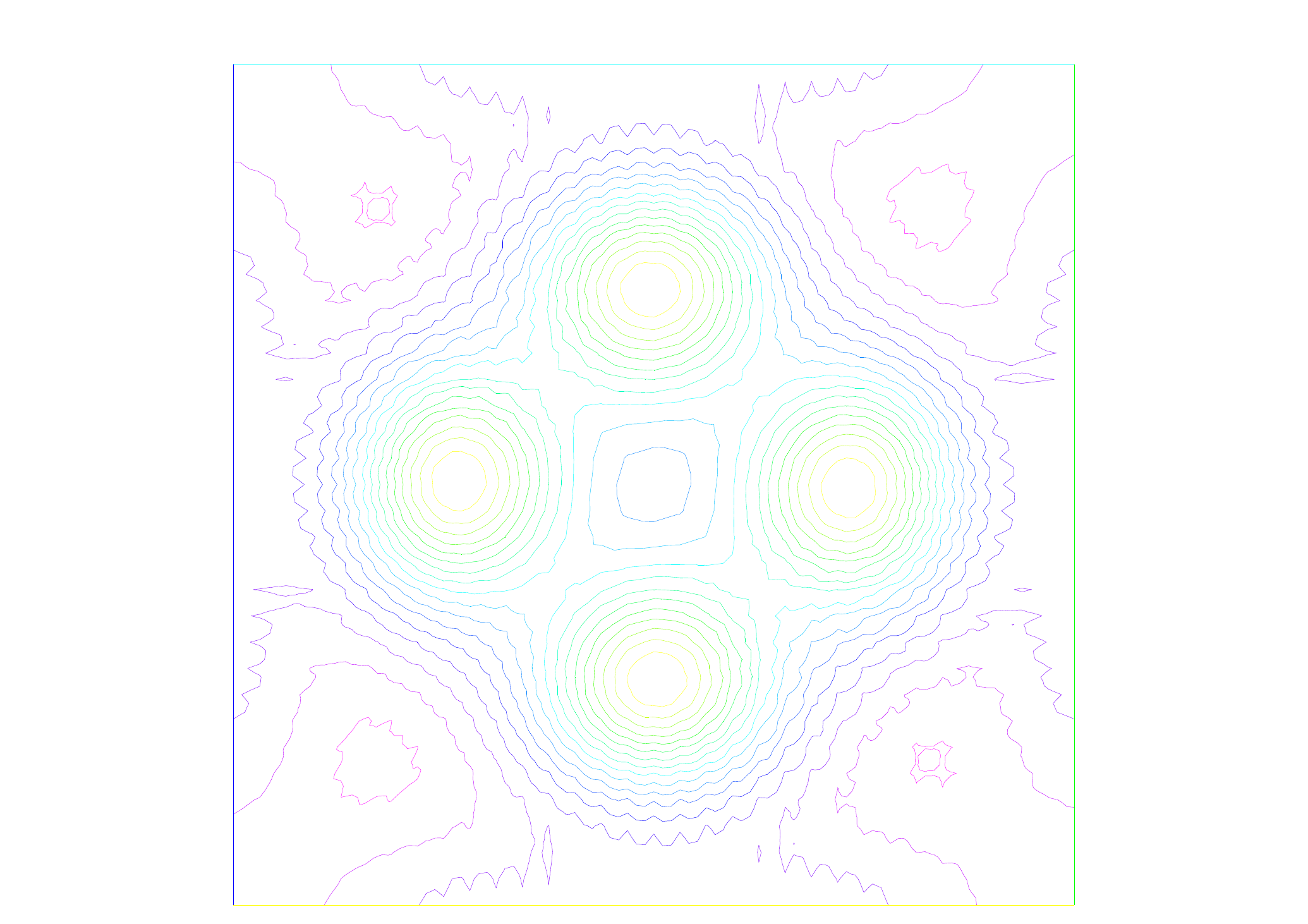}}
\subfigure[$t=40$] 
{\includegraphics[height=1.15in,width=1.6in]{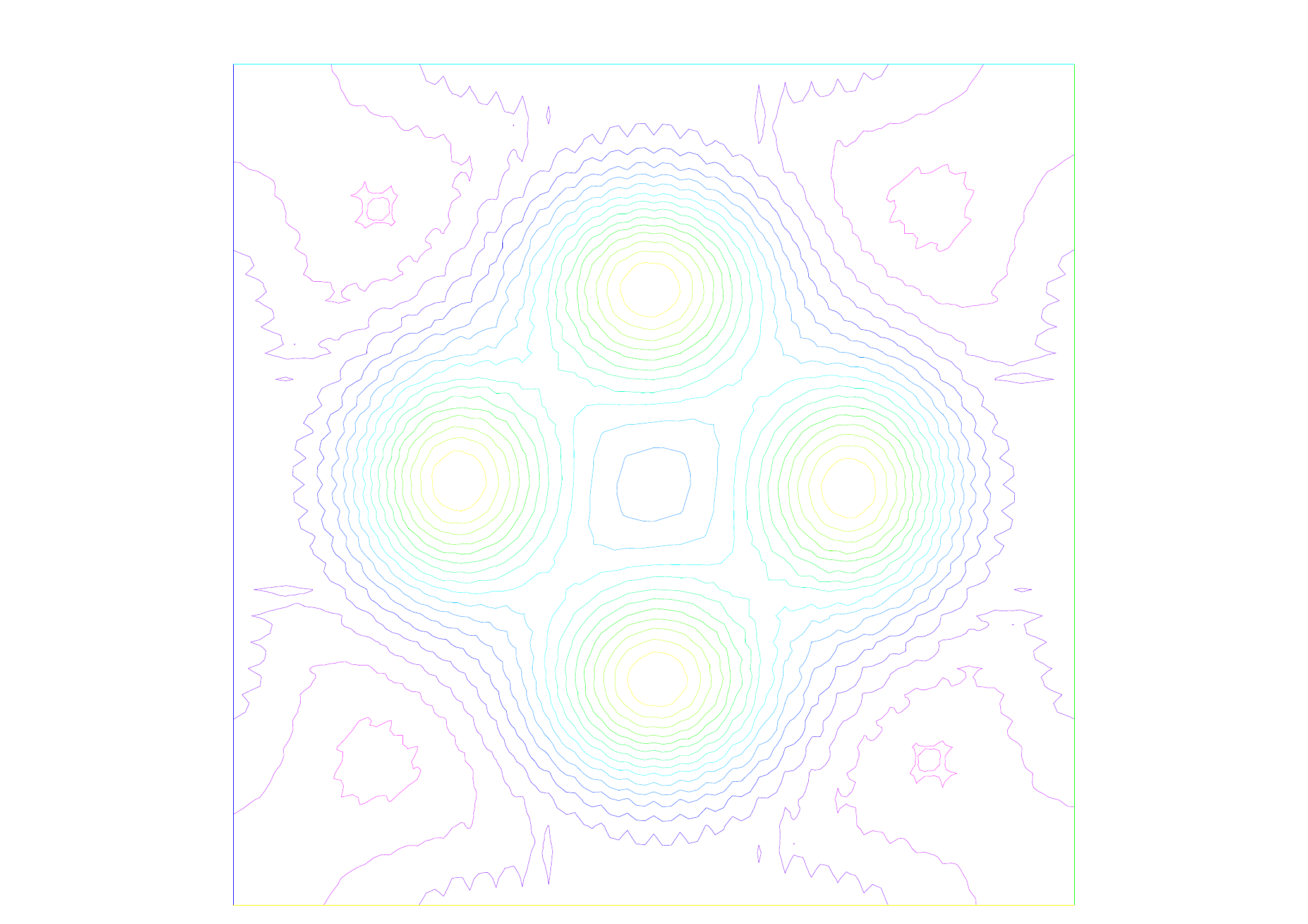}}
\caption{Contour of $|\psi|^2$ computed by the new approach.}
\label{FigR3}
\end{figure}

For $H=2.02$, there are a larger number of vortices and the 
problem becomes more difficult.
We solve \refe{TPDEq1}-\refe{TBDC} and \refe{RFPDE1}-\refe{RFinit}
with quadratic finite elements subject to a common locally refined mesh;
see Figure \ref{CircularD}-(b).
The numerical results are presented in Figure \ref{FigS8}-\ref{FigS9},
where we see that the numerical solution of 
the TDGL under the temporal gauge 
looks strange, while the numerical solution given by 
our new approach looks reasonable. 
This shows that our new approach is also superior 
with locally refined mesh and high-order finite elements.\bigskip

{\it Example 3.3}$\quad$
Finally, we solve the three systems in a convex domain 
$\Omega = (0,1) \times (0,1)$ with $\kappa=10$ and 
\begin{align*}
&\psi_0=0.6+0.8\, i,\qquad {\bf A}_0=(0,0),\qquad H=5  .
\end{align*}
This example was tested before in \cite{CD01} 
by solving \refe{PDE1}-\refe{TBDC2}, and tested in \cite{Yang}
by solving \refe{TPDEq1}-\refe{TBDC}.
We triangulate the domain $\Omega$ into uniform right
triangles with 32 points on each side. 
The contour plots of $|\psi|^2$ at different time levels are 
presented in Figure \ref{FigR1}-\ref{FigR3} by solving the
equations with the time-step size $\tau=0.1$,
which show that solving the three systems gives almost the same solution
in the domain alway from the corners.  
Although there is a little difference between Figure \ref{FigR2}
and Figure \ref{FigR1}-\ref{FigR3} near the corners,
this difference can be eliminated by using a smaller mesh size.
Roughly speaking, the three systems are equivalent 
in a domain without reentrant corners,
both theoretically and computationally. \bigskip

\section{Conclusions}
$\quad$
We have introduced a new approach
for the numerical simulation of the time-dependent
Ginzburg--Landau model of superconductivity
in a general curved polygon which may contain
reentrant corners,
by reformulating the equations under
the Lorentz gauge into an equivalent system
of equations.
Mathematically speaking, this new approach is more suitable for
Ginzburg--Landau equations with strong corner singularities. Indeed, numerical
simulations demonstrate that it is more stable
and accurate than the traditional approaches
in the presence of a reentrant corner,
and comparably accurate as the
traditional approaches  in a convex domain.
\bigskip

\noindent{\bf Acknowledgement.}$\quad$
We would like to thank Professor Qiang Du
for helpful discussions.\bigskip

\newpage
$\,$

\end{document}